\documentclass{amsart}

\usepackage{lipsum}
\usepackage{amsfonts}
\usepackage{graphicx}
\usepackage{epstopdf}

\usepackage{algorithm}

\usepackage{amsmath,amscd,amssymb}
\usepackage{xfrac}
\usepackage{url}
\usepackage{mathtools}
\usepackage{color}
\usepackage{subcaption}
\usepackage[pagebackref,colorlinks,citecolor=blue,linkcolor=magenta]{hyperref}

\usepackage[capitalize]{cleveref}
\usepackage{tikz}
\usetikzlibrary{calc}
\usepackage{color}
\usepackage{centernot}
\usepackage{enumitem}
\usepackage[utf8]{inputenc}
\usepackage{algpseudocode}

\newcommand{\GG}{\mathcal{G}}
\newcommand{\HH}{\mathcal{H}}
\newcommand{\DD}{\mathcal{D}}
\newcommand{\SSS}{\mathcal{S}}
\newcommand{\RR}{\mathbb{R}}
\newcommand{\spans}{\operatorname{span}}

\renewcommand{\ne}{\operatorname{ne}}
\newcommand{\child}{\mathfrak{c}}
\newcommand{\parent}{\mathfrak{p}}

\makeatletter
\DeclareMathAlphabet{\@mymathbb}{U}{bbold}{m}{n}
\newcommand{\zero}{\@mymathbb{0}}
\newcommand{\one}{\@mymathbb{1}}
\makeatother

\DeclareMathOperator{\pa}{pa}
\DeclareMathOperator{\ch}{ch}
\DeclareMathOperator{\de}{de}
\DeclareMathOperator{\an}{an}
\DeclareMathOperator{\nd}{nd}
\DeclareMathOperator{\cl}{cl}

\DeclareMathOperator{\BIC}{BIC}

\DeclareMathOperator{\argmax}{argmax}

\newcommand\independent{\protect\mathpalette{\protect\independenT}{\perp}}
\def\independenT#1#2{\mathrel{\rlap{$#1#2$}\mkern2mu{#1#2}}}

\def\newop#1{\expandafter\def\csname #1\endcsname{\mathop{\rm
#1}\nolimits}}
\newop{STAB}
\newop{conv}
\newop{CIM}
\newop{odd}
\newop{even}

\usepackage{enumitem}
\setlist[enumerate]{leftmargin=.5in}
\setlist[itemize]{leftmargin=.5in}

\newtheorem{theorem}{Theorem}[section]
\newtheorem{proposition}[theorem]{Proposition}

\newtheorem{lemma}[theorem]{Lemma}

\theoremstyle{definition}
\newtheorem{definition}[theorem]{Definition}

\theoremstyle{remark}

\theoremstyle{plain}
\newtheorem{example}[theorem]{Example} 

\title{On the Edges of Characteristic Imset Polytopes}

\author{Svante Linusson% \thanks{Matematik, KTH, SE-100 44 Stockholm, Sweden
\and Petter Restadh%\thanks{Matematik, KTH, SE-100 44 Stockholm, Sweden
\and Liam Solus%\thanks{Matematik, KTH, SE-100 44 Stockholm, Sweden
}

\email[Svante Linusson]{linusson@math.kth.se}
\email[Petter Restadh]{petterre@kth.se}
\email[Liam Solus]{solus@kth.se}

\address{Department of Mathematics\\
    KTH Royal Institute of Technology\\
    SE-100 44 Stockholm, Sweden}

\begin{document}
\maketitle

%%%%%%%%%%%%%%%%%%%%%%%%%%
%---ABSTRACT
\begin{abstract}
The edges of the characteristic imset polytope, $\CIM_p$, were recently shown 
to have strong connections to causal discovery as many algorithms could be 
interpreted as greedy restricted edge-walks, even though only a strict subset 
of the edges are known. 
To better understand the general edge structure of the polytope we describe
the edge structure of faces with a clear combinatorial interpretation: for any undirected graph $G$ we have the face $\CIM_G$, the convex hull of the characteristic imsets of DAGs with skeleton $G$. 
We give a full edge-description of $\CIM_G$ when $G$ is a tree,
leading to interesting connections to other polytopes. 
In particular the well-studied stable set polytope can be recovered as a face 
of $\CIM_G$ when $G$ is a tree. Building on this connection we are also able to give a description of all 
edges of $\CIM_G$ when $G$ is a cycle, 
suggesting possible inroads for generalization.
We then introduce an algorithm for learning directed trees from data, 
utilizing our newly discovered edges, that outperforms classical methods on 
simulated Gaussian data. 
\end{abstract}

%%%%%%%%%%%%%%%%%%%%%%%%%%
%---SECTION: Intro
\section{Introduction}
\label{sec: intro}
Let $[p]\coloneqq\{1,\dots,p\}$. 
For a directed acyclic graph (DAG) $\GG=([p], E)$ Studen\'y, Hemmecke and Lindner introduced the \emph{characteristic imset} \cite{SHL10}  defined as the function
$c_\GG\colon \{S\subseteq [p]\colon |S|\geq 2\}\to \{0,1\}$ where
\[
c_\GG(S) \coloneqq
\begin{cases}
1	&	\text{ if there exists $i\in S$ such that $S\setminus\{i\}\subseteq \pa_\GG(i)$},	\\
0	&	\text{ otherwise}.	\\
\end{cases}
\]
Since $c_\GG$ is a function from a finite set into $\mathbb{R}$ we identify it with a vector of length $2^p-p-1$, and thus we can consider the \emph{characteristic imset polytope}:
\[
\CIM_p\coloneqq \conv\left(c_\GG\colon \GG=([p], E)\textrm{ a DAG}\right).
\] 
The polytope $\CIM_p$ is a full-dimensional 0/1-polytope 
% , the vertices have only 0 or 1 as entries, 
whose vertices are exactly the characteristic imsets of DAGs \cite{S05, SHL10}.  
A \emph{v-structure} is an induced subgraph of the form $i\to j\leftarrow k$ and the \emph{skeleton} of a directed graph $\GG$ is the undirected graph that has the same adjacencies as $\GG$.
A unique DAG $\GG$ cannot in general be recovered from the characteristic imset $c_\GG$.
We can, however, easily recover both the skeleton and the v-structures. 
Vice versa, any two graphs that share the same skeleton and v-structures have the same characteristic imset \cite{L12, SHL10} (see \Cref{lem: imset structure} and \ref{lem: lindner}). 
Two DAGs that have the same skeleton and the same v-structures are known as \emph{Markov equivalent} (\Cref{thm: verma pearl}), and hence belong to the same \emph{Markov equivalence class (MEC)}.
That is, each characteristic imset corresponds to a unique MEC. 
A graphical representation of a MEC, called the \emph{essential graph}, was given by Andersson et. al. \cite{AMP97}. % as a partially directed graph. 
The essential graph of a DAG $\GG$ is a partially directed graph that has the same skeleton as $\GG$ and directed edges being exactly those edges that have the same direction in all DAGs in the MEC of $\GG$.

The polyhedral geometry of $\CIM_p$ has been previously studied.
Cussens, Haws and Studen{\`y} were able to obtain classes of facets of $\CIM_p$ \cite{CHM16}. 
These facets are however not exhaustive and, due to the high dimensionality of $\CIM_p$, a complete facet description is only available for small $p$ ($p\leq 4$).
In \cite{LRS20}, lower-dimensional faces of $\CIM_p$ with a more direct combinatorial interpretation were identified.
% In \cite{LRS20} we were given faces with a more direct combinatorial interpretation.
For example, given any undirected graph $G$ we have the face %\cite{LRS20}
\[
\CIM_G\coloneqq \conv\left(c_\GG\colon \GG=([p], E)\textrm{ a DAG with skeleton } G\right).
\]
It was also shown in \cite{LRS20} that generalizations of reversing and adding in edges of DAGs constitute edges of $\CIM_p$.
This raises the question as to whether there is a graphical explanation of other edges as well.
While a complete characterization of the edges of $\CIM_p$ appears challenging, a complete characterization of those edges corresponding to edges of $\CIM_G$ for well-chosen families of $G$ is achievable.

The focus of this paper is on $\CIM_G$ when $G$ is a tree. 
First, in \Cref{subsec: two examples} we will consider the case when $G$ is a star.
Using some standard techniques from the theory of partially ordered sets (posets) we will see that, in this case, $\CIM_G$ is a simplex (see \Cref{col: star graph simplex}).
As trees are locally stars we can utilize this local structure to impose more global conditions for adjacency in the edge graph of $\CIM_G$, resulting in a complete characterization of all edges of $\CIM_G$ for $G$ a tree, both in terms of DAGs (\Cref{thm: subtree condition}) and in terms of essential graphs (\Cref{thm: edges of trees characterization}). 
For completeness, we include \Cref{subsec: disjoint graphs}, specifically \Cref{prop: disjoint graphs}, which generalizes these results to forests.

In \Cref{subsec: stable set polytopes} we also observe a connection between $\CIM_G$ and another well-studied polytope; namely the \emph{stable set polytope}, $\STAB(G)$.
For a given graph $G$, $\STAB(G)$ has vertices corresponding to the stable sets of $G$. 
Chv\'atal gave a characterization of the edges of $\STAB(G)$ in terms of the stable sets \cite{C75}. 
Here, we present a unimodular equivalence between a certain face of the characteristic imset polytope $\CIM_G$ for $G$ a tree and the stable set polytope of a tree (\Cref{prop: stab as face of cim}).
Hence, we recover Chv\'atal's results in these cases.
The connection between the stable set polytope and the characteristic imset polytope also allows us to describe the edges of $\CIM_{C_p}$, where $C_p$ is the cycle with $p$ nodes (\Cref{thm: paths and cycles}), suggesting that there may be a more general connection. 

In \Cref{sec: applications}, we discuss the application of these geometric observations to the problem of causal discovery. 
Given that we have all edges of $\CIM_G$, when $G$ is a tree, we define a new algorithm, \Cref{alg: EFT}, for learning a directed tree (also known as a \emph{polytree}) from data. 
We prove that this algorithm is asymptotically consistent, and the same follows for another algorithm defined in \cite{LRS20} (\Cref{prop: skeletal greedy cim consistent}).
We observe that the additional search capabilities provided by our geometric observations on the edge structure of $\CIM_G$ results in improved performance over classical methods on simulated Gaussian data (\Cref{subsec: computations}).

%%%%%%%%%%%%%%%%%%%%%%%%%%%%%%%%%%%%%%%%%%%%%%%
%---SUBSECTION: PRELIMINARIES -----
\subsection{Preliminaries}
\label{subsec: preliminaries} 

Let $\GG=([p], E)$ be a directed graph. 
We will denote with $i\to j\in \GG$ if $(i,j)\in E$, and if this is the case we say that $i$ is a \emph{parent} of $j$ and $j$ is a \emph{child} of $i$. 
The set of parents and children of $i$ in $\GG$ is denoted with $\pa_\GG(i)$ and $\ch_\GG(i)$, respectively. 
Analogously, if $G=([p], E)$ is undirected we denote $i-j\in G$ if $\{i,j\}\in E$, and say that $i$ and $j$ are \emph{neighbors} in $G$. 
The set of all neighbours of $i$ in $G$ is denoted by $\ne_G(i)$ and the closure of a node $i$ is defined as $\cl_G(i)=\{i\}\cup\ne_G(i)$.

The \emph{skeleton} of a directed graph $\GG$ is the undirected graph $G$ with the same set of nodes, but each edge $i\to j\in \GG$ is replaced with $i-j\in G$. 
A \emph{path} in an undirected graph $G$ is a sequence of distinct vertices $P=(v_0, v_1, \dots, v_n)$ such that $v_{i-1}-v_{i}\in G$ for all $i\in[n]$, and any path in which we allow $v_0=v_1$ is called a \emph{cycle}.
A graph $G$ is \emph{connected} if for every pair of nodes $i$ and $j$ there exists a path $P$ such that $v_0=i$ and $v_n=j$, and $G$ is a \emph{tree} if it is connected and does not contain a cycle. 
A \emph{path} in a directed graph $\GG$ is a path in the skeleton, and the path is \emph{directed} if $v_{i-1}\to v_{i}\in G$ for all $i\in[n]$. 
A directed path in $\GG$ is a \emph{directed cycle} if $v_0=v_n$, and $\GG$ is \emph{connected} if the skeleton is connected. 
For a directed graph $\GG$ we say that $i$ and $j$ are \emph{neighbors} if $i$ and $j$ are neighbors in the skeleton of $\GG$ and denote the corresponding set with $\ne_\GG(i)$.
If there is a directed path $i\to\dots\to j$ in $\GG$ we say that $i$ is an \emph{ancestor} of $j$ and that $j$ is a \emph{descendent} of $i$. 
We denote the set of all ancestors of $i$ in $\GG$ with $\an_\GG(i)$, the set of all descendants with $\de_\GG(i)$, and the set of all \emph{non-descendants} with $\nd_\GG(i)$. 

If $A$ is a subset of the nodes of $\GG$ we let $\GG|_A$ denote the \emph{induced graph} on $A$; that is, $\GG|_A$ has the nodes $A$ and for any $i, j\in A$ we have $i\to j\in \GG|_A$ if and only if $i\to j\in\GG$. 
An induced subgraph of $\GG$ of the form $i\to j\leftarrow k$ is called a \emph{v-structure}. 
A node $i$ is a \emph{leaf} in a tree $G$ if it has a unique neighbor.
A node that is not a leaf is called an \emph{interior node}.
The induced subgraph of $G$ on all interior nodes of $G$ is denoted $G^\circ$.
% An edge $i-j$ is in the interior of a tree if both $i$ and $j$ are interior nodes.
If $G$ is a tree and $D$ is a subset of the vertices, then $\spans_G(D)$ denotes the unique minimal spanning tree of $D$. 

We will also consider \emph{partially directed graphs}; i.e., graphs whose edges can be either directed or undirected. 
Assume $P=(v_0, v_1, \dots, v_k)$ is a path in the skeleton of a partially directed graph, then we will say that $P\to v_k$ in $\GG$ if we have $v_0\to v_1\to \dots\to v_k$ in $\GG$. 
A path $P$ in the skeleton of a DAG $\GG$ from $v_0$ to $v_n$ is \emph{d-connecting given $C$} if for every $v_i$, $0<i<k$ such that $v_{i-1}\to v_i\leftarrow v_{i+1}$ we have $C\cap (\de_\GG(v_i)\cup\{v_i\})\neq \emptyset$, and for every other $v_i$ we have $v_i\notin C$. 
We say that two subsets $A,B\subseteq[p]$ are \emph{d-connected} given a third subset $C\subseteq[p]$ if there exists $i\in A$ and $j\in B$ such that there is a d-connecting path from $i$ to $j$ given $C$ in $\GG$. 
Otherwise, we say $A$ and $B$ are \emph{d-separated} given $C$ in $\GG$. 
% We will regularly confuse any path $P$ with the set of vertices of $P$, as is common practice. 

For any DAG $\GG=([p], E)$ we associate a set of random variables $X_1,\dots, X_p$ and we say that the joint distribution $\mathbb{P}$ of $(X_1,\ldots, X_p)$ is \emph{Markov} to $\GG$ if $\mathbb{P}$ entails $X_i\independent X_{\nd_\GG(i)\setminus \pa_\GG(i)}|X_{\pa_\GG(i)}$ for all $i\in [p]$.
Equivalently $\mathbb{P}$ is Markov to $\GG$ if we have $X_A\independent X_B|X_C$ whenever $A$ and $B$ are d-separated given $C$ in $\GG$ \cite{lauritzen1996}.
Furthermore, $\mathbb{P}$ is \emph{faithful} to $\GG$ if $P$ entails exactly the CI statements encoded by $\GG$. 
It can happen that two DAGs encode the same set of CI statements.
In this case, we say that the DAGs are \emph{Markov equivalent} and that they belong to the same \emph{Markov equivalence class} (MEC). 
The graphical side of Markov equivalence is well-understood with this classical result from Verma and Pearl. 
\begin{theorem}
\label{thm: verma pearl}
\cite{VP92}
Two DAGs are Markov equivalent if and only if they have the same skeleton and the same v-structures.
\end{theorem}
Using DAGs to model complex systems is nowadays common within many areas \cite{FLNP00, P00, BHR00, S01}.
Algorithms for inferring a MEC from data is a well studied area and many algorithms for doing so have been proposed \cite{C02, HB12, SUW20, TBA06, VP92, WSYU17}. 
Many of these algorithms are \emph{score-based}; i.e., given a scoring criterion $S(-,\mathbb{D})$, where $\mathbb{D}$ is a random sample from the joint distribution of $(X_1, \dots, X_p)$, we aim to find the DAG $\GG$ maximizing $S(\GG,\mathbb{D})$. 
A commonly used scoring criterion is the \emph{Bayesian information criterion} (BIC), which is defined in \Cref{sec: applications}. % for a definition, that has been shown to work well in both practice and theory. 

To get a unique representative for each MEC Andersson, Madigan and Perlman introduced essential graphs \cite{AMP97}, a special family of partially directed graphs. 
The \emph{essential graph} of an MEC is the graph that has the same skeleton as each DAG in the MEC, and a directed edge $i\to j$ if and only if we have $i\to j\in \GG$ for all DAGs $\GG$ in the MEC. 
Using this definition the authors gave a complete classification of the directed edges of essential graphs, namely $i\to j$ is directed in an essential graph $\GG$ if and only if it is \emph{strongly protected} in $\GG$ \cite[Theorem 4.1]{AMP97}.  

As an alternative to working with DAGs or essential graphs, Studen\'y described how to use vector encodings to represent CI models \cite{S05}. 
Studen\'y, Lindner, and Hemmecke developed this idea in \cite{SHL10} and introduced the characteristic imset (as defined in the introduction), from which the skeleton and the v-structures are easily recovered. 
%---LEMMA: Imset Structure
\begin{lemma}
\label{lem: imset structure}
\cite{SHL10}
Let $\GG$ be a DAG with nodes $[p]$.
Then for any distinct nodes $i$, $j$, and $k$ we have
\begin{enumerate}[label=(\arabic*)]
\item{$i\leftarrow j$ or $j\rightarrow i$ in $\GG$ if and only if $c_\GG(\{i, j\})=1$.}
\item{$i\rightarrow j \leftarrow k$ is a v-structure in $\GG$ if and only of $c_\GG(\{i,j,k\})=1$ and $c_\GG(\{i,k\})=0$.}
\end{enumerate}
\end{lemma}
Moreover, as the characteristic imset encodes the CI statements we get the following alternative characterization of Markov equivalence.
%----LEMMA: Studeny
\begin{theorem}
\label{lem: studeny}
\cite{SHL10}
Two DAGs $\GG$ and $\HH$ are Markov equivalent if and only if $c_\GG = c_\HH$.
\end{theorem}
As is noted in \cite{SHL10}, any score equivalent, decomposable scoring criterion can be seen as an affine function over the characteristic imsets, motivating the definition of the characteristic imset polytope. 
The following is direct from Theorem~\ref{thm: verma pearl} and \Cref{lem: imset structure}. % the following is direct. 

%---LEMMA: 2/3 sets characterize the imset
\begin{lemma}
\label{lem: lindner}
\cite[Corollary 2.2.6]{L12} %Corollary 2.2.6
Two characteristic imsets $c_\GG$ and $c_\HH$ are equal if and only if $c_\GG(S)=c_\HH(S)$ for all sets $|S|\in\{2,3\}$.
\end{lemma}

%%%%%%%%%%%%%%%%%%%%%%%%%%
%---SUBSECTION: Two examples
\subsection{Two examples}
\label{subsec: two examples}
The general question of characterizing all edges of $\CIM_p$ or $\CIM_G$ seems to be hard in general.
However, in certain cases it can be done. 
Here we will give two such examples that will be relevant for later results. 
As $\CIM_p$ is a $0/1$ polytope, that is every coordinate of the vertices are either $0$ or $1$, it makes sense to consider a change in a single coordinate.

%---DEFINITION: Addition------
\begin{definition}
[Addition]
\label{def: addition}
Let $\GG$ and $\HH$ be two DAGs on node set $[p]$.  
We say the pair $\{\GG,\HH\}$ is an \emph{addition} if $c_\HH = c_\GG+e_{S^\ast}$ for some $S^\ast\subseteq[p]$ with $|S|\geq2$.  
We further say that $\{\GG,\HH\}$ is an \emph{edge addition} if $|S^\ast| = 2$ and a \emph{v-structure addition} if $|S^\ast| = 3$.
\end{definition}

One can equivalently define an addition as $c_\GG$ and $c_\HH$ having Hamiltonian distance $1$ from one another.
As $\CIM_p$ is a 0/1-polytope this also implies that $\conv(c_\GG, c_\HH)$ is an edge whenever $\{\GG, \HH\}$ is an addition. 
The partition into edge- and v-structure additions becomes clearer in light of the following proposition. 
We see that not only are all additions either an edge or v-structure addition, but we also characterize them in terms of the underlying graphs, $\GG$ and $\HH$.

%---PROPOSITION: Addition Characterization-----
\begin{proposition}
\label{prop: addition characterization}
Let $\GG$ and $\HH$ be DAGs on node set $[p]$, and suppose that the pair $\{\GG,\HH\}$ is an addition such that
\[
c_\HH = c_\GG+e_{S^\ast},
\]
where $|S^\ast|\geq 2$. 
Then either $S^\ast = \{i,j\}$ for some $i,j\in[p]$, all v-structures of $\HH$ are present in $\GG$, and the skeletons of $\GG$ and $\HH$ differ by the presence of the edge $\{i,j\}$, or $S^\ast = \{i,j,k\}$ and $\GG$ and $\HH$ have the same skeleton but differ by a single v-structure $i\rightarrow k\leftarrow j$.  
\end{proposition}

\begin{proof}
Suppose first that $|S^\ast| = 2$.  
Then $S^\ast = \{i,j\}$ for some $i,j\in[p]$ for which $c_\GG(\{i,j\}) = 0$ and $c_\HH(\{i,j\}) = 1$, and so the skeletons of $\GG$ and $\HH$ differ by the presence of the edge $\{i,j\}$.    
Since $c_\HH = c_\GG + e_{S^\ast}$, it follows that $c_\HH(S) = c_\GG(S)$ for all $S\subseteq[p]$ with $|S| = 3$.  
By \Cref{lem: imset structure} we have $c_\HH(S) = 1$ for such an $S$ if and only if $S$ is complete in $\HH$ or if the induced subgraph $\HH_S$ is a v-structure.  
Moreover, if the set $S = \{i,j,k\}$ is not complete in the skeleton of $\HH$, then there is at most one way to orient its edges on $\{i,j,k\}$ to produce a v-structure.  
Thus, since $c_\HH(S) = c_\GG(S)$ for all $S\subseteq[p]$ with $|S| = 3$, it follows that $\GG$ and $\HH$ have the same v-structures except for possibly v-structures of the form $i\to k\leftarrow j\in \GG$ for some $k\in[n]$. 

Suppose now that $|S^\ast| = 3$.  
Then $S^\ast = \{i,j,k\}$ for some $i,j,k\in [p]$.  
Since $c_\HH = c_\GG + e_{S^\ast}$, it follows that $\GG$ and $\HH$ have the same skeleton and that $c_\GG(S^\ast) = 0$.  
Hence, we know that $S^\ast$ is not complete in either of $\GG$ or $\HH$.  
From this, and the fact that $c_\HH(S^\ast) = 1$, it follows that the induced subgraph $\HH|_{S^\ast}$ of $\HH$ is a v-structure.
Since $c_\GG(S) = c_\HH(S)$ for all $S\neq S^\ast$ with $|S| = 3$, a similar argument to the previous case shows that all other $v$-stuctures in $\GG$ and $\HH$ are the same. 

Finally suppose that $|S^\ast| \geq 4$. 
That gives us $c_\GG(S)=c_\HH(S)$ for all $|S|\leq 3$ and by \Cref{lem: lindner} we have $c_\GG=c_\HH$, a contradiction.
\end{proof}
From this it also follows that edge additions are a special type of \emph{edge pair}, and that v-structure additions are a special type of \emph{turn pair}, as defined in \cite{LRS20}. 
We will also see how v-structure additions show up naturally in \Cref{subsec: stable set polytopes}.

In the next case we instead show that restricting the skeleton, as opposed to restricting a relation between the characteristic imsets, can also lead to cases where we can describe the edges. 
This next example will also be useful later when we consider trees. 
First we will give one lemma. 
Given a finite set $P$, let $\mathbb{R}[P]$ denote the vector space of dimension $|P|$ where the basis vectors are indexed by the elements of $P$. 
%---LEMMA: Basis in PO vector space ------
\begin{lemma}
\label{lem: basis in po vector space}
Given a finite partial order $(P,\preceq)$. 
Define $b_p$ for any $p\in P$ as
\[
b_p\coloneqq \sum_{q\preceq p}e_q 
\]
where $e_q$ denotes the standard basis of $\mathbb{R}[P]$.
Then $\{b_p\}_{p\in P}$ is a basis for $\mathbb{R}[P]$.
\end{lemma}

\begin{proof}
By the Möbius inversion formula $e_p = \sum_{q\preceq p}\mu_P(q,p) b_q$. 
Thus $b_p$ is the image of $e_p$ under an invertible linear transformation. 
The result follows.
\end{proof}

Then we consider the following case:

%---PROPOSITION: Star-like graph ------
\begin{proposition}
\label{prop: star like graph simplex}
Assume $G'=([p-1], E')$ is a vertex disjoint union of complete graphs $K_{p_1},\dots,K_{p_m}$.
Let $G=([p], E)$ be the graph with $\ne_G(p)=[p-1]$ and the induced graph $G|_{[p-1]} = G'$. 
Then $\CIM_G$ is a $d$-simplex with $d= 2^{p-1}-1-\sum_{s\in[m]}\left(2^{p_s}-1\right)$.
\end{proposition}

\begin{proof}[Proof of \Cref{prop: star like graph simplex}]
We begin by counting the number of MECs.
As the MEC is determined by the skeleton and v-structures, and we have fixed the skeleton, we only need to consider sets $\{i,j,k\}$ such that $G|_{\{i,j,k\}} = i-j-k$.
From the assumptions on $G$ all triples of this form have $j=p$, $i\in K_{p_s}$, and $k\in K_{p_{s'}}$ where $s\neq s'$.

Thus, the MEC of $\GG$ is completely determined by $\pa_\GG(p)$. 
If $|\pa_\GG(p)|\leq 1$ we have no v-structures, and all such DAGs are Markov equivalent. 
Likewise, if $\pa_\GG(p)\subseteq K_{p_{s}}$ for some $s\in[m]$.
Hence, the number of MECs are counted as $1$ plus the number of ways to choose $\pa_\GG(p)$ such that $|\pa_\GG(p)|\geq 2$ and $\pa_\GG(p)\not\subseteq K_{p_s}$ for any $s\in[m]$. 
Using the fact that $\sum_{s\in[m]}p_s=p-1$, this quantity is
\[
1+2^{p-1}-(p-1) -1 -\sum_{s\in[m]}\left(2^{p_s}-p_s-1\right)=2^{p-1}-\sum_{s\in[m]}\left(2^{p_s}-1\right).
\]

Next we wish to show that all $c_\GG$ are affinely independent.
Let $\DD$ be a DAG such that $\pa_\DD(p)=\emptyset$ and skeleton $G$, then $\DD$ has no v-structures.
Using the above, a straightforward calculation then gives us that for any DAG $\GG$ with skeleton $G$ we get 
\[
c_\GG = c_\DD + \sum_{S\in \SSS}  e_S
\]
where $\SSS = \left\{T \cup \{p\}\colon T\subseteq \pa_\GG(p), \centernot\exists s\in[m]\textrm{ such that }T\subseteq K_{p_s}\right\}$. 
Notice that we do not have any sets of size $2$ or less in $\SSS$.
Moreover, the parents of $p$ uniquely determine the MEC unless we are in the class containing $\DD$. 
Thus, 
\[
P\coloneqq \left\{T\cup\{p\}\colon |S|\geq 2, \centernot\exists s\in[m]\textrm{ such that } T\subseteq K_{p_s}\right\}\cup \{\emptyset\}
\]
can be partially ordered by inclusion and for every DAG $\GG$ there exists a unique $S\in P$ such that 
$
c_\GG=c_\DD+\sum_{T\in P\colon T\subseteq S} e_T
$. 
Thus we can apply \Cref{lem: basis in po vector space} to $\left\{c_\GG-c_\DD\colon \GG\text{ a DAG with skeleton $G$}\right\}$, and the result follows.
\end{proof}
Applying the above proposition with $p_i=1$ for all $i$ we get the following lemma:

%---COROLLARY: Star-like graph ------
\begin{lemma}
\label{col: star graph simplex}
If $G=([p], E)$ is a star, then $\CIM_G$ is a $2^p-p-1$-simplex.
\end{lemma}

The proof of \Cref{prop: star like graph simplex} is possible due to the fact that the skeleton imposes a one-to-one correspondence between the MECs and the principal order ideals of a poset. 
In \Cref{subsec: almost complete graphs} we will see another example of when this method can be applied, however, it seems that this is a rather rare quality of $G$.

All edges we have encountered thus far are, in a sense, local in the graph; i.e., all known edges only depend on vertices that are, or will become, neighbours. 
The question then arises whether all edges of $\CIM_p$ are local, or if there are edges $\conv(c_\GG, c_\HH)$ where the relation between $\GG$ and $\HH$ depends on some global structure.

%%%%%%%%%%%%%%%%%%%%%%%%%%%%%%%%%%%%%%%%%%
%---SECTION: Trees ------
\section{Trees}
\label{sec: trees}
In this section we will examine the edge structure of $\CIM_G$ when $G$ is a tree.
As we saw in \Cref{col: star graph simplex}, $\CIM_G$ is a simplex when $G$ is a star.
In this section we will see how this can be used to impose a local structure on essential graphs whose skeleton is a tree. 
With this local structure we then impose more global structures and are able to recover all edges of $\CIM_G$. 
Our first characterization will be in terms if essential graphs. 
As this characterization might be hard to deal with in practice, we will give a characterization in terms of DAGs as well (see \Cref{subsec: essential flips as dags}).

%---SUBSECTION: Essential Flips ---
\subsection{The Essential Side of the Trees}
\label{subsec: essenstial flips}
The authors of \cite{AMP97} introduced the essential graph as a unique representative of an MEC. 
As edges of a polytope, in a sense, represent a minimal change between vertices, the question of finding edges of $\CIM_G$ becomes ``what does a minimal change of the essential graph look like?".
To answer this question we introduce \emph{essential flips} (\Cref{def: essential flip}), and show that these indeed characterize the edges of $\CIM_G$ when $G$ is a tree. 

As previoously mentioned, in \cite{AMP97} the authors gave a complete characterization of the directed edges of essential graphs.
However, in the case of trees, the criterion simplifies significantly. 
In the later proofs we will frequently make use of the following proposition.
%---Proposition: Characterization of Essential Edges in Trees---
\begin{proposition}
\label{prop: characterization of essential edges}
Let $\GG$ be a DAG whose skeleton is a tree.  An arrow $i\rightarrow j$ is essential in $\GG$ if and only if either
\begin{enumerate}
    \item there is a node $k$ such that $i \rightarrow j \leftarrow k$ is an induced subgraph of $\GG$, or 
    \item $i$ is a descendant of a node with (at least) two parents. 
\end{enumerate}
\end{proposition}

\begin{proof}
Suppose first that there exists a node $k$ such that $i \rightarrow j \leftarrow k$ is an induced subgraph of $\GG$.  
As this forms a v-structure, it follows that $i\rightarrow j$ is essential.  
Suppose next that there is an ancestor $k$ of $i$ such that $k$ has two parents, say $p_1$ and $p_2$.  
Assume for the sake of contradiction that $i\rightarrow j$ is non-essential.  
Since $i$ is a descendant of $k$, there exists a directed path from $k$ to $j$ in $\GG$: $k \rightarrow k_1\rightarrow \cdots \rightarrow k_m\rightarrow i \rightarrow j$. 
If $i\rightarrow j$ is non-essential then there exists an element of the Markov equivalence class of $\GG$, say $\HH$, in which $i\leftarrow j$.  
Since $\GG$ has skeleton a tree and $k_m \rightarrow i \rightarrow j$ in $\GG$ is not a v-structure, then it is also not a v-structure in $\HH$.  
Since $\HH$ must have the same skeleton as $\GG$, it must be that $k_m\leftarrow i$ in $\HH$.  
Iterating this argument up the directed path from $k$ to $j$ in $\GG$, we get that $\HH$ contains the directed path $k \leftarrow k_1\leftarrow \cdots \leftarrow k_m\leftarrow i \leftarrow j$.  
However, $k$ has two parents, $p_1$ and $p_2$, and, as the skeleton of $\GG$ was a tree, hence $\GG$ contains the v-structure $p_1\rightarrow k\leftarrow p_2$.
Since the triples $p_1,k,k_1$ and $p_2,k,k_1$ do not form v-structures in $\GG$ they must also not form v-structures in $\HH$.  
Hence, it must be that $p_1\leftarrow k$ and $p_2\leftarrow k$ in $\HH$, meaning that $\HH$ lacks a v-structure that is in $\GG$, a contradiction.  

Conversely, suppose that $i\rightarrow j$ is essential in $\GG$. 
Then $i\rightarrow j$ is strongly protected in the essential graph $D$ of $\GG$.  
Since $\GG$ a tree it follows that wither (a) $i\rightarrow j\leftarrow k_0$ or (b) $k_0\rightarrow i\rightarrow j$ is an induced subgraph of $D$ for some node $k_0$.  
In the case that we have (a), we are done. 
In the case that we have only (b), we know that $k_0\rightarrow i$ is essential in $D$ (as it is directed in $D$).  
Hence, iterating this argument assuming that the new edge is always strongly protected as in (b), we find a directed path in $D$: 
$k_m\rightarrow k_{m-1}\rightarrow \cdots \rightarrow k_0 \rightarrow i \rightarrow j$ where $\pa_D(k_m) = \emptyset$ and $\pa_D(k_i) = \{k_{i+1}\}$ for $i = 0,\ldots, m-1$.
In this case, $k_m\rightarrow k_{m-1}$ is not strongly protected in $D$, which is a contradiction.
\end{proof}
We then have the following lemma.
%---LEMMA: Essential arrows at a node---
\begin{lemma}
\label{lem: essential arrows at a node}
Let $\GG$ be a DAG whose skeleton is a tree. If $i\rightarrow j$ is an essential arrow of $\GG$ then every edge with $j$ as an endpoint is essential in $\GG$.
\end{lemma}

\begin{proof}
If $i\rightarrow j$ is essential and $k$ is adjacent to $j$ then either $i\rightarrow j \leftarrow k$ is a v-structure in $\GG$ or $i\rightarrow j \rightarrow k$ is an induced subgraph of $\GG$.  
In the latter case, by Proposition~\ref{prop: characterization of essential edges}, since $i\rightarrow j$ is essential, either we have a v-structure at $j$ or $i$ is a descendant of a node with two non-adjacent parents in $\GG$.  
Regardless of the case, it follows from Proposition~\ref{prop: characterization of essential edges}~(2) that the arrow $i\rightarrow k$ is also essential in $\GG$. 
\end{proof}

Based on \Cref{lem: essential arrows at a node}, given $j\in\ne_\GG(i)$, we say that $j$ is an \emph{essential parent (resp., child)} of $i$ if $j\rightarrow i$ (resp., $j\leftarrow i$) and the edge $j\rightarrow i$ (resp., $j\leftarrow i$) is essential.
Recall that we are interested in relations between DAGs or essential graphs that give rise to edges of the $\CIM_G$ polytope. 
The following definition is then natural in this context. 

\begin{definition}
\label{def: delta set}
Let $\GG$ and $\HH$ be two essential graphs with skeleton $G=([p], E)$ and assume $G$ is a tree. 
Let $N_i\coloneqq \left\{S\subseteq [p]\colon i\in S\subseteq \ne_G(i)\cup\{i\}, |S|\geq 3\right\}$.
Then we define 
\[
\Delta(\GG, \HH)\coloneqq \left\{i\in[p]\colon c_\GG|_{N_i}\neq c_\HH|_{N_i}\right\}.
\]
\end{definition}
Notice that $\Delta(\GG, \HH)$ is empty if and only of $\GG$ and $\HH$ are Markov equivalent. 
The set $\Delta(\GG, \HH)$ is natural to consider since it can equivalently be defined as all nodes $j$ such that we have a v-structure $i\to j\leftarrow k$ in $\GG$ but not in $\HH$, or vice versa. 

\begin{lemma}
\label{lem: delta set}
Let $\GG$ and $\HH$ be two essential graphs on node set $[p]$.  Suppose $\GG$ and $\HH$ both have the same skeleton $G$ and that $G$ is a tree.  
Then $j\in \Delta(\GG, \HH)$ if and only if $\GG$ and $\HH$ contain different sets of v-structures centered at $j$.  
\end{lemma}

\begin{proof}
Suppose there is a v-structure centered at $j$ in $\GG$ that is not in $\HH$, say $i\rightarrow j\leftarrow k$.  
Then $\{k,j,i\}\in N_j$ and 
$
c_\GG(\{k,j,i\}) = 1 \neq 0 = c_\HH(\{k,j,i\}).
$
Thus, $j\in \Delta$.

Conversely, suppose that $j\in \Delta$.  
Then there exists $S\in N_j$ such that $c_\GG(S)\neq c_\HH(S)$. 
Without loss of generality, assume $c_\GG(S) = 1$ and $c_\HH(S) = 0$.  
Since $j\in S$ and $G$ is a tree, then $c_\GG(S) = 1$ if and only if $k\rightarrow j$ for all $i\in S\setminus\{j\}$.  
Thus, since $|S|\geq 3$, $\GG$ contains a v-structure that is not in $\HH$.  
\end{proof}

Now we have an understanding of both the essential graphs and how a change in the characteristic imset affects changes in v-structures and vice versa. 
With this we can present the relation of particular interest of this section. 

%---DEFINITION: Essential flip----
\begin{definition}
[Essential flip]
\label{def: essential flip}
Let $\GG$ and $\HH$ be two non-Markov equivalent essential graphs with skeleton $G$, a tree, and denote $\Delta=\Delta(\GG, \HH)$. 
Assume that both $\GG|_{\spans(\Delta)}$ and $\HH|_{\spans(\Delta)}$ do not contain any undirected edges. 
Assume moreover that each edge of $\GG$ and $\HH$ differ on $G|_{\spans(\Delta)}$. 
Then we say that the pair $\{\GG, \HH\}$ is an essential flip. 
\end{definition}

For convenience we will say that a pair of DAGs constitutes an essential flip if their essential graphs constitute an essential flip.  
Let us consider an example of an essential flip. 

\begin{example}
\label{ex: essential flip}
In \Cref{fig: essential flip ex} we give two graphs $\GG$ and $\HH$ such that $\{\GG, \HH\}$ constitutes an essential flip.
Here we have $\Delta(\GG, \HH)=\{\delta_i\}_{i=1}^{7}$. 
We see that edges outside $G|_{\spans(\Delta(\GG, \HH))}$ can change both direction and essentiality. 
Moreover, the inclusion $\Delta(\GG, \HH)\subseteq \spans(\Delta(\GG, \HH))$ can be strict, here with $n_3\in  \spans(\Delta(\GG, \HH)) \setminus\Delta(\GG, \HH)$. 
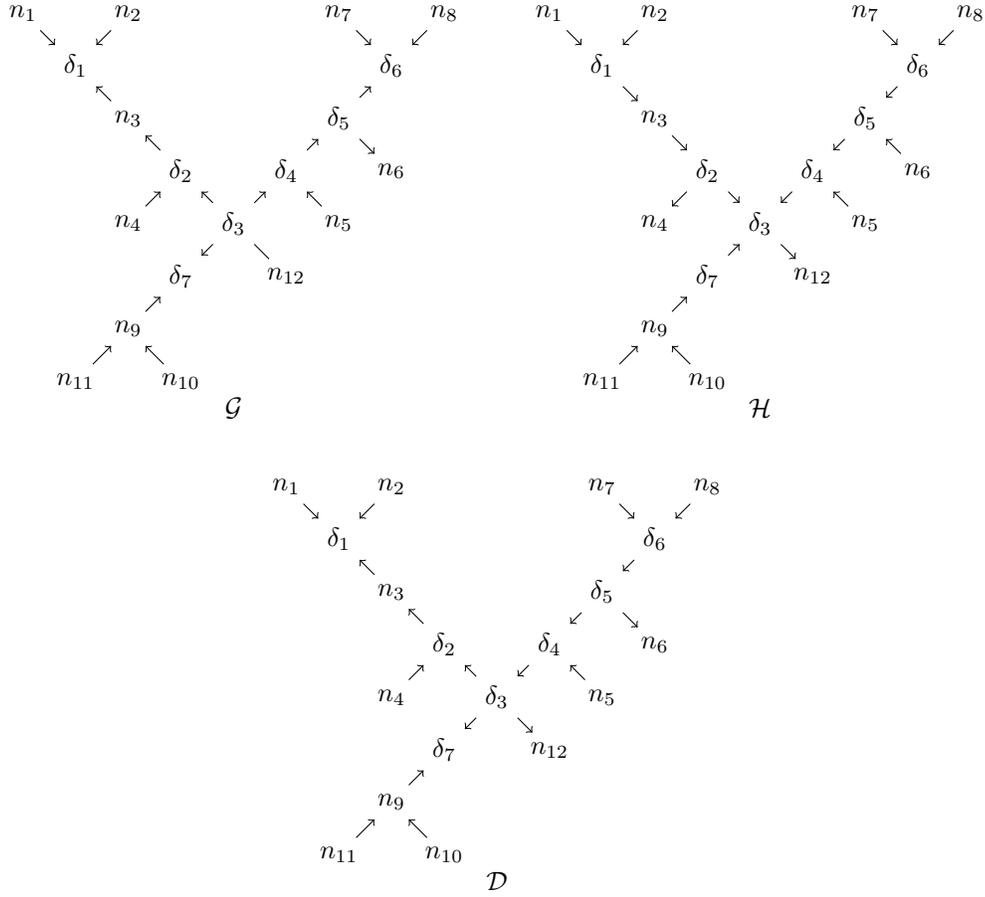
\begin{figure}
\[
\begin{tikzpicture}[scale=.7]
\begin{scope}[]
\node (d1) at (-3,3) {$\delta_1$};
\node (d2) at (-1,1) {$\delta_2$};
\node (d3) at (0,0) {$\delta_3$};
\node (d4) at (1,1) {$\delta_4$};
\node (d5) at (2,2) {$\delta_5$};
\node (d6) at (3,3) {$\delta_6$};
\node (d7) at (-1,-1) {$\delta_7$};

\node (n1) at (-4,4) {$n_1$};
\node (n2) at (-2,4) {$n_2$};
\node (n3) at (-2,2) {$n_3$};
\node (n4) at (-2,0) {$n_4$};
\node (n5) at (2,0) {$n_5$};
\node (n6) at (3,1) {$n_6$};
\node (n7) at (2,4) {$n_7$};
\node (n8) at (4,4) {$n_8$};
\node (n9) at (-2,-2) {$n_9$};
\node (n10) at (-1,-3) {$n_{10}$};
\node (n11) at (-3,-3) {$n_{11}$};
\node (n12) at (1, -1) {$n_{12}$};

\node at (0, -3.5) {$\GG$};

\foreach \from/\to in {n3/d1, d2/n3, d3/d2, n4/d2, d3/d4, d4/d5, d3/d7, d5/d6, d5/n6}{
  \draw[->] (\from) -- (\to);
}

\foreach \from/\to in {n1/d1, n2/d1, n5/d4, n7/d6, n8/d6, n9/d7, n10/n9, n11/n9}{
  \draw[->] (\from) -- (\to);
}

\draw (n12) -- (d3);
\end{scope}

\begin{scope}[shift={(10,0)}]
\node (d1) at (-3,3) {$\delta_1$};
\node (d2) at (-1,1) {$\delta_2$};
\node (d3) at (0,0) {$\delta_3$};
\node (d4) at (1,1) {$\delta_4$};
\node (d5) at (2,2) {$\delta_5$};
\node (d6) at (3,3) {$\delta_6$};
\node (d7) at (-1,-1) {$\delta_7$};

\node (n1) at (-4,4) {$n_1$};
\node (n2) at (-2,4) {$n_2$};
\node (n3) at (-2,2) {$n_3$};
\node (n4) at (-2,0) {$n_4$};
\node (n5) at (2,0) {$n_5$};
\node (n6) at (3,1) {$n_6$};
\node (n7) at (2,4) {$n_7$};
\node (n8) at (4,4) {$n_8$};
\node (n9) at (-2,-2) {$n_9$};
\node (n10) at (-1,-3) {$n_{10}$};
\node (n11) at (-3,-3) {$n_{11}$};
\node (n12) at (1, -1) {$n_{12}$};

\node at (0, -3.5) {$\HH$};

\foreach \from/\to in {n3/d1, d2/n3, d3/d2, n4/d2, d3/d4, d4/d5, d3/d7, d5/d6, d5/n6}{
  \draw[<-] (\from) -- (\to);
}

\foreach \from/\to in {n1/d1, n2/d1, n5/d4, n7/d6, n8/d6, n9/d7, n10/n9, n11/n9}{
  \draw[->] (\from) -- (\to);
}

\draw[<-] (n12) -- (d3);
\end{scope}

\begin{scope}[shift={(5, -9)}]
\node (d1) at (-3,3) {$\delta_1$};
\node (d2) at (-1,1) {$\delta_2$};
\node (d3) at (0,0) {$\delta_3$};
\node (d4) at (1,1) {$\delta_4$};
\node (d5) at (2,2) {$\delta_5$};
\node (d6) at (3,3) {$\delta_6$};
\node (d7) at (-1,-1) {$\delta_7$};

\node (n1) at (-4,4) {$n_1$};
\node (n2) at (-2,4) {$n_2$};
\node (n3) at (-2,2) {$n_3$};
\node (n4) at (-2,0) {$n_4$};
\node (n5) at (2,0) {$n_5$};
\node (n6) at (3,1) {$n_6$};
\node (n7) at (2,4) {$n_7$};
\node (n8) at (4,4) {$n_8$};
\node (n9) at (-2,-2) {$n_9$};
\node (n10) at (-1,-3) {$n_{10}$};
\node (n11) at (-3,-3) {$n_{11}$};
\node (n12) at (1, -1) {$n_{12}$};

\node at (0, -3.5) {$\DD$};

\foreach \from/\to in {n3/d1, d2/n3, d3/d2, n4/d2, d4/d3, d5/d4, d3/d7, d6/d5, d5/n6}{
  \draw[->] (\from) -- (\to);
}

\foreach \from/\to in {n1/d1, n2/d1, n5/d4, n7/d6, n8/d6, n9/d7, n10/n9, n11/n9}{
  \draw[->] (\from) -- (\to);
}

\draw[<-] (n12) -- (d3);
\end{scope}
\end{tikzpicture}
\]
\caption{An example of an essential flip $\{\GG, \HH\}$. 
Here we have $\Delta(\GG, \HH)=\{\delta_i\}_{i=1}^7$.
Here we also give $\DD$ as an example for the proof of \Cref{prop: edge of tree}, it can however also be checked that both $\{\GG, \DD\}$ and $\{\HH, \DD\}$ also constitute essential flips. 
See \Cref{ex: essential flip} and \Cref{ex: proof essential flip ex} for more details. 
}
\label{fig: essential flip ex}
\end{figure}
\end{example}

While the definition of essential flips may seem unintuitive they give us a complete characterization of edges of the $\CIM_G$ polytope, whenever $G$ is a tree. 

\begin{theorem}
\label{thm: edges of trees characterization}
If $G$ is a tree, then $\conv(c_\GG, c_\HH)$ is an edge of $\CIM_G$ if and only if the pair $\{\GG, \HH\}$ is an essential flip. 
\end{theorem}
This is the result of \Cref{prop: edge of tree} and \Cref{prop: not edge of tree} below, each showing one way of the equivalence.

\begin{proposition}
\label{prop: edge of tree}
If $\{\GG, \HH\}$ is an essential flip where $\GG$ and $\HH$ are essential graphs with skeleton $G$, and $G$ is a tree, then $\conv(c_\GG, c_\HH)$ is an edge of $\CIM_G$. 
\end{proposition}
As the following proof is rather technical we advise the reader to keep \Cref{ex: essential flip} and \Cref{ex: proof essential flip ex} in mind.

\begin{proof}
Throughout this proof we will say that an essential graph $\DD$ \emph{looks like} $\GG$ (or $\HH$) at $i$ if we have $c_\DD|_{N_i}=c_\GG|_{N_i}$ (or $c_\DD|_{N_i}=c_\HH|_{N_i}$). 
Equivalently, via \Cref{lem: delta set}, $\DD$ looks like $\GG$ at $i$ if and only if $\DD$ and $\GG$ share the same v-structures centered at $i$. 
We will find successively smaller faces of $\CIM_G$, each containing both $c_\GG$ and $c_\HH$. 
Each of these faces will be defined via cost functions $W_1, \dots, W_6$ that successively imposes more restrictions on the essential graphs encoding for vectors in the faces. 
In some sense, $W_1$ will fix all coordinates that $c_\GG$ and $c_\HH$ share, $W_2$ will make sure that around each node we will look like either $\GG$ or $\HH$, and $W_3$ will direct all paths in $G|_{\spans(\Delta)}$.
$W_4$ will make sure that for each non-endpoint of paths, $i\in\spans(\Delta)$ we have that if we locally look like $\GG$ around $i$, we have a path $P\subseteq \spans(\Delta)$ with $i\in P$ that is directed as in $\GG$, and similarly if we look like $\HH$.
Using $W_5$ we will make all paths be directed as in either $\GG$ or $\HH$, as opposed to some directed as in $\GG$ and some as in $\HH$, and $W_6$ will take care of the endpoints of paths in a similar fashion as $W_4$ took care of non-endpoints. 
% After we have constructed $W_1,\dots, W_6$ we will apply \Cref{lem: weight lemma} and define $W$. 
Crucially, in our construction we will have score functions $g_i$ and $h_i$ that are indicator functions of the v-structures of $\GG$ and $\HH$, respectively, around a node $i$. 
These are constructed via \Cref{col: star graph simplex}.
We will continuously use \Cref{lem: essential arrows at a node}, which states that if we have an essential edge $i\gets j$ in $\GG$, then $\GG|_{\cl(i)}$ is fully directed.
We will also say that an essential graph $\DD$ maximizes a cost function $W$ if $c_\DD$ maximizes $W^Tc_\DD$. 

We begin by defining
\[
W_1(S)\coloneqq\left\{
\begin{array}{ll}
1&\text{if }c_\GG(S)=c_\HH(S)=1,\\
-1&\text{if }c_\GG(S)=c_\HH(S)=0, \text{ and}\\
0&\text{otherwise.}\\
\end{array}
\right.
\]
Then any essential graph $\DD$ maximizing $W_1$ must have $c_\DD(S)=c_\GG(S)$ whenever $c_\GG(S)=c_\HH(S)$. 
% Then we note that by \Cref{prop: case 2 3 simplex} $\CIM_T$ is a simplex whenever $T$ is a star. 

Restricting $\CIM_G$, via projection, to $\spans(e_S\colon S\in N_i)$ is the same as considering $\CIM_{G|_{\cl(i)}}$, as the value of the characteristic imset $c_\DD(S)$ only depends on $\DD|_S$. 
As $G$ is a tree, $G|_{\cl(i)}$ is a star and thus, by \Cref{col: star graph simplex}, $\CIM_{G|_{\cl(i)}}$ is a simplex. 
Hence, for each node $i\in [p]$ there exists an affine function $w_i$ such that $w_i$ only depends on $N_i$ and maximizes on $\conv(c_\GG|_{N_i}, c_\HH|_{N_i})$. 
Moreover, as $\CIM_{G|_{\cl(i)}}$ is a simplex, we can assume that $w_i^Tc_\GG=w_i^Tc_\HH=1$ and $w_i^Tc_\DD=-1$ if $c_\DD|_{N_i}\neq c_\GG|_{N_i}$ and $c_\DD|_{N_i}\neq c_\GG|_{N_i}$.
Then let $W_2=\sum_{i\in[p]}w_i$.
It follows that $W_2^Tc_\GG=W_2^Tc_\HH=p$ and that $W_2^Tc_\DD < p-1$ if we do not have $c_\DD|_{N_i}=c_\GG|_{N_i}$ or $c_\DD|_{N_i}=c_\HH|_{N_i}$ for every $i\in[p]$. 
Note that if $N_i=\emptyset$, that is $|\ne_G(i)|\leq 1$, for some $i\in[p]$ we get that $w_i$ is the constant function $1$.

Similarly, for every $i\in\Delta$ we can let $g_i$ and $h_i$ be affine functions that only depend on $N_i$ and satisfy 
\[
h_i^Tc_\HH=g_i^Tc_\GG=1,\ g_i^Tc_\HH=h_i^Tc_\GG=0,\]
and $g_i^Tc_\DD=h_i^Tc_\DD=0$, if $c_\DD|_{N_i}\neq c_\GG|_{N_i}$ and $c_\DD|_{N_i}\neq c_\HH|_{N_i}$.
For every $i\in[p]\setminus\Delta$ we let $g_i = h_i$ be affine functions, again only dependent on $N_i$, such that $h_i^Tc_\HH=g_i^Tc_\GG=1$ and $g_i^Tc_\DD=0$ if $c_\DD|_{N_i}\neq c_\GG|_{N_i}= c_\HH|_{N_i}$.
Then $g_i$ and $h_i$ work as indicator functions for an essential graph looking like $\GG$ or $\HH$ around $i$, in terms of v-structures.

Let $\mathcal{C}$ be the set of all directed paths in $\GG|_{\spans(\Delta)}$. 
Since $\GG$ and $\HH$ differ on every edge of $G|_{\spans(\Delta)}$, this is also the set of directed paths in $\HH|_{\spans(\Delta)}$.
Let $\mathcal{C}$ be partially ordered by inclusion.
Later in the proof, a maximal path will refer to a maximal path in $\mathcal{C}$ with respect to this partial order. 
Take $P=(v_0, v_1, \dots, v_n)$ to be a maximal element of $\mathcal{C}$, by symmetry we can assume that we have $v_0\to v_1\to\dots\to v_n$ in $\GG$.
Notice that by definition of essential flip we have $v_0\leftarrow v_1\leftarrow\dots\leftarrow v_n$ in $\HH$. 

We will now construct $W_3$ and show that $c_\GG$ and $c_\HH$ maximize $W_3$ under the assumption that we maximize $W_1$ and $W_2$, which both $c_\GG$ and $c_\HH$ do. 
Therefore, assume that $\DD$ is an essential graph maximizing $W_1$ and $W_2$.
% At least one such graph exists since both $\GG$ and $\HH$ have this property. 
By construction of $W_1$ and $W_2$ this is the same as either $c_\DD|_{N_k}=c_\GG|_{N_k}$ or $c_\DD|_{N_k}=c_\HH|_{N_k}$ for all  $k\in[p]$. 
Since the edge $v_0\to v_1$ is directed in $\GG$, by \Cref{prop: characterization of essential edges}, we have three cases;
{\bf Case I:} there exists a v-structure $i\to v_0\leftarrow k$ in $\GG$;
{\bf Case II:} we have a v-structure $i\to j \leftarrow k$ with $j\to \dots\to v_0$ in $\GG$; 
{\bf Case III:} we have a v-structure $v_0\to v_1\leftarrow k$ for some $k\in[p]$.
These cases are not exclusive, however at least one must be true. 
\begin{enumerate}[label={\bf Case \Roman*:}]
\item{%As $P$ was maximal we must have $i,k\notin \Delta$. 
%Hence $\GG$, $\HH$ and $\DD$ share the v-structure, as $\DD$ was assumed to maximize $W_1$. 
If $c_\DD|_{N_{v_0}}=c_\GG|_{N_{v_0}}$ we must have $c_\DD(\{i, v_0, k\})=c_\GG(\{i, v_0, k\})=1$ and hence $\DD$ must have the v-structure $i\to v_0\leftarrow k$ as well. 
If we have $v_0\leftarrow v_1$ in $\DD$ we have the v-structure $i\to v_0\leftarrow v_1$ in $\DD$ but not in $\GG$, hence $c_\DD|_{N_{v_0}}\neq c_\GG|_{N_{v_0}}$. 
Thus we conclude that if $c_\DD|_{N_{v_0}}=c_\GG|_{N_{v_0}}$ then $v_0\to v_1$ in $\DD$.
}
\item{If $j\in \Delta$ we then have a directed path $j\to \dots \to v_0\to\dots\to v_n$ bigger than $P$ in $\GG|_{\spans(\Delta)}$, contradicting the maximally of $P$ in $\mathcal{C}$. 
In the same way we get that $j'\notin \Delta$ for all $j'$ on the path between $j$ and $v_0$. 
Hence $c_\GG|_{N_{i'}}=c_\HH|_{N_{i'}}=c_\DD|_{N_{i'}}$ for all $i'$ between $i$ and $v_0$, including $i$ but not $v_0$ as $\DD$ maximized $W_1$.
Thus we have an essential edge $\to v_0$ in $\GG$, $\HH$ and $\DD$.
We can by \Cref{lem: essential arrows at a node} say that if $c_\DD|_{N_{v_0}}=c_\GG|_{N_{v_0}}$ we must have $v_0\to v_1$ in $\DD$. 
}
\item{If $c_\DD|_{N_{v_1}}=c_\GG|_{N_{v_1}}$, then $\DD$ must have the v-structure $v_0\to v_1\leftarrow k$. 
Hence if $c_\DD|_{N_{v_1}}=c_\GG|_{N_{v_1}}$ we must have $v_0\to v_1$ in $\DD$. 
Notice the difference in index from the previous cases, $v_1$ instead of $v_0$. 
}
\end{enumerate}

Hence no matter the case there exists a node $\alpha\in\Delta$ such that $c_\DD|_{N_{\alpha}}=c_\GG|_{N_{\alpha}}$ implies that $v_0\to v_1\in\DD$, under the assumption that $\DD$ maximizes both $W_1$ and $W_2$. 
Similarly there exists a node $\beta\in\Delta$ such that $c_\DD|_{N_{\beta}}=c_\HH|_{N_{\beta}}$ implies that $v_{n-1}\leftarrow  v_n\in\DD$. 
For each maximal path $P$, we fix a pair of these nodes as $\alpha_P$ and $\beta_P$. 
If possible, as it is in {\bf Case  I} and {\bf Case II} above, we choose $\alpha_P$ and $\beta_P$ as endpoints of $P$.
In \Cref{ex: proof essential flip ex} we have given an example how these vertices are chosen.  
% This will help us in constructing $W_6$. 

For each maximal path in $P=v_0\to\dots\to v_n$ in $\GG|_{\spans(\Delta)}$ we let $w_P=-(g_{\beta_P}+h_{\alpha_P})$. 
Define $W_3=\sum w_P$ where we sum over all maximal paths in $\mathcal{C}$.
% Notice that $|\mathcal{C}|\leq \binom{p}{2}<\frac{p^2}{2}$ and hence $\frac{-1}2 <  W_3^Tc_\DD\leq 0$. 
By construction of $w_P$ we have $w_P^Tc_\GG=w_P^Tc_\HH=-1$. 
Again assuming $\DD$ maximizes $W_1$ and $W_2$, we want to prove that $W_3^Tc_\DD\leq W_3^Tc_\GG=W_3^Tc_\HH$.
Then it is enough to show that $w_P^Tc_\DD\leq w_P^Tc_\GG=w_P^Tc_\HH$ for every $P$. 
If this is not the case, by construction of $w_p$, we have $g_{\beta_P}^Tc_\DD=h_{\alpha_P}^Tc_\DD=0$. 
As $\DD$ maximized $W_1$ and $W_2$ we must have $c_\DD|_{N_{\beta_P}}=c_\HH|_{N_{\beta_P}}$ and $c_\DD|_{N_{\alpha_P}}=c_\GG|_{N_{\alpha_P}}$. 
If the length of $P$ is greater than $1$ ($n\geq 2$) we must have $v_0\to \dots \leftarrow v_n$ in $\DD$, but as neither $\GG$ or $\HH$ has a v-structure along the path $v_0-\dots-v_n$, neither can $\DD$, a contradiction. 
If $n=1$ the edge $v_0-v_1$ must be directed in both directions in $\DD$, which cannot happen. 
Hence, any essential graph $\DD$ such that $c_\DD$ maximizes $W_1$ and $W_2$ has $W_3^Tc_\DD\leq W_3^Tc_\GG= W_3^Tc_\HH$. 
If $w_P^Tc_\DD = -2$ for some $P$, we get $ W_3^Tc_\DD\leq W_3^Tc_\GG-1<W_3^Tc_\GG$. 
For all maximal paths $P\in\mathcal{C}$, as neither $\GG$ or $\HH$ has a v-structure along $P$, neither can $\DD$. 
Moreover, if $w_P^Tc_\DD = -1$ we have we must have $g_{\beta_P}^Tc_\DD=0$ or $h_{\alpha_P}^Tc_\DD=0$, and as $\DD$ maximized $W_2$ we get $h_{\beta_P}^Tc_\DD=1$ or $g_{\alpha_P}^Tc_\DD=1$, respectively. 
Hence we get $v_0\to v_1$ or $v_{n-1}\leftarrow v_n$ in $\DD$. 
Thus $\DD|_P$ must be a directed path, assuming $c_\DD$ maximizes $W_1$, $W_2$ and has $W_3^Tc_\DD=W_3^Tc_\GG$. 
This is what we meant when we said that $W_3$ directs all paths in $G|_{\spans(\Delta)}$ either as in $\GG$ or $\HH$.

We will now construct $W_4$ that will make sure that if $i$ is not an endpoint of any maximal path in $\mathcal{C}$ and a graph $\DD$ looks like $\GG$ at $i$, then there is a maximal path $P$ containing $i$ such that $\DD|_P=\GG|_P$, and similar for $\HH$. 
This will be similar to how $W_3$ directed all paths.  % agrees with the direction of a path it is in. 
For example $\delta_5$ in $\DD$, \Cref{fig: essential flip ex}, does not have this property; 
Despite the maximal path $(\delta_3, \delta_4, \delta_5, \delta_6)$ being directed as in $\HH$, see \Cref{fig: essential flip ex}, we have $c_{\DD_1}|_{N_{\delta_5}}=c_{\GG}|_{N_{\delta_5}}$.

To avoid this, we consider all nodes $i\in\Delta$ such that $i$ is not the endpoint of any maximal path $P\in \mathcal{C}$. 
As $i\in\Delta$ there must exist a v-structure $k\to i\leftarrow j$ in $\GG$ or $\HH$.
Assume it is in $\GG$. 
Fix a maximal path $P$ with $i\in P$. 
Notice as $i\in P$, and $i$ has at least one v-structure in $\GG$, there is at least one v-structure at $i$ containing an edge of $P$ in $\GG$. 
This v-structure cannot be present if $P$ is directed as in $\HH$. 
Let $w'_i=-(h_i+g_{\beta_P})$. 
We claim that we cannot have $w_{i}{'}^{T}c_\DD=0$ if $c_\DD$ maximizes $W_1$, $W_2$ and has $W_3^Tc_\DD=W_3^Tc_\GG$.
Then  we would have $h_i^Tc_\DD=0$ and $g_{\beta_P}^Tc_\DD=0$, or, as $c_\DD$ maximized $W_1$ and $W_2$, $g_i^Tc_\DD=1$ and $h_{\beta_P}^Tc_\DD=1$. 
By the above argument we than get $P|_\DD=P|_\HH$, as well as the v-structure in $\GG$, but as mentioned above this cannot happen.
% Hence $P$ must be directed as in $\HH$ in $\DD$ as well, contradicting that we have $c_\DD|_{S_i}=c_\GG|_{S_i}$. 
If there was no v-structure at $i$ in $\GG$, there must have been a v-structure at $i$ in $\HH$, and hence we let $w'_i=-(g_i+h_{\alpha_P})$ and proceed similarly. % with $j$ chosen similarly.
Let $W_4=\sum w'_{i}$ where we sum over all $i\in\Delta$ that are not endpoints of any maximal paths in $\GG|_{\spans(\Delta)}$. 
Thus we conclude that any essential graph $\DD$ that maximizes $W_1$, $W_2$, has $W_3^Tc_\DD=W_3^Tc_\GG$, and $W_4^Tc_\DD=W_4^Tc_\GG$, must have for each maximal element $P\in\mathcal{C}$ that either $\DD|_P=\GG|_P$ and for each node $i\in P\cap \Delta$, that is not the beginning or end of $P$, $c_\DD|_{N_i}=c_\GG|_{N_i}$, or similar for $\HH$. 

The summands of $W_5$ are a bit more intricate than the ones we have seen before, therefore we will begin to construct these. 
Assume that $\DD$ maximizes $W_1$, $W_2$, has $W_3^Tc_\DD=W_3^Tc_\GG$, and $W_4^Tc_\DD=W_4^Tc_\GG$.
Then, by our above reasoning, for any two maximal paths $P_1, P_2\in \mathcal{C}$ we have $\DD|_{P_i}$ are directed paths and coincides with either $\GG|_{P_i}$ or $\HH|_{P_i}$. 
Assume $P_1\cap P_2$ is non-empty. 
Then if $|P_1\cap P_2|\geq 2$, $P_1$ and $P_2$ must share an edge, as $G$ was assumed to be a tree. 
Then there can only exist two ways of fully directing $P_1\cup P_2$, which must coincide with either $\GG$ or $\HH$. 
If $P_1\cap P_2=\{i\}$ then we either have $c_\DD|_{N_i}=c_\GG|_{N_i}$ or $c_\DD|_{N_i}=c_\HH|_{N_i}$. 
By maximality of $P_1$ and $P_2$ we get 3 cases, either $i$ is the source of both $P_1$ and $P_2$ in $\GG$, $i$ is the endpoint of both $P_1$ and $P_2$ in $\GG$, or $i$ is a node in the middle of both $P_1$ and $P_2$. 
We note that if we have at least one v-structure in $\cl_G(i)$ in both $\GG$ and $\HH$, then $\cl_G(i)$ must be fully directed whenever we have $c_\DD|_{N_i}=c_\GG|_{N_i}$ or $c_\DD|_{N_i}=c_\HH|_{N_i}$. 
Thus we must have $\DD|_{P_1\cup P_2}=\GG|_{P_1\cup P_2}$ or $\DD|_{P_1\cup P_2}=\HH|_{P_1\cup P_2}$. 
In the third case we must have v-structures in both $\GG$ and $\HH$. 
Hence the interesting case is when $i$ is the source of both $P_1$ and $P_2$ in either $\GG$ or $\HH$. 
If it is the source in $\GG$, let $w_{P_1, P_2} = (2h_i-h_{\beta_{P_1}}-h_{\beta_{P_2}})$. 
Otherwise $i$ is the source of $P_i$ and $P_2$ in  $\HH$, and we let $w_{P_1, P_2} = (2g_i-g_{\alpha_{P_1}}-g_{\alpha_{P_2}})$. %, again with $j_1$ and $j_2$ as above. 
Notice that $w_{P_1, P_2}^Tc_\GG=w_{P_1, P_2}^Tc_\HH=0$ regardless of which case. 

Now that the summands are constructed we can show that these behave as intended,
we do however advise the reader to again consider \Cref{ex: proof essential flip ex}.

Let us begin to show that $w_{P_1, P_2}^Tc_\DD\leq w_{P_1, P_2}^Tc_\GG=0$ for all appropiate $\DD$.
If we were to have $w_{P_1, P_2}^Tc_\DD > 0$ we must have $h_i^Tc_\DD=1$ (or $g_i^Tc_\DD=1$, but it is the same up to symmetry).
That is, we have $P_1\to i \leftarrow P_2$ in $\DD$. 
%Let $v_1$ be the node that directs $P_1$ as $\GG$.
As $\DD$ was assumed to have $W_3^Tc_\DD=W_3^Tc_\GG$ we must have $(g_{\beta_{P_1}}+h_{\alpha_{P_1}})^Tc_\DD=1$.  
Since we had that $P_1$ is directed as in $\HH$, we cannot have $g_{\alpha_{P_1}}^Tc_\DD=1$.
By assumption $c_\DD$ maximized $W_2$, hence we then must have $h_{\alpha_{P_1}}^Tc_\DD=1$, and by the above $g_{\beta_{P_1}}^Tc_\DD=0$.
Again by the assumption that $c_\DD$ maximizes $W_2$ we get $h_{\beta_{P_1}}^Tc_\DD=1$. 
We can repeat this argument with $P_2$ and get $h_{\beta_{P_2}}^Tc_\DD=1$. 
Hence $w_{P_1, P_2}^Tc_\DD=0$, a contradiction. 
Therefore we let $W_5=\sum w_{P, Q}$ where we sum over all pairs of maximal paths in $\mathcal{C}$ whose intersection is a single point. 
Then, by the above, we get $W_5^Tc_\DD\leq W_5^Tc_\GG=W_5^Tc_\HH$, assuming that $\DD$ maximizes $W_1$, $W_2$, has $W_3^Tc_\DD=W_3^Tc_\GG$, and $W_4^Tc_\DD=W_4^Tc_\GG$.
If $\DD$ would have two maximal paths $P_1$ and $P_2$ such that $P_1\cap P_2\neq \emptyset$ but $\DD|_{P_1}=\GG|_{P_1}$ and $\DD|_{P_2}=\HH|_{P_2}$, then by the above we must have that $P_1$ and $P_2$ intersect in a single point $i$ such that $i$ is the middle of a v-structure in either $\GG$ or $\HH$. 
We can without loss of generality assume it is a v-structure in $\HH$ and thus a source in $\GG$. 
Then we must have $P_1\to i\to P_2$ in $\DD$. 
We defined $w_{P_1, P_2}=(2h_i-h_{\beta_{P_1}}-h_{\beta_{P_2}})$.
As we do not have the v-structure present in $\HH$ at $i$ in $\DD$ we must have $h_i^Tc_\DD=0$, but as argued above we must have $h_{\beta_{P_2}}^Tc_\DD=1$.
Hence, $w_{P_1, P_2}^Tc_\DD< w_{P_1, P_2}^Tc_\GG$. 

So far we have constructed $W_1,\dots,W_5$ that direct all paths in $\mathcal{C}$ either as in $\GG$ or $\HH$. % and makes sure that the characteristic imsets agree with $\GG$ or $\HH$ for all internal nodes of all paths in $\mathcal{C}$. 
With $W_4$, we made sure that if $c_\DD|_{N_i}=c_\GG|_{N_i}$ for any node $i$ that is not the endpoint of a maximal path in $\mathcal{C}$ there exists a path $P$ with $i\in P$ such that $\DD|_P=\GG|_P$. 
Now we want similar guarantees for the endpoints as well. 
Due to our choice to have $\alpha_P$ and $\beta_P$ as endpoints whenever possible, $W_3$ does exactly that in {\bf Cases I} and {\bf II}.
We will see this later in the proof. 
Thus let us consider the endpoints of maximal paths that had $\alpha_P$ or $\beta_P$ from {\bf Case III} above.
Assume there exists a node $i\in\Delta$ such that $i$ is the endpoint of a path $i\to P$ in $\GG$ (or $\HH$) and we do not have any v-structures in $\GG|_{\cl(i)}$. 
Then we let $w''_P=-(g_i+h_{\beta_P})$ and $W_6=\sum w''_P$. 
Assume $\DD$ is an essential graph such that $c_\DD$ maximizes $W_1$, $W_2$, has $W_3^Tc_\DD=W_3^Tc_\GG$, $W_4^Tc_\DD=W_4^Tc_\GG$, and $W_5^Tc_\DD= W_5^Tc_\GG$. 
By assumption we must have $i\leftarrow P$ in $\HH$ and we must moreover have a v-structure in $\HH|_{\cl(i)}$ as $i\in\Delta$. 
If $g_i^Tc_\DD=0$ by $W_2$ we get $h_i^Tc_\DD=1$.
Thus the first edge of $P$ must be directed as in $\HH$ in $\DD$.
Hence, we cannot have $g_{\alpha_P}^Tc_\DD=1$ as that would direct $P$ as in $\GG$. 
By $W_2$ we get $g_{\alpha_P}^Tc_\DD=0$, and thus $h_{\alpha_P}^Tc_\DD=1$.
By $W_3$ we get $g_{\beta_P}^Tc_\DD=0$ and again by $W_2$ we get $h_{\beta_P}^Tc_\DD=1$. 
Hence $W_6^Tc_\DD\leq W_6^Tc_\GG=W_6^Tc_\HH$.

Now that we have constructed $W_1, \dots, W_6$ we are in a position to finish the proof. 
So far we have proved that there is a face $F$ of $\CIM_G$ such that $c_\DD\in F$ if and only if $\DD$ maximizes $W_1$, $W_2$, has $W_3^Tc_\DD=W_3^Tc_\GG$, $W_4^Tc_\DD=W_4^Tc_\GG$, $W_5^Tc_\DD= W_5^Tc_\GG$, and $W_6^Tc_\DD= W_6^Tc_\GG$. 
Now we want to show $F$ only contains the vertices $c_\GG$ and $c_\HH$, that is $F=\conv(c_\GG, c_\HH)$. 
By construction of $W_1,\dots, W_6$ it follows that $\conv(c_\GG, c_\HH)\subseteq F$. 

Let $\DD$ be an essential graph such that $c_\DD\in F$; that is, $\DD$ maximizes $W_1$, $W_2$, has $W_3^Tc_\DD=W_3^Tc_\GG$, $W_4^Tc_\DD=W_4^Tc_\GG$, $W_5^Tc_\DD= W_5^Tc_\GG$, and $W_6^Tc_\DD= W_6^Tc_\GG$. 
%Let $\DD$ be such an essential graph. 
Then for all $i\notin\Delta$ we have $c_\DD|_{N_i}=c_\GG|_{N_i}=c_\HH|_{N_i}$, as $c_\DD$ maximized $W_1$. 
For any fixed $i\in\Delta$, as $c_\DD$ maximized $W_2$ we have that either $c_\DD|_{N_i}=c_\GG|_{N_i}$ or $c_\DD|_{N_i}=c_\HH|_{N_i}$. 
By symmetry we can assume that $c_\DD|_{N_i}=c_\GG|_{N_i}$. 

We want to show that there exists a maximal path $P\in\mathcal{C}$ with $i\in P$ such that $\DD|_P=\GG|_P$. 
If $i$ is not the endpoint of a maximal path in $\mathcal{C}$, we can choose $P$ to be the same path as when we constructed $w'_i$ (see the construction of $W_4$).
As $c_\DD$ maximized $W_1$, $W_2$, had $W_3^Tc_\DD=W_3^Tc_\GG$, and $W_4^Tc_\DD=W_4^Tc_\GG$ we can apply the same argument as there and obtain $P$. 

Otherwise $i$ is the endpoint of some maximal path $P\in\mathcal{C}$. 
If we have a v-structure in $\GG|_{\cl(i)}$, then $\GG|_{\cl(i)}$ must be fully directed. 
As $c_\DD|_{N_i}=c_\GG|_{N_i}$ we get that $\DD|_{\cl(i)}$ is fully directed, specifically the edge in $\cl(i)$ belonging to $P$. 
Thus $P$ must be directed the same way in both $\GG$ and $\DD$. 
If we do not have a v-structure in $\GG|_{\cl(i)}$ we must have a v-structure in $\HH|_{\cl(i)}$, as $i\in\Delta$. 
Then we can divide this into two cases, either $P\to i$ in $\GG$ or $P\leftarrow i$ in $\GG$. 
In the case of $P\to i$ in $\GG$ we have $P\leftarrow i$ in $\HH$, as $P\in\mathcal{C}$ and thus $P\subseteq\spans(\Delta)$.
Thus following our cases we must have chosen $\beta_{P}=i$. 
As $W_3^Tc_\DD=W_3^Tc_\GG$ we have $(g_{i} + h_{\alpha_{P}})^Tc_\DD=1$ and since $c_\DD|_{N_i}=c_\GG|_{N_i}$ we must have $h_{\alpha_{P}}^Tc_\DD=0$. 
Then we get that $g_{\alpha_{P}}c_\DD=1$, thus $\DD|_{P}=\GG|_{P}$ and hence $P$ is a valid choice. 
In the case of $P\leftarrow i$ in $\GG$ we are exactly in the case where we created $w''_P$. 
Then we have $(g_i+h_{\beta_P})^Tc_\DD=1$, which gives us $h_{\beta_P}^Tc_\DD=0$, by $W_2$ we get $g_{\beta_P}^Tc_\DD=1$, by $W_3$ we get $h_{\alpha_P}^Tc_\DD=0$, and again by $W_2$ we get $g_{\alpha_P}^Tc_\DD=1$. 
Hence $P$ must be directed as in $\GG$ in $\DD$. 

Now the proof reduces to taking another node $j$ and showing that $c_\DD|_{N_j}=c_\GG|_{N_j}$. 
As $c_\DD|_{N_j}=c_\GG|_{N_j}=c_\HH|_{N_j}$ for all nodes $j\notin\Delta$, as $\DD$ maximized $W_1$ and $W_2$, we let $j\in\Delta\setminus\{i\}$. 
For the sake of contradiction, if $c_\DD|_{N_j}=c_\HH|_{N_j}$, by the above, we have two maximal paths, in $\mathcal{C}$, $P_i$ and $P_j$, containing $i$ and $j$ respectively, such that $P_i$ is directed as in $\GG$ and $P_j$ is directed as in $\HH$. 
As $G$ is connected we can find maximal paths $\{Q_k\}_{k=1}^t\in\mathcal{C}$ such that $P_i=Q_1$, $P_j=Q_t$ and $Q_{k}\cap Q_{k+1}\neq \emptyset$ for all $1\leq k\leq t-1$. 
As shown above when constructing $W_5$, since $c_\DD\in F$ we must have that $Q_{k}\cup Q_{k+1}$ is directed as in either $\GG$ or $\HH$.
However, we have $Q_1$ directed as in $\GG$ and $Q_{k}\cap Q_{k+1}\neq \emptyset$. 
Inductively we get that all $Q_k$ are directed as in $\GG$, a contradiction since $Q_t=P_j$ is directed as in $\HH$.

Thus the only characteristic imsets in $F$ are $c_\GG$ and $c_\HH$. 
Hence $\conv(c_\GG, c_\HH)$ is an edge of $\CIM_G$.  
\end{proof}

\begin{example}
\label{ex: proof essential flip ex}
In \Cref{ex: essential flip} we gave an example of an essential flip $\{\GG, \HH\}$, see \Cref{fig: essential flip ex}.
Following the proof of \Cref{prop: edge of tree} we have three maximal paths in $\mathcal{C}$, $P_1=(\delta_3, \delta_2, n_3, \delta_1)$, $P_2=(\delta_3, \delta_4, \delta_5, \delta_6)$, and $P_3=(\delta_3, \delta_7)$ with $\alpha_{P_1}=\delta_2$ (Case III), $\beta_{P_1}=\delta_1$ (Case I), $\alpha_{P_2}=\delta_4$ (Case III), $\beta_{P_2}=\delta_6$ (Case I), and $\alpha_{P_3}=\beta_{P_3}=\delta_7$ (Case III and Case II, respectively). 
Notice that if we would not have the convention of choosing endpoints whenever possible we could have chosen $\beta_{P_2}=\delta_5$. 

The essential graph $\DD$ looks like $\GG$ or $\HH$ at every node, but is not Markov equivalent to either. It is straightforward to see that $\DD$ maximizes $W_1$, $W_2$, and has $W_3^Tc_\DD=W_3^Tc_\GG$. 
However, as $P_2$ is directed as in $\HH$ in $\DD$ but $c_\DD|_{N_{\delta_5}}=c_\GG|_{N_{\delta_5}}$, and $\delta_5$ is only in one maximal path, $c_\DD$ does not maximize $W_4$. 

Moreover, in $\DD$ the path $P_1$ is directed as in $\GG$ but the path $P_2$ is directed as in $\HH$. 
As the intersection between $P_1$ and $P_2$ is a single point, $\DD$ cannot maximize $W_5$ either. 
Indeed, as $\delta_3$ is the source of $P_1$ and $P_2$ in $\GG$ we defined $w_{P_1, P_2}=2h_{\delta_3}-h_{\delta_1}-h_{\delta_6}$. 
Then it follows that $w_{P_1, P_2}^Tc_\DD=-1<0=w_{P_1, P_2}^Tc_\GG=w_{P_1, P_2}^Tc_\HH$.

In our construction of $W_6$ when considering $\delta_3$ we could have chosen any one path of $P_1$, $P_2$, or $P_3$. 
Then, following the construction of $W_6$ have  $w_{P_1}''=-(g_{\delta_3}+h_{\delta_1})$, $w_{P_2}''=-(g_{\delta_3}+h_{\delta_6})$, and $w_{P_3}''=-(g_{\delta_3}+h_{\delta_7})$ as summands in $W_6$.
As $\DD$ looks like $\GG$ at $\delta_3$, both have no v-structures, but the path $P_2$ is directed as in $\HH$ we will have $w_{P_1}''^Tc_\DD=w_{P_3}''^Tc_\DD=-1$ and $w_{P_2}''^Tc_\DD=-2$. 
Notice that $w_{P_i}''^Tc_\GG=w_{P_i}''^Tc_\HH=-1$ for all maximal paths $P_i$.
\end{example}

Thus essential flips give rise to edges of $\CIM_G$, and in fact they give us a complete characterization. 
To show this we will use the following well-known fact. 
% The following is well known.%\todo{at least I think so}.
\begin{lemma}
\label{lem: square not edge}
Let $P$ be a polytope and let $v$ be a vertex of $P$. 
If there exists non-zero vectors $u_1$ and $u_2$ such that $v+u_1$, $v+u_2$, and $v+u_1+u_2$ are all vertices of $P$, then $\conv(v, v+u_1+u_2)$ is not an edge of $P$. 
\end{lemma}
Then the converse of \Cref{prop: edge of tree} follows as well.
\begin{proposition}
\label{prop: not edge of tree}
If $\GG$ and $\HH$ are essential graphs such that $\{\GG, \HH\}$ is not an essential flip, then $\conv(c_\GG, c_\HH)$ is not an edge of $\CIM_G$, where $G$ is a tree. 
\end{proposition}

\begin{proof}
Recall that $\Delta(\GG, \HH)= \{i\in[p]\colon c_\GG|_{N_i}\neq c_\HH|_{N_i}\}$. 
Since $\{\GG,\HH\}$ is not an essential flip, by symmetry there either exists an undirected edge in $i-j\in\GG|_{\spans(\Delta)}$ or we have an edge $i\to j\in\GG|_{\spans(\Delta)}$ such that $i\to j\in\HH|_{\spans(\Delta)}$.

{\bf Case I, $i-j\in\GG|_{\spans(\Delta)}$:} 
Take any DAG $\HH'$ in the Markov equivalence class of $\HH$. 
By symmetry we can assume that $i\to j\in\HH'$. 
Since $i-j$ was undirected in $\GG$ there exists a DAG $\GG'$ Markov equivalent to $\GG$ with $i\to j\in \GG'$. 
Let $C_1$ be the nodes in the connected component in $\GG'\setminus\{i\to j\}$ containing $i$ and $C_2$ be the complement of $C_1$. 
Let $\DD_1$ be the DAG such that $\DD_1|_{C_1}=\HH'|_{C_1}$, $\DD_1|_{C_2}=\GG'|_{C_2}$, and $i\to j\in \DD_1$. 
Let $\DD_2$ be the DAG such that $\DD_2|_{C_1}=\GG'|_{C_1}$, $\DD_2|_{C_2}=\HH'|_{C_2}$, and $i\to j\in \DD_2$. 
Letting $u_1=c_{\DD_1}-c_\HH$ and $u_2=c_{\DD_2}-c_\HH$, we claim that $c_\GG=c_\HH+u_1+u_2$, or, equivalently, $c_{\GG'}=c_{\HH'}+u_1+u_2$.
By \Cref{lem: lindner} it is enough to show that $c_\GG(S)=c_{\DD_1}(S)+c_{\DD_2}(S)-c_\HH(S)$ holds for all sets $S$ of size $2$ and $3$. 
Since all $\GG$, $\DD_1$, $\DD_2$, and $\HH$ all share the same skeleton the equality is true for all sets of size $2$.

Then, for any set $S=\{s_1, s_2, s_3\}$, if we do not have that $G|_S\simeq P_3$, we have $c_\GG(S)=c_{\DD_1}(S)=c_{\DD_2}(S)=c_\HH(S)=0$ and thus the equality holds. 
If $S\subseteq C_1$ or $S\subseteq C_2$ then either $c_{\DD_1}(S)=c_{\HH'}(S)$ and $c_{\DD_2}(S)=c_{\GG'}(S)$, or $c_{\DD_1}(S)=c_{\GG'}(S)$ and $c_{\DD_2}(S)=c_{\HH'}(S)$, respectively.
This follows from the construction of $\DD_1$ and $\DD_2$. 
The final case is when $\{i,j\}\subseteq S$, which follows similarly as $\GG'$, $\DD_1$, $\DD_2$, and $\HH'$ all have the same direction of $i\to j$, again by construction. 

All that is left to check in  \Cref{lem: square not edge} is that $u_1$ and $u_2$ are non-zero.
Equivalently we can say that $\DD_1$ is not Markov equivalent to either $\GG$ or $\HH$. 
As $i-j\in \spans_G(\Delta)$ and $G$ was a tree we must have $\Delta\cap C_1\neq \emptyset$ and $\Delta\cap C_2\neq \emptyset$. 
Thus there exists nodes $\alpha\in \Delta\cap C_1$ and $\beta\in\Delta\cap C_2$ such that $c_{\DD_1}|_{S_\alpha}=c_{\HH}|_{S_\alpha}\neq c_{\HH}|_{S_\alpha}$ and $c_{\DD_1}|_{S_\beta}=c_{\HH}|_{S_\beta}\neq c_{\HH}|_{S_\beta}$. 
Hence $u_1\neq 0$ and similarly for $u_2$. 

{\bf Case II, $i\to j\in\GG|_{\spans(\Delta)}$ and $i\to j\in\HH|_{\spans(\Delta)}$:} 
In this case we can take any two DAGs $\GG'$ Markov equivalent to $\GG$ and $\HH'$ Markov equivalent to $\HH$. 
Notice that $i\to j\in\GG'$ and $i\to j\in\HH'$. 
Then we can repeat the exact same construction as in {\bf Case I}. 
\end{proof}

A remarkable fact about these polytopes is that every non-edge is of the form of \Cref{lem: square not edge}. 
This is something rather unusual even for $0/1$-polytopes and fails already in dimension $3$. % , for example if $P=\conv(0, e_1, e_2, e_3, e_1+e_2+e_3)$. 
The question if this is true for all $\CIM_G$ polytopes is open.

%---SUBSECTION: The directed side of the story---
\subsection{The Directed Side of the Trees}
\label{subsec: essential flips as dags}
In the previous subsection we gave a description of the edges of $\CIM_G$ in terms of essential graphs. 
However, we are interested in describing transformations of a DAG $\GG$ that produce a DAG $\HH$ such that the essential graphs of $\GG$ and $\HH$ constitute an essential flip.  
By definition of essential flips, the difference between two essential graphs $\GG$ and $\HH$ is a connected subtree. 
Thus the difference between two DAGs that constitute an essential flip can only change v-structures in one unique subtree, and all other differences cannot change the essential graph. 
Thus if two DAGs consitute an essential flip we can assume they differ on a subtree $T$. 
The following theorem gives a characterization in terms of every internal node $i$ of $T$. Note that we have a symmetry between $\GG$ and $\HH$ given by e.g. $T\cap \pa_\GG(i)=T\cap \ch_\HH(i)$ and $ \pa_\GG(i)\setminus T= \pa_\HH(i)\setminus T$.

\begin{theorem}
\label{thm: subtree condition}
Suppose that $\GG$ and $\HH$ are DAGs with the same skeleton $G$ that is a tree. Assume the edges that differ between $\GG$ and $\HH$ form a subtree $T$ of $G$. 
Suppose further that $\Delta(\GG,\HH)\neq\emptyset$. Then the essential graphs of $\HH$ and $\GG$ form an essential flip if and only if each internal node $i$ of $T$ satisfy the conditions given below. 
We use notation $\{\child_i\}=T\cap \ch_\GG(i)$ and $\{\parent_i\}=T\cap \pa_\GG(i)$, when these sets are singletons.

\smallskip
\noindent
\begin{tabular}{| l | l | l | l |}
\hline
& $| T\cap \pa_\GG(i)| $ & $|T\cap \ch_\GG(i)| $ & Local criteria for $\GG$ and $\HH$ to form essential flip\\
\hline
\hline
I & $\ge 2$ &$\ge 2$ & \\
\hline
II & $\ge 2$ &$0$ & \\
\hline
III& $0$ & $\ge 2$& \\
\hline
IV & $\ge 2$ &$1$ & $ | \pa_\GG(i)\setminus T|\ge 1$ or \\
 & & & if $\exists$ v-structure at $\child_i$ in $\GG$, then \\
 &&& $\child_i$ has essential parent in $\HH$ \\
\hline
V & $1$ &$\ge 2$ & $ | \pa_\GG(i)\setminus T|\ge 1$ or \\
 & & & if $\exists$ v-structure at $\parent_i$ in $\HH$, then \\
 &&& $\parent_i$ has essential parent in $\GG$ \\
\hline
VI & $1$ &$1$ & if there are nodes of  $\Delta$ in both\\ 
&&& connected components of $T\setminus \{i\}$\\
& & & then $|\pa_\GG(i)\setminus T|\ge 1$ or  $\child_i$\\
%& & & $\exists$ v-structure at $j$ in $\GG$ or v-structure at $k$ in $\HH$ but if \\
& & & has essential parent in $\HH$ and $\parent_i$ \\ 
&&& has essential parent in $\GG$.\\
\hline
\end{tabular}
\end{theorem}

\begin{proof} Note that all vertices of $T$ will be of exactly one of the types I-VI or a leaf of $T$.
We will make extensive use of \Cref{lem: delta set} in this proof. In particular, it implies that $\Delta$ is a subset of the vertices of $T$, as vertices in $T$ are the only spots in which $\GG$ and $\HH$ could differ in the presence of v-structures.

By Definition \ref{def: essential flip}
of essential flip we must prove the condition that every edge on a path between two nodes in $\Delta$ is essential in both DAGs. We start be proving that this is true for all edges of the form $i\leftarrow j$, where $i$ is of type I-V. 
If $\left| T\cap \pa_\GG(i)\right|\ge 2 $, then there is a v-structure at $i$ in $\GG$ not in $\HH$ so $i\in\Delta$ and 
 by Lemma~\ref{lem: essential arrows at a node}, it follows that all edges incident to $i$ in $\GG$ are essential. Symmetrically, if $\left| T\cap \ch_\GG(i)\right|\ge 2 $ then $i\in\Delta$ and all edges incident to $i$ are essential in $\HH$. 
This implies that vertices of type I-V are always in $\Delta$. It also implies that all edges incident to a node of type I are essential in both $\GG$ and $\HH$, and hence that case is settled.

Let $i-j$ be any edge in $T$ and assume by symmetry it is directed $i\leftarrow j$ in $\GG$ and assume first that $i$ is of type II or IV and hence essential in $\GG$. 
We now go through all possibilities for $j$. Note that $j$ cannot be of type II since $i\in\ch_\GG(j)$. If $j$ is of type III or V, then $i\to j$ is essential also in $\HH$ by the previous paragraph. If instead $j$ is a leaf in $T$ then $j\in \Delta$ if and only if $i\to j$ is part of 
v-structure in $\HH$. But then $i\to j$ is essential in $\HH$ so the condition in Definition  \ref{def: essential flip} is true for the edge $i-j$ in either case. The remaining possibility is that $ T\cap\ch_\GG(j)=\{i\}$, that is $j$ is of type IV or VI, with $\child_j=i$. If the first local criterion of type IV is true, that is $| \pa_\GG(j)\setminus T|\ge 1$ then $i\to j$ is part of a 
v-structure at $j$ in $\HH$ and hence essential. If $j$ is of type IV and $| \pa_\GG(j)\setminus T|=0$, the second local criterion for type IV gives that $i$ has an essential parent in $\HH$ and thus $i\to j$ is an essential edge. 
%Please note that now $i$ and $j$ has switched roles on the condition for type 
%IV.  
If $j$ is of type VI, and there is a node of $\Delta$ in both connected components of $T\backslash\{j\}$ then the same reasoning as for type IV is valid. If there is no node
of $\Delta$ in the other part of  $T\backslash\{j\}$, then $i\to j$ does not need to be essential.

Note that $i$ cannot be of type III, since $i\leftarrow j$. Assume now $i$ is of type V, then we know that $i$ has at least two parents in $\HH$ and thus all edges incident to $i$ in $\HH$, including $i\to j$, are essential in $\HH$.  By the local criteria in type V, we have either a parent of $i$ outside $T$ making $i\leftarrow j$ essential in $\GG$ or if there exist a v-structure at $j=\parent_i$ in $\HH$ there is an essential parent of $j$ in $\GG$, which again makes $i\leftarrow j$ essential in $\GG$. We will again go through all possibilities for $j$. If $j$ is a leaf in $T$ then $j$ is in $\Delta$ if and only if 
there is a v-structure at $j$ in $\HH$ and thus $i \leftarrow j$ is essential in $\GG$ as desired.
If $j$ is of type I or IV, then every edge incident with
$j$ is essential in $\GG$ including $i\leftarrow j$. Note that $j$ has a child in $\GG$ and cannot be of type II. If $j$ is of type III or V, 
%{\color{red} FIGURE} 
then $\ch_\GG(j)=\pa_\HH(j)$ has at least two elements and thus $i\to j$ is part of a v-structure at $j$ in $\HH$ and by the local criterion for type V at $i$ we thus have an essential parent of $j$ in $\GG$
which implies $i \leftarrow j$ is essential in $\GG$. The last possibility is that $j$ is of type VI. If $|\pa_\GG(j)\setminus T |\ge 1$ then there is a v-structure at $j$ in $\GG$ and thus $i\leftarrow j$ is essential in $\GG$. If $| \pa_\GG(j)\setminus T|=0$ then $j\notin\Delta$ and
$i\leftarrow j$ need to be essential only if there are nodes of $\Delta$ on both sides of $T\setminus \{j\}$, in which case the 
last local criterion for type VI is a reformulation of the conditions in Proposition \ref{prop: characterization of essential edges}.

Finally, assume $i$ is of type  VI. If $| \pa_\GG(i)\setminus T|\ge 1$ then there is a v-structure at $i$ in both DAGs and all edges incident to $i$ are essential. If $| \pa_\GG(i)\setminus T|=0$, then $i$ is not in $\Delta$. The edge 
$i-j$ needs to be essential only if there are nodes of $\Delta$ on both sides of $T\setminus \{i\}$, in which case the 
conditions for type VI is a reformulation of the conditions in Proposition \ref{prop: characterization of essential edges}. We have thus proved that the local criteria are sufficient.

For necessity consider first a vertex $i$ of type IV that does not fulfill the local criteria. That is, $| \pa_\GG(i)\setminus T|=0$, and there is a v-structure at $\child_i$ in $\GG$ but no essential parent of $\child_i$ in $\HH$. The v-structure at $\child_i$ implies that $\child_i\in\Delta$ and $| T\cap\pa_\GG(i)|\ge 2$ means $i\in\Delta$. But since $\pa_\HH(i)=\{\child_i\}$ and $\child_i$ has no essential parent the edge $i\leftarrow \child_i$ is not essential in $\HH$. Thus $\GG$ and $\HH$ do not form an essential flip. Type V is symmetric to type IV with the roles of $\GG$ and $\HH$ interchanged.

The remaining case is if $i$ is of type VI and does not fulfill the local criteria. That is, first that there are nodes of $\Delta$ in both connected components of $T\setminus\{i\}$ which means that the two edges incident to $i$ in $T$ must both be essential for $\GG$ and $\HH$ to form an essential flip. And secondly that $| \pa_\GG(i)\setminus T|=0$, and there is no essential parent of $\child_i$ in $\HH$ or no essential parent of $\parent_i$ in $\GG$. By symmetry we can assume the former and then as in type IV conclude that the edge $i\leftarrow \child_i$ is not essential in $\HH$. Thus also the local criteria for type VI are necessary.
\end{proof}

\renewcommand{\algorithmicrequire}{\textbf{Input:}}
\renewcommand{\algorithmicensure}{\textbf{Output:}}

%---SECTION: Examples ------
\section{Examples}
\label{sec: examples}

We will now consider a few examples of how \Cref{thm: edges of trees characterization} can be applied. 
We will begin by extending the results on trees to forests via an observation about disjoint graphs. 
Then we will show two examples of essential flips that appear naturally when considering paths and cycles. 
The observations made for these examples lead to a connection between the characteristic imset polytope and the stable set polytope.

We have previously mentioned the turn pairs defined in \cite{LRS20} which strictly generalize edge reversals for arbitrary skeletons $G$. 
However, for trees these concepts coincide, as we will show in \Cref{subsec: graphical side of turn pairs}. 
We end this section with another example of when \Cref{lem: basis in po vector space} can be applied.

%%%%%%%%%%%%%%%%%%%%%%%%%%%%%%%%%%%%%%%%
%---SUBSECTION: Disjoint graphs
\subsection{Forests}
\label{subsec: disjoint graphs}
Apart from this subsection we have and will only consider skeletons $G$ that are connected, so as to make the results more compact. 
However, for completeness, we will show how to generalize the results to disjoint graphs. %, from the perspective of $\CIM_G$. 
The relevant result is the following:

\begin{proposition}
\label{prop: disjoint graphs}
Let $G=G_1\cup G_2$ be a disjoint union of $G_1=(V_1, E_1)$ and $G_2=(V_2, E_2)$.
Then $\CIM_G$ is affinely equivalent to $\CIM_{G_1}\times\CIM_{G_2}$. 
\end{proposition}

\begin{proof}
Assume we have a DAG $\GG$ with skeleton $G$. 
If $c_\GG(S)=1$, for some set $S\subseteq [p]$ then $G|_S$ is connected.
Indeed, since $c_\GG(S)=1$ we have a node $i$ such that $k\to i$ for all $k\in S\setminus\{i\}$, and thus we have $k-i\in G$ for all $k\in S\setminus\{i\}$. 
Hence for any set $S\subseteq [p]$ we have three cases, $S\subseteq V_1$, $S\subseteq V_2$, or none of the above. 
In the third case we immediately get $c_\GG(S)=0$ for any DAG $\GG$ with skeleton $G$. 
It follows that, up to a permutation of indices, $\CIM_G=\CIM_{G_1}\times\CIM_{G_2}\times\zero$. 
Here $\zero$ denotes the $0$-vector of appropriate length corresponding to the third case above. 
The result follows. 
\end{proof}

Furthermore, it is well-known that if $P$ and $Q$ are two polytopes then any face of $P\times Q$ has the form $F_P\times F_Q$ where $F_P$ and $F_Q$ are faces of $P$ and $Q$ respectively. 
Hence, if we can characterize all faces of $\CIM_G$ where $G$ is connected, we will also characterize all faces of $\CIM_G$ when $G$ is not necessarily connected. 
This is especially true in the case of edges as $\dim(F_P\times F_Q)=\dim(F_P)+\dim(F_Q)$, for any non-empty faces $F_P$ and $F_Q$.
Thus any edge of $P\times Q$ is of the form $e_P\times v_Q$ or $v_P\times e_Q$, where $e_P$ and $e_Q$ are edges of $P$ and $Q$, and $v_P$ and $v_Q$ are vertices of $P$ and $Q$, respectively. 
Thus an edge-walk along $P\times Q$ can be done via two simultaneous edge-walks on $P$ and $Q$. 
The following proposition is direct consequence of \Cref{thm: edges of trees characterization} and \Cref{prop: disjoint graphs}.

\begin{proposition}
Let $G$ be a forest and let $\GG$ and $\HH$ be essential graphs with skeleton $G$. 
Then $\conv(c_\GG, c_\HH)$ is an edge of $\CIM_G$ if and only if there is a unique subtree $T$ of $G$ such that $\{\GG|_T, \HH|_T\}$ is an essential flip and $\GG|_{G\setminus T}=\HH|_{G\setminus T}$.
\end{proposition}

% %---SUBSECTION: Splits and Shifts----
% \subsection{Splits and Shifts}
% \label{subsec: splits and shifts}

%%%%%%%%%%%%%%%%%%%%%%%%%%%%%%%%%%%%%%%%%%
%---SUBSECTION: STABLE SET POLYTOPES ------
  \subsection{Splits, Shifts, and a Connection to Stable Set Polytopes}
\label{subsec: stable set polytopes}
In \Cref{sec: trees} we considered the edges of $\CIM_G$ polytopes directly.
The methods used were specialized towards trees and lack a straightforward generalization.
% Let $I_p$ denote the path of length $p-1$ and $C_p$ denote the cycle of length $p$. 
In the following section we will show a connection between another well studied polytope and $\CIM_{G}$, where $G$ is either $I_p$ or $C_p$.
This will allow us to characterize the edges of the path and the cycle. 
To motivate this, we first mention two examples of essential flips.
% We will however begin to mention two examples of essential flips.
The idea is that we shift and add v-structures along a path in $G$.

%---DEFINITION: Shift----
\begin{definition}
[Shift]
\label{def: shift}
Let $\GG$ and $\HH$ be two DAGs on node set $[p]$.  
We say the pair $\{\GG,\HH\}$ is a \emph{shift} if there exists a path $\pi  = \langle i_0,i_1,\ldots, i_{2m+1}\rangle$ in $\GG$ and $\HH$ such that
\[
c_\HH - \sum_{S\in\SSS_{\pi,\odd}}e_S = c_\GG - \sum_{S\in\SSS_{\pi,\even}}e_S,
\]
where 
\[
\SSS_{\pi,\odd} \coloneqq \big\{\{i_{j-1},i_j,i_{j+1}\} \colon j = 1,3,\ldots, 2m-1\big\} 
\]
and 
\[
\SSS_{\pi,\even} \coloneqq \big\{\{i_{j-1},i_j,i_{j+1}\}\colon j = 2,4,\ldots, 2m\big\}. 
\]
\end{definition}

%--DEFINITION: Spilt----
\begin{definition}
[Split]
\label{def: split}
Let $\GG$ and $\HH$ be two DAGs on node set $[p]$.  
We say the pair $\{\GG,\HH\}$ is a \emph{split} if there exists a path $\pi = \langle i_0,i_1,\ldots,i_{2m}\rangle$ such that 
\[
c_\HH - \sum_{S\in \SSS_{\pi,\odd}}e_S = c_\GG - \sum_{S\in\SSS_{\pi,\even}}e_S,
\]
where 
\[
\SSS_{\pi,\odd} \coloneqq \big\{\{i_{j-1},i_j,i_{j+1}\}\colon j = 1,3,\ldots, 2m-1\big\}
\]
and 
\[
\SSS_{\pi,\even} \coloneqq \big\{\{i_{j-1},i_j,i_{j+1}\big\}\colon j = 2,4,\ldots, 2m-2\}. 
\]
\end{definition}

Both shifts and splits correspond to edges of $\CIM_G$ when $G$ is a tree. 

%--- PROPOSITION: SHIFTS AND SPLITS
\begin{proposition}
\label{prop: splits and shifts}
Let $\GG$ and $\HH$ be two DAGs with the same skeleton $G$, and suppose that $G$ is a tree. 
If $\{\GG,\HH\}$ is a shift or a split then $\conv(c_\GG,c_\HH)$ is an edge of $\CIM_G$.  
\end{proposition}

\begin{proof}
It suffices to show that both shifts and splits constitute essential flips as then the result follows from \Cref{thm: edges of trees characterization}. 
% Using Theorem \ref{thm: subtree condition} Shifts and splits are essential flips since every internal node of the path is of type II or III.
We first show this for shifts. 

By definition of shift we have that $\Delta(\GG, \HH)=\{i_1,\dots,i_{2m}\}$. 
Then all we need to show is that $\GG|_{\spans(\Delta)}$ and $\HH|_{\spans(\Delta)}$ is fully directed, that is every edge $i_k-i_{k+1}$, for $1\leq k\leq 2m-1$ is directed in $\GG$ and $\HH$. 
If $k$ is even we have $c_\GG(\{i_{k-1}, i_k, i_{k+1}\})=1$ and $c_\HH(\{i_k, i_{k+1}, i_{k+2}\})=1$, thus $i_k-i_{k+1}$ is part of a v-structure in both $\GG$ and $\HH$ and thus directed. 
If $k$ is odd, we have $c_\GG(\{i_k, i_{k+1}, i_{k+2}\})=1$ and $c_\HH(\{i_{k-1}, i_k, i_{k+1}\})=1$, and by the same reasoning $i_k-i_{k+1}$ is directed. 
In both cases $i_k-i_{k+1}$ is part of a v-structure, and hence directed, but in different directions.
Thus shifts are essential flips. 

It follows from the definition of splits that $\Delta(\GG, \HH)=\{i_1, i_2,\dots,2m-1\}$. 
Then a split constitutes an essential flip if we can show that every edge $i_k-i_{k+1}$, for $1\leq k\leq 2m-2$, is directed in both $\GG$ and $\HH$. 
Then we can divide into cases depending on whether $k$ is even or odd and proceed exactly as in the case of shifts. 
\end{proof}

While the above proposition only applies to $\CIM_G$ when $G$ is a tree, we believe that shifts and splits can be generalized to edges of $\CIM_G$ for more arbitrary $G$. 
Doing this could be a first step toward extending essential flips to non-trees, and possibly characterizing more edges of $\CIM_G$, regardless of $G$. 

As mentioned, the methods we have used thus far are closely linked to the tree structure of the skeleto, and it is unclear whether they generalize to arbitrary skeletons. 
To this end we will now consider a relation to another well-studied polytope. 
Let $G = ([p],E)$ be an undirected graph. 
A set of nodes $S\subseteq[p]$ is called \emph{stable} (or \emph{independent}) if no nodes of $S$ are joined by an edge.  
Given a set $S\subseteq[p]$, its \emph{incidence vector} $\chi_S\in\{0,1\}^p$ is defined by 
\begin{equation*}
\chi_S(u) = 
\begin{cases}
1 	&	\mbox{if $u\in S$}, \\
0	&	\mbox{if $u\notin S$}.\\
\end{cases}
\end{equation*}
The \emph{stable set polytope} of $G$ is the convex hull
\begin{equation*}
\STAB(G) \coloneqq \conv\left(\chi_S : \mbox{ $S$ is a stable set in $G$}\right)\subseteq\RR^p.
\end{equation*}

%----PROPOSITION: CIM and STAB for Path----
\begin{proposition}
\label{prop: CIM and STAB for path}
For the path $I_p$ on $p\geq 2$ nodes we have $\CIM_{I_p} = \STAB(I_{p-2})$, and for the cycle $C_p$ on $p\geq4$ nodes, 
\[
\CIM_{C_p} = \conv\left(\chi_S : \mbox{ $S$ is a non-empty stable set in $G$}\right).
\]%{\color{red} CORRECT THIS THEOREM}
\end{proposition}

\begin{proof}
Since $I_p$ is the path on $p$ nodes with edges $\{i,i+1\}$ for all $i\in[p-1]$ then for every pair of characteristic imsets $c_\GG$, $c_\HH$ where $\GG$ and $\HH$ are DAGs with skeleton $I_p$, we have that $c_\GG(S) = c_\HH(S)$ for all $S\subseteq [p]$ with $|S|\geq2$ and $S\neq\{i-1,i,i+1\}$ for any $i\in\{2,\ldots,p-1\}$. 
Hence, $\CIM_{I_p}$ is affinely equivalent to its projection into $\RR^{p-2}$, where we associate the standard basis vector $e_{i-1}$ with the standard basis vector $e_{\{i-1,i,i+1\}}$ in $\RR^{2^p-p-1}$.  
We let $\rho(c_\GG)$ denote the image of $c_\GG$ under this projection.  
The natural bijection between MECs with skeleton $I_p$ and stable sets of $I_{p-2}$ described in \cite[Theorem~2.1]{RSU18} then implies that 
\begin{align*}
\{\rho(c_\GG)\in\RR^{p-2} \colon \GG&\text{ is a DAG with skeleton $I_p$} \} =\\ &\left\{\chi_S\in\RR^{p-2} \colon S \text{ a stable set in $I_{p-2}$}\right\}.
\end{align*}
Similarly, since $C_p$ is the cycle on $p$ nodes with edges $\{i,i+1\}$ for $i\in[p-1]$ and $\{1,p\}$, and $p\geq4$, then for every pair of imsets $c_\GG$ and $c_\HH$ where $\GG$ and $\HH$ are DAGs with skeleton $C_p$, we have that $c_\GG(S) = c_\HH(S)$ for all $S\subseteq[p]$ with $|S|\geq 2$ and $S\neq\{i-1,i,i+1\}$ for every $i\in[p]$ (where we treat addition modulo $p$).  
Hence similar to the case of $I_p$, the characteristic imset polytope $\CIM_p$ is affinely equivalent to its projection into $\RR^p$, where we map the standard basis vector $e_{i-1,i,i+1}$ in $\RR^{2^p-p-1}$ to the standard basis vector $e_{i-1}$ in $\RR^p$, for all $i\in[p]$.  
Again, applying the bijection in \cite[Theorem 2.1]{RSU18} between MECs with skeleton $C_p$ and stable sets in $C_p$ implies that 
\[
\{\rho(c_\GG)\in\RR^{p} : c_\GG\in\CIM_{C_p} \} = \{\chi_S\in\RR^{p} : S \text{ a non-empty stable set in $C_{p}$}\},
\]
which completes the proof.  
\end{proof}

Since $\CIM_{I_p}$ and $\CIM_{C_p}$ are, almost, affinely equivalent to stable set polytopes, we can apply a result of Chv\'atal \cite{C75} to give a complete characterization of the edges of $\CIM_{I_p}$ and $\CIM_{C_p}$ in terms of splits, shifts, and v-structure additions. 
%the families of edges of $\CIM_G$ characterized in \Cref{subsec: splits and shifts}.
%---THEOREM: Chvatal Characterization----
\begin{theorem}
\cite[Theorem 6.2]{C75}
\label{thm: chvatal characterization}
Let $G = ([p],E)$ be an undirected graph, let $a,b\in\RR^p$ be two vertices of $\STAB(G)$, and let $A$ and $B$ be their corresponding stable sets.  
Then $\conv(a,b)$ is an edge of $\STAB(G)$ if and only if the subgraph of $G$ induced by $A\setminus B\cup B\setminus A$ is connected. 
\end{theorem}

Hence, to show the desired characterization of the edges of $\CIM_{I_p}$ and $\CIM_{C_p}$, it suffices to characterize the pairs of stable sets in $I_p$ and $C_p$, respectively, for which the subgraph of $I_p$ (or $C_p$) induced by the symmetric difference $A\setminus B \cup B\setminus A$ is connected (i.e., a path or the full cycle).  
\begin{lemma}
\label{lem: stable sets}
Let $A$ and $B$ be stable sets in $I_p$ (or $C_p$). 
Then $\conv(\chi_A,\chi_B)$ is an edge of $\STAB(I_p)$ (or $\STAB(C_p)$) if and only if 
\begin{enumerate}[label=(\arabic*)]
	\item $A\setminus B = \{i,i+2,\ldots,i+2j\}$ for some $i$ and $j$, and 
	\item $B\setminus A = \{i+1,i+3,\ldots,i+2j-1\}$ or $B\setminus A = \{i+1,i+3,\ldots,i+2j+1\}$.
\end{enumerate}
\end{lemma}

\begin{proof}

The proof for $C_p$ works the same as for $I_p$ by taking addition modulo $p$, so we only state it for $I_p$.  
Suppose that $A$ and $B$ are two stable sets in $I_p$ such that $A \setminus B$ and $B\setminus A$ satisfy conditions $(1)$ and $(2)$. 
Then $A \setminus B\cup B\setminus A = \{i,i+1,i+2,\ldots,i+2j+1\}$ for some $i$ and $j$, and hence the induced subgraph on this set is connected.  
It follows from \Cref{thm: chvatal characterization} that $\conv(\chi_A,\chi_B)$ is an edge of $\STAB(I_p)$.

Conversely, if $A$ and $B$ are stable sets of $I_p$ such that $\conv(\chi_A,\chi_B)$ is an edge of $\STAB(I_p)$, then by \Cref{thm: chvatal characterization}, we know that the subgraph of $I_p$ induced by $A\setminus B \cup B\setminus A$ is connected, and hence must be a subpath of $I_p$, say $\langle i,i+1,\ldots, i+t\rangle$.  
Since $A$ and $B$ are stable, then two neighbors on this path cannot belong to the same set $A$ or $B$.  
It follows that $A\setminus B$ and $B\setminus A$ must satisfy conditions $(1)$ and $(2)$, completing the proof. 
\end{proof}

As a consequence of \Cref{lem: stable sets}, we can characterize all edges of $\CIM_{I_p}$ and $\CIM_{C_p}$.
%---THEOREM: Paths and Cycles---
\begin{theorem}
\label{thm: paths and cycles}
%Let $I_p$ and $C_p$ denote the path and cycle on $p$ nodes, respectively, and
Let $c_\GG$ and $c_\HH$ be two characteristic imsets for DAGs both having skeleton the path $I_p$ (or the cycle $C_p$).  
Then $\conv(c_\GG,c_\HH)$ is an edge of $\CIM_{I_p}$, for $p\geq 2$ if and only if $\{\GG,\HH\}$ is a v-structure addition, shift, or split.
Moreover $\conv(c_\GG, c_\HH)$ is an edge of $\CIM_{C_p}$ if and only if $\{\GG, \HH\}$ is a v-structure addition, shift, split, or both $\GG$ and $\HH$ contain exactly one v-structure. 
\end{theorem}

\begin{proof}
By \Cref{prop: CIM and STAB for path}, we know that $\CIM_{I_p} = \STAB(I_{p-2})\subseteq\RR^{p-2}$, for $p\geq 2$, and $\CIM_{C_p} = \{\chi_S\in\RR^{p} : S \text{ a non-empty stable set in $C_{p}$}\}\subset\RR^p$, for $p\geq 4$, where we have identified the standard basis vector $e_{\{i-1,i,i+1\}}\in\RR^{2^p-p-1}$ with the standard basis vector $e_{i-1}\in\RR^{p-2}$ and $e_{i-1}\in\RR^p$, respectively.  
Here, we again consider addition modulo $p$ in the case of the cycle $C_p$.  
Since the hyperplane $\sum_{i=1}^px_i = 1$ is facet-defining for $\CIM_{C_p}\subset\RR^p$, the edges of $\CIM_{C_p}$ are precisely the edges of $\STAB(C_p)$ and the standard $p-1$-simplex in $\RR^p$, minus those between the origin and the standard basis vectors in $\RR^p$.  
As the same proof works for both the cycle and the path, apart from the previous sentence, in the following we state it only for the path $I_p$.  

By \Cref{lem: stable sets}, we get that for two stable sets $A$ and $B$ in $I_p$, $\conv(\chi_A,\chi_B)$ is an edge of $\STAB(I_{p-2})$ if and only 
\begin{enumerate}[label=(\arabic*)]
	\item $A\setminus B = \{i,i+2,\ldots,i+2j\}$ for some $i$ and $j$, and 
	\item $B\setminus A = \{i+1,i+3,\ldots,i+2j-1\}$ or $B\setminus A = \{i+1,i+3,\ldots,i+2j+1\}$.
\end{enumerate}
As $A\setminus B \cup B \setminus A = \{i,i+1,\ldots, i+m\}$, where $m = 2j$ or $m = 2j+1$ accordingly, it follows that $\conv(\chi_A,\chi_B)$ is an edge of $\STAB(I_{p-2})$ if and only if $\pi = \langle i,i+1,i+2,\ldots,i+m\rangle$ is a path in $I_{p-2}$, and 
\[
\chi_A - \sum_{k\in A\setminus B}e_k = \chi_B - \sum_{k\in B\setminus A}e_k.
\] 
Following the correspondence between vertices of $\CIM_{I_p}$ and vertices of $\STAB(I_{p-2})$ established above, $\conv(c_\GG,c_\HH)$ is an edge of $\CIM_{I_p}$ if and only if $\pi = \langle i,i+1,i+2,i+3,\ldots,i+m+2\rangle$ is a path in $I_p$ and 
\[
c_\GG - \sum_{\substack{\{k-1,k,k+1\}: \\ i+\ell=k\leq i+m+1\\ \text{ and $\ell$ is odd}} }e_{\{k-1,k,k+1\}} = c_\HH - \sum_{\substack{\{k-1,k,k+1\}: \\ i+\ell=k\leq i+m+1, \\ \ell\geq 2 \text{ and $\ell$ is even} }}e_{\{k-1,k,k+1\}}.
\] 
When $m = 0$, then $\{\GG,\HH\}$ is a v-structure addition, when $m>0$ and odd $\{\GG,\HH\}$ is a shift, and when $m>0$ and even $\{\GG,\HH\}$ is a split.
Hence, $\conv(c_\GG,c_\HH)$ is an edge of $\CIM_{I_p}$ if and only if $\{\GG,\HH\}$ is a shift, split, or v-structure addition. 
\end{proof}

We remark that $I_p^\circ=I_{p-2}$, hence \Cref{prop: CIM and STAB for path} says that $\STAB(I_p^\circ)=\CIM_{I_p}$.
More generally, the stable set polytope of a tree can always be realized as a face of $\CIM_G$ for an appropriately chosen graph $G$.

%--- PROPOSITION: STAB AS FACE OF CIM ----
\begin{proposition}
\label{prop: stab as face of cim}
Let $G$ be a tree. 
Then $\STAB(G^\circ)$ is unimodularly equivalent to a face of $\CIM_G$. 
\end{proposition}

\begin{proof}
Let $i$ be an internal node of $G$ and let $\GG_i$ be a DAG with skeleton $G$, $\pa_{\GG_i}=\ne_G(i)$, and no other v-structures.
Take $\HH$ to be a DAG without any v-structure. 
Such DAGs exists since $G$ is a tree. 
Similar to the proof of \Cref{prop: edge of tree}, there exists affine functions $w_i$ that only depend on $N_i$ such that $w_i^Tc_{\GG_i}=w_i^Tc_{\HH}=0$ and $w_i^Tc_\DD=-1$ if $c_\DD|_{N_i}\neq c_{\GG_i}|_{N_i}$ and $c_\DD|_{N_i}\neq c_{\HH}|_{N_i}$. 
Notice that since $i$ was taken to be an internal node we have $c_{\GG_i}\neq c_\HH$.
Then we let $W=\sum w_i$ where we sum over all internal nodes of $G$.
For convenience we also let $s_i$ be affine functions similar to $w_i$ but $s_i^Tc_{\GG_i}=1$ and $s_i^Tc_{\DD}= 0$ for all $\DD$ such that $c_\DD|_{N_i}\neq c_{\GG_i}|_{N_i}$. 
Then we can define $S(\DD)=\{i\in G^\circ \colon s_i^Tc_\DD=1\}$.

Assume $\DD$ maximizes $W^Tc_\DD$, then we claim that $S(\DD)$ is a stable set in $\DD$. 
It is clear that $W^Tc_\DD\leq 0$ as each $w_i$ only assumes non-positive values and as $W^Tc_\HH=0$ it follows that we must have $W^Tc_\DD= 0$. 
Thus by the definition of $W$ we must have $w_i^Tc_\DD=0$ for every internal node $i$. 
As every internal node has degree at least $2$, we must have at least one v-structure at every $i\in S(\DD)$. 
Then by \Cref{lem: essential arrows at a node} we must have that $\DD|_{\cl_G(i)}$ is fully directed and as $s_i^Tc_\DD=1$ we must have that every node in $\ne_G(i)$ is an essential parent of $i$ in $\DD$. 
Then if we would have two neighboring nodes $i, j\in S(\DD)$ then they would each be a parent of each other, which cannot happen. 
Hence $S(\DD)$ is indeed a stable set in $G^\circ$. 

Conversely, given a stable set $T$ in $G$, we want to construct a DAG $\DD$ with $S(\DD)=T$. 
For each $i\in T$ we let $\pa_{\DD}=\ne_G(i)$. 
This is possible since $T$ is a stable set.
% , as $T$ was stable this is possible. 
Then we can direct every part in $G|_{[p]\setminus T}$ without v-structures.
This is possible as $G|_{[p]\setminus T}$ is a forest. 
Left to check is that $\DD$ indeed maximizes $W$, which follows since $w_i^Tc_\DD=0$.

Hence $S$ is a bijection between essential graphs $\DD$ maximizing $W$ and stable sets of $G^\circ$.
Let $F$ denote the face of $\CIM_G$ maximizing $W$. 
Then it is direct that 
\[
\chi_{S(\DD)}=\left(\chi_{S(\DD)}(1),\dots,\chi_{S(\DD)}(m)\right)^T=\left(s_1^Tc_\DD, \dots, s_m^Tc_\DD\right)^T 
\]
for any essential graph $\DD\in F$, where $m$ is the number of internal nodes of $G$. 
Thus we have a bijective map between the vertices of $F$ and the vertices of $\STAB(G^\circ)$ given by $\rho\colon c_\DD\mapsto \left(s_1^Tc_\DD, \dots, s_m^Tc_\DD\right)^T$. 
As $\rho$ is affine in each coordinate we can extend it to an affine map from the affine subspace containing $F$ to the subspace containing $\STAB(G^\circ)$. 
% this gives us an affine bijective map as $c_\DD\mapsto \left(s_1^Tc_\DD, \dots, s_m^Tc_\DD\right)^T$ from the vertices of $F$ to the vertices of $\STAB(G^\circ)$.  
What is left to show is that this affine map is invertible.
It can however be checked that 
\[
\chi_T\mapsto c_\HH+\sum_{i\in G^\circ} \chi_T(i)c_{\GG_i}|_{N_i}
\]
is the inverse of $\rho$, and it is also affine. 
Hence we have a bijective affine correspondence between $F$ and $\STAB(G^\circ)$.
The map $\rho$ is also unimodular since the lattice of the affine subspace of $F$ is spanned by $\left\{c_{\GG_i}|_{N_i}\right\}_{i\in G^\circ}=\left\{c_{\GG_i}-c_\HH\right\}_{i\in G^\circ}$.
\end{proof}

In the case of $I_p$ and $C_p$ it can be checked that all essential graphs maximize the cost function $W$, as constructed in the proof above. 
This is due to the fact that every internal vertex has degree exactly two, and thus $|N_i|=1$. 
Hence the face of $\CIM_{I_p}$ isomorphic to $\STAB(I_{p}^\circ)$ is the polytope itself.
This isomorphism fails for the cycle, as the empty set (which is stable) is not mapped to a characteristic imset of any DAG. 
This is also why \Cref{prop: CIM and STAB for path} required non-empty stable sets for the cycle.

%%%%%%%%%%%%%%%%%%%%%%%%%%%%%%%%%%%%%%%%
%---SUBSECTION: GRAPHICAL SIDE OF TURN PAIRS
\subsection{The Graphical Side of Turn Pairs}
\label{subsec: graphical side of turn pairs}
It was shown in \cite{LRS20} that turn pairs strictly generalize reversing an edge of a DAG.
However, for trees, this is not true. 
That is, the converse of \cite[Proposition 3.2]{LRS20} holds when the underlying skeleton $G$ is a tree.  
Let us first recall the definition of a turn pair.

%---DEFINITION: Turn pair------
\begin{definition}
[Turn pair]
\label{def: turn pair}
\cite{LRS20}
Let $\GG$ and $\HH$ be two DAGs on node set $[p]$ and with skeleton $G$.  
Suppose there exist $i$,  $j$, $S_i\subseteq[p]\backslash \{i,j\}$ and $S_j\subseteq[p]\backslash \{i,j\}$ such that 
\begin{enumerate}[label=(\arabic*)]
\item{$c_\GG(\{i,j\})=1$;}
\item{$c_\GG(S\cup \{i\})=1$ for all $S\subseteq S_i$ with $|S|\geq 1$;}
\item{$c_\GG(S\cup \{j\})=1$ for all $S\subseteq S_j$ with $|S|\geq 1$;}
\item{either $S_i\not\subseteq\ne_G(j)$ or $S_j\not\subseteq\ne_G(i)$.}
\end{enumerate}
Then we say that $\{\GG, \HH\}$ is a \emph{turn pair} with respect to $(i,j,S_i, S_j)$ if 
\[
c_{\HH} = c_\GG +\sum_{S\in \SSS^+}e_S - \sum_{S\in \SSS^-}e_S
\]
where  $\SSS^+\coloneqq \{T\cup \{i,j\}\colon T\subseteq S_i,\  T \not\subseteq\ne_G(j)\}$ and $\SSS^-\coloneqq \{T\cup \{i,j\}\colon T\subseteq S_j,\ T \not\subseteq\ne_G(i)\}$.
\end{definition}

Then we get the following graphical characterization. 

%---PROPOSITION: Characterizing Turn Pairs in Trees----
\begin{proposition}
\label{prop: turn pair trees}
Assume $\{\GG,\HH\}$ is a turn pair where $\GG$ and $\HH$ both have skeleton $G$, a tree.
Then there exists a DAG $\DD$ and nodes $i'$ and $j'$ such that
\begin{enumerate}[label=(\arabic*)]
\item{$\GG$ and $\DD$ are Markov equivalent,}
\item{$i'\to j'\in\DD$,}
\item{$\DD_{i'\gets j'}$ is a DAG, and}
\item{$\HH$ and $\DD_{i'\gets j'}$ are Markov equivalent,}
\end{enumerate}
\end{proposition}
% We will use the same notation as in \cite{LRS20}.

\begin{proof}
% This is the combined result of \Cref{lem: addition trees}, \Cref{lem: bud trees} and \Cref{lem: flip trees}.
Notice that condition (3) is immediate as $G$ is a tree. 
By definition of turn pair $\{\GG, \HH\}$ is with respect to some $(i, j,  S_i, S_j)$ such that either $S_i\neq\emptyset$ or $S_j\neq\emptyset$.
Hence we can, by symmetry, assume that $|S_i|\geq |S_j|$ and that $|S_i|\geq 1$. 
Then we have three cases, $|S_i|=1$ and $|S_j|=0$, $|S_i|=1$ and $|S_j|=1$, and $|S_i|\geq2$. 

\textbf{Case I, $|S_i|=1$ and $|S_j|=0$:} It follows by definition of turn pair that $|\SSS^+|=1$ and thus $\{\GG, \HH\}$ is an addition with respect to $\{i,j\}\cup S_i$. 

Since $\GG$ and $\HH$ have the same skeleton it follows by \Cref{prop: addition characterization} that they differ by a single v-structure.
That is, the induced graph of $\HH$ on $\{i,j\}\cup S_i$ is a v-structure, $\alpha\to \beta\leftarrow \gamma$.
As $c_\GG(\{i,j\}\cup S_i)=0$ we have that $\alpha\leftarrow \beta\in \GG$ or $\beta\to \gamma\in \GG$. 
We can assume by symmetry that $\beta\to\gamma \in\GG$.
Assume we have a node $s\in \pa_{\GG}(\gamma)\setminus\{\beta\}$.
Since $G$ is a tree it follows that $s\neq \alpha$.
Then $c_\GG(\{s,\beta,\gamma\})=1$ and, as $\{s,\beta,\alpha\}\neq \{i,j\}\cup S_i$, we get $c_\HH(\{s,\beta,\gamma\})=1$. 
Again since $G$ is a tree we get that $\beta\to \gamma\leftarrow s\in \HH$, a contradiction.
Thus $\pa_\GG(\gamma)=\{\beta\}$.

Assume we have a node $s\in\pa_\HH(\beta)\setminus\{\alpha, \gamma\}$. 
Then we have $s\to \beta\leftarrow \gamma\in\HH$ and thus $c_\HH(\{s,\beta,\gamma\})=1$.
However, as $\{s,\beta,\gamma\}\neq \{i,j\}\cup S_i$ we get $s\to \beta\leftarrow \gamma\in\GG$, a contradiction.
Thus $\pa_\HH(\beta)=\{\alpha, \gamma\}$.

Now let us construct $\DD$.
Let $G'$ denote the graph identical to $G$ with the exception that it does not contain the edge $\beta - \gamma$.
Let $G_\beta$ be the connected part of $G'$ containing $\beta$ and let $G_\gamma$ be the rest of $G'$.
Define $\DD$ as the orientation of $G$ where we direct $G_\beta$ as in $\HH$, direct $G_\gamma$ according to $\GG$ and direct $\beta-\gamma$ as in $\beta\to \gamma$.

First we want to show that $\DD$ is Markov equivalent to $\GG$. 
As $\GG$ and $\HH$ only differ by a single v-structure and $\DD$ is everywhere directed as in either $\GG$ or in $\HH$ the only place where $\DD$ can differ by a v-structure from $\GG$ is around $\beta$. 
However, since $\pa_\DD(\beta)=\pa_\HH(\beta)\setminus\{\gamma\} =\{\alpha\}$ neither has any v-structures around $\beta$.
Thus, they must be Markov equivalent.

It follows that $\DD_{\beta\leftarrow \gamma}$ will have the v-structure $\alpha \to \beta\leftarrow \gamma$ but no other v-structures are created or destroyed when we reverse $\beta\to\gamma$, as $\pa_\DD(\beta)=\{\alpha\}$. 
Hence $\DD_{\beta\gets \gamma}$ is Markov equivalent to $\HH$.
In conclusion, we can choose $\DD$ with $i'=\beta$ and $j'=\gamma$. 

\textbf{Case II, $S_i=\{\alpha_i\}$ and $S_j=\{\alpha_j\}$:} 
Since $G$ is a tree we must have the v-structure $\alpha_i\to i\leftarrow j$ in $\HH$ and the v-structure $i\to j\leftarrow \alpha_j$ in $\GG$.  
Then we can define $G_i$ and $G_j$ similar to how we defined $G_\beta$ and $G_\gamma$ in {\bf Case I} above. 
The rest is completely analogous and we get $i'=i$ and $j'=j$.  

\textbf{Case III, $|S_i|\geq 2$:} 
As $c_\GG(S_i\cup \{i\})=1$ we have at least one v-structure $\alpha\to i\leftarrow \gamma$ for some $\alpha, \gamma\in S_i$ in $\GG$. 
Thus we must have $\beta\to i$ in $\GG$ for all $\beta\in S_i$.
The same is true for $\HH$. 
Then we have $c_\GG(S_i\cup \{i,j\})=0$ as we cannot have $S_i\subseteq \ne_G(j)$ as $G$ is a tree. 
Hence $i\to j$ is essential in $\GG$, and similar reasoning gives us $i\leftarrow j$ is essential in $\HH$. 
Then we can choose $\DD=\GG$. 
From a reasoning similar to {\bf Case I} we get that $\pa_\GG(j)=S_j$, and that $\GG_{i\leftarrow j}$ is Markov equivalent to $\HH$ is straightforward.
Thus we get $i'=i$ and $j'=j$ in this case as well. 
\end{proof}

Thus turn pairs exactly correspond to reversing an edge in some element of the MEC. 

%%%%%%%%%%%%%%%%%%%%%%%%%%%%%%%%%%%%%%%%%%%%%%%
%---SUBSECTION: ALMOST COMPLETE GRAPHS
\subsection{Almost Complete Graphs}
\label{subsec: almost complete graphs}

When we proved \Cref{prop: star like graph simplex} we utilized the fact that $G$ induced a partial order on the MECs with skeleton $G$. 
Here we present another case where this happens as well. 
% We can prove this in a similar manner, using \Cref{lem: basis in po vector space}.
%---LEMMA: Dense graph simplex ------
\begin{proposition}
\label{prop: case 4 simplex}
Let $G=([p], E)$ be the complete graph missing only the edge $1-2$, that is $K_p\setminus \{1,2\}$. 
%That is $\{i,j\}\in E$ if and only if $\{i,j\}\neq\{1,2\}$.
Then $\CIM_G$ is a simplex of dimension $2^{p-2}-1$.
\end{proposition}

\begin{proof}
By the structure of $G$ the only possible v-structures are of the form $1\rightarrow i \leftarrow 2$ for some $2<i\leq p$. 
We wish to show that for any subset $S\subseteq \{3,\dots p\}$ we have a DAG $\DD_S$ such that $1\rightarrow i \leftarrow 2\in \DD_S$ for all $i\in S$ but for no $i\in \{3,\dots p\}\setminus S$. 
As every MEC is characterized by the v-structures we will then have a representative from each MEC.

As for the construction of $\DD_S$, for every $i\in S$ we direct the edges $1-i-2$ as $1\rightarrow i \leftarrow 2$ and for every $i\in \{3,\dots p\}\setminus S$ we direct the edges $1-i-2$ as $1\leftarrow i \rightarrow 2$. 
Notice that we have thus far not created any cycles. 
Hence, the currently directed edges induce a partial order on $[p]$. 
Extending this order to a total order and directing the remaining edges of $\DD_S$ according to this order gives us a DAG.
Notice that no new v-structures can have been created in this last step, because all possible triples that could have been v-structures were already directed.
Then we can calculate the characteristic imset of $\DD_S$ as 
\[
c_{\DD_S}=c_{\DD_\emptyset} + \sum_{T\in \SSS}e_T
\]
where $\SSS=\{U\cup \{1,2\} \colon U\subseteq \{3,\dots, p\}, U\cap S\neq \emptyset\}$. Equivalently we have 
$\mathbf{1}-e_{\{1,2\}}-c_{\DD_S}=\sum_{T\in\SSS'} e_T$
where $\SSS'=\{ U\cup \{1,2\} \colon U\subseteq [p]\setminus S\}$ where $\mathbf{1}$ is the constant $1$-vector. 
By \Cref{lem: basis in po vector space} the right-hand-side constitutes a basis and since translations preserve this property the result follows. 
\end{proof}

We have now seen several cases where posets can be of use in describing the geometry of the $\CIM_G$ polytope. 
Whether there is a more general connection or not is left to future work.

%%%%%%%%%%%%%%%%%%%%%%%%%%%%%%%%%%%%%%%%%%%%%%%%%%%%%
%---SECTION: Applications ------
\section{Applications}
\label{sec: applications}
A fundamental problem in modern data science and artificial intelligence is the problem of \emph{causal inference} \cite{pearl2009causality}, in which one is interested in estimating the cause-effect relations between jointly distributed random variables.  
This task is typically broken into two well-studied subproblems: (1) the task of inferring the strength and nature of the causal effect of one variable on another, and (2) the task of estimating which variables have direct causal effects on another.  
The latter of the two problems is referred to as the problem of \emph{causal discovery}.  
In its most basic form, we assume that we have a random sample $\mathbb{D}$ drawn from the joint distribution $\mathbb{P}$ of the random variables $(X_1,\ldots, X_p)$, and we would like to infer a DAG $\GG = ([p], E)$ in which the arrows $i\rightarrow j$ correspond to the direct cause-effect relations in the system; i.e., $\GG$ contains the arrow $i\rightarrow j$ if and only if $i$ is a direct cause of $j$. 
As the data $\mathbb{D}$ is drawn from a probability distribution, and correlation does not imply causation, one must impose some assumption that associates the data-generating distribution $\mathbb{P}$ to the underlying (unknown) causal structure $\GG$ we wish to infer.  
We say a distribution $\mathbb{P}$ over $(X_1,\ldots, X_p)$ is \emph{Markov} to a DAG $\GG = ([p], E)$ if there is a linear extension (topological ordering) $\pi = \pi_1\cdots\pi_p$ of $\GG$ such that for all $i\in[p]$, $\mathbb{P}$ entails the CI relation 
\[
X_{\pi_i} \independent X_{\{\pi_1,\ldots, \pi_{i-1}\}\setminus \pa_{\GG}(\pi_i)}\mid X_{\pa_{\GG}(\pi_i)}.
\]
The above relations capture rudimentary causal information encoded in the distribution $\mathbb{P}$; namely, that each variable (event) is independent of all preceding variables (events) given its direct causes. 
Given a DAG $\GG$, its associated \emph{DAG model} is the collection of distributions
\[
\mathcal{M}(\GG) = \{\mathbb{P} : \mathbb{P} \mbox{ is Markov to } \GG\}.
\]
Since $\mathcal{M}(\GG) = \mathcal{M}(\HH)$ if and only if $\GG$ and $\HH$ are Markov equivalent \cite{lauritzen1996, pearl2009causality}, then, given only the data $\mathbb{D}$, the best we can hope for is to recover $\GG$ up to Markov equivalence.  
Hence, the basic problem of causal discovery is to learn the essential graph of the causal system of the data-generating distribution based on the random sample $\mathbb{D}$. 

Proposed causal discovery algorithms typically come in one of three forms: \emph{constraint-based algorithms}, such as the PC-algorithm \cite{SG91, T19}, that recover an essential graph via statistical tests for conditional independence, \emph{(greedy) score-based algorithms}, such as the Greedy Equivalence Search (GES) \cite{C02}, that assign a score to each DAG (or essential graph) based on the data and then search for the optimal scoring DAG, and \emph{hybrid algorithms}, such as the Max-Min Hill Climbing Algorithm (MMHC) \cite{TBA06}, that use a mixture of CI-testing and optimization. 
Most classic causal discovery algorithms rely solely on the combinatorics of DAGs.
However, recent advancements have introduced discrete-geometric methods with promising results. 
These include the hybrid algorithm \emph{GreedySP} \cite{SUW20} and score-based linear optimization approaches that aim to solve the LP associated to the polytope $\CIM_p$ \cite{SHL10}.  
Most recently, \cite{LRS20} introduced a hybrid algorithm, \emph{skeletal Greedy CIM}, that first estimates the skeleton $G$ of the essential graph using CI-tests and then uses the edges of $\CIM_G$ identified in \cite{LRS20} to search over the polytope for the optimal essential graph.  

Skeletal greedy CIM performed quite well on simulated data, despite using CI-tests (which are prone to error propagation) to learn the skeleton $G$ and only the subset of the edges of $\CIM_G$ corresponding to turn pairs (see Definition~\ref{def: turn pair}). 
In the case that the unknown causal system is a directed tree (often called a \emph{polytree}), the results of \Cref{thm: edges of trees characterization} can be used to overcome the latter of these two limiting factors. 
The CI-tests can similarly be avoided by taking advantage of the fact that we wish to learn only a polytree.  
In place of CI-tests, we can instead learn the skeleton $T = ([p], E)$ of the polytree by learning a \emph{minimum weight spanning tree (MWST)} of a complete graph $K_p$ where the weight assigned to each edge $i - j$  is the negation of the \emph{mutual information} of $X_i$ and $X_j$:
\[
I(X_i; X_j) = \int\int f(x_i,x_j)\log\left(\frac{f(x_i,x_j)}{f(x_i)f(x_j)}\right)dx_idx_j.
\]
Given the inferred skeleton $\hat{T}$, via the edge characterization of $\CIM_T$ for $T$ presented in \Cref{thm: edges of trees characterization}, we can estimate the causal structure $\GG$ by performing an edge walk along $\CIM_{\hat{T}}$, where at each step we walk to the neighbor of the current essential graph that optimally increases the \emph{Bayesian Information Criterion (BIC)} score of the model:
\[
\BIC(\GG,\mathbb{D}) = \log f\left(\mathbb{D}|\hat{\theta}, \GG^h\right)-\frac{d}2\log(|\mathbb{D}|).
\]
Here, $\hat{\theta}$ is the maximum likelihood estimate for the model parameters, $d$ denotes the number of free parameters, and $\GG^h$ denotes the hypothesis that $\mathbb{D}$ is a random sample from a distribution entailing the CI statements encoded by the given DAG $\GG$.
Since $\BIC(\GG, \mathbb{D})$ is known to be an (affine) linear function over $\CIM_p$ \cite{SHL10}, the algorithm terminates once no neighboring vertex of the current characteristic imset increases in $\BIC$. 
We call this algorithm for learning polytrees, presented in \Cref{alg: EFT}, the \emph{Essential Flip Tree Search (EFT)}. 

%---ALGORITHM: EFT----
\begin{algorithm}
  \caption{Essential Flip Tree Search (EFT)}
  \label{alg: EFT}
  \raggedright
  \hspace*{\algorithmicindent} \textbf{Input:} a random sample $\mathbb{D}$ from $(X_1,\ldots,X_p)$.\\
%   $G = ([p],E)$ an undirected tree, $\mathbb{D}$ iid samples from a distribution on $(X_1,\ldots,X_p)$, $S_\mathbb{D}$ a decomposable scoring criterion.\\
  \hspace*{\algorithmicindent} \textbf{Output:} $\GG = ([p],E)$ a polytree
  \begin{algorithmic}[1]
    \State $K_p \gets$ {a weighted complete graph on $[p]$ where $i - j$ is assigned the weight $-I(X_i;X_j)$.}
    \State $G = ([p],E) \gets \mbox{a minimum weight spanning tree of $K_p$}$
    \State $\GG = ([p],A) \gets \mbox{a random polytree with skeleton $G$}$
    \State $\HH \gets \GG$
    \State $M \gets$ all subtrees of $G$
    \While{TRUE}
        \State $\GG_{\overleftarrow{s}}\gets \GG$ with edges in $s\in M$ reversed if $\{\GG_{\overleftarrow{s}}, \GG\}$ constitutes an essential flip. 
        \State $\GG_M \gets \{\GG_{\overleftarrow{s}} : s\in M\}$
        \State $\HH \gets \argmax_{\DD\in \GG_M}\BIC(\GG_{\overleftarrow{s}},\mathbb{D})$
        \If{$\BIC(\GG,\mathbb{D}) \geq \BIC(\HH,\mathbb{D})$}
            \State $\mathbf{break}$
        \EndIf
        \State $\GG\gets \HH$
    \EndWhile\\
    \Return $\GG$
  \end{algorithmic}
\end{algorithm}

An implementation of EFT is available at \cite{github}.  
In the remainder of this section, we consider the aspects of causal discovery algorithms arising from our discrete-geometric understanding of $\CIM_p$ and its faces. 
In subsection~\ref{subsec: consistency}, we consider the asymptotic consistency of causal discovery algorithms that perform edge walks along faces of the $\CIM_p$, including the hybrid algorithms EFT and skeletal Greedy CIM, as well as the purely score-based algorithm \emph{Greedy~CIM} presented in \cite{LRS20}. 
In subsection~\ref{subsec: computations}, we analyze the performance of EFT on simulated and real data, comparing it with classic approaches for learning polytrees such as the hybrid algorithm of Rebane and Pearl \cite{RP87}, which we will call the \emph{RP-algorithm}.  
We observe that EFT outperforms such classic polytree learning algorithms, and can even outperform general causal discovery algorithms such as GES and GreedySP on real data.

%---SUBSECTION: Consistency ------
\subsection{A Sufficient Criterion for Consistency}
\label{subsec: consistency}
It was shown in \cite{C02} that GES is asymptotically consistent when the data used to identify a BIC-optimal MEC is drawn from a distribution faithful to a DAG in the MEC; that is, GES will return the correct MEC as the sample-size goes to infinity.  
In short, we say that GES is \emph{consistent under the faithfulness assumption}.
To show that GES is consistent under faithfulness Chickering \cite{C02} defined the notion of \emph{local consistency}, which informally says that adding in missing edges that should be present and removing the existing edges that should be missing increases the score. 
Here we will show a new version of consistency for reversing edges similar to that of local consistency (\Cref{prop: local turn consitency}). 
With this we get consistency of skeletal greedy CIM, as defined in \cite{LRS20}, and EFT (i.e., \Cref{alg: EFT}). 
These results will all assume \emph{consistency} of the scoring criterion. 

\begin{definition}[Consistent scoring Criterion]
Let $\mathbb{D}$ be a random sample of size $m$ from some distribution
$\mathbb{P}$. A scoring criterion $S$ is \emph{consistent} if in the limit as $m$ grows large, the following two
properties hold:
\begin{enumerate}
    \item If $\mathcal{M}(\HH)$ contains $\mathbb{P}$ and $\mathcal{M}(\GG)$ does not contain $\mathbb{P}$ then $S(\HH, \mathbb{D}) > S(\GG, \mathbb{D})$.
    \item If $\HH$ and $\GG$ both contain $\mathbb{P}$ and $\GG$ contains fewer parameters than $\HH$ then $S(\GG, \mathbb{D})> S(\HH, \mathbb{D})$.
\end{enumerate} 
\end{definition}

% \label{prop: bic consistent}
We also say that a scoring criterion $S( -,\mathbb{D})$ is \emph{score equivalent} if $S(\GG, \mathbb{D})=S(\HH, \mathbb{D})$ whenever $\GG$ and $\HH$ are Markov equivalent. 
In \cite{H88} the author showed that $\BIC$ is a consistent and score equivalent scoring criterion for a large class of models which include DAG models. 
Moreover, in \cite{C02, SHL10} the authors remark that $\BIC$ is \emph{decomposable}, a condition making it a suitable scoring criterion for Skeletal Greedy CIM and Greedy CIM. 
Hence, the following proposition is indeed applies to $\BIC$. 

\begin{proposition}
\label{prop: local turn consitency}
Let $\GG$ be a DAG with $i\to j\in\GG$, $\mathbb{P}$ a (positive) distribution faithful to $\HH$, $\mathbb{D}$ a random sample from $\mathbb{P}$, $S(-,\mathbb{D})$ a score equivalent, decomposable and consistent scoring criterion, and assume that $i$ and $j$ are neighbours in $\HH$.  
Assume that $\GG_{i\leftarrow j}$ is a DAG not Markov equivalent to $\GG$. 
Then if $\HH$ has all v-structures $k\to i\leftarrow j$ for $k\in\pa_\GG(i)\setminus \pa_\GG(j)$ and none of the v-structures $i\to j\leftarrow k$ for $k\in\pa_\GG(j)\setminus(\pa_\GG(i)\cup\{i\})$ 
we have $S(\GG, \mathbb{D})<S(\GG_{i\leftarrow j}, \mathbb{D})$. 
\end{proposition}

\begin{proof}
The proof is similar in nature as the proof of local consistency, see \cite[Lemma 7]{C02}. 
Since $S(-,\mathbb{D})$ is decomposable the difference in score solely depends on the structure of $\GG|_{\pa(i)\cup \pa(j)\cup \{i,j\}}$. % and the structure therein. 
More specifically, $\{\GG, \GG_{i\leftarrow j}\}$ is a turn pair, with $S_i=\pa_\GG(i)$ and $S_j=\pa_\GG(j)\setminus\{i\}$, and thus by the equality in \Cref{def: turn pair},
% as can be seen in computations done in \cite{LRS20}, 
the difference in score only depends on the edges between $j$ and vertices in $\pa_\GG(i)$, and edges between $i$ and vertices in $\pa_\GG(j)$. 
Hence we can assume that we have an edge between every pair of vertices that is not of the described form in $\GG$. 
Let $\tilde{\HH}$ be a DAG with all edges of $\HH$ and for each pair of nodes $\alpha, \beta$, such that we do not have $\alpha= i$ and $\beta\in\pa_\HH(j)$ or $\alpha= j$ and $\beta\in\pa_\HH(i)$, or vice versa, we add in the edge $\alpha -\beta$ with an acyclic orientation. 
This is possible as $\HH$ is a DAG. 
Notice that as $\mathbb{P}$ is faithful to $\HH$ it must be Markov to $\tilde\HH$, as $\tilde\HH$ imposes less restrictions on $\mathbb{P}$. % {\color{red} maybe add in another sentence here, or find nice reference}.
The only v-structures possible in $\tilde\HH$ are of the form $k\to i\leftarrow j$ for $k\in\pa_\GG(i)\setminus \pa_\GG(j)$ or $i\to j\leftarrow k$ for $k\in\pa_\GG(j)\setminus(\pa_\GG(i)\cup\{i\})$. 
Hence, by assumption, $\GG_{i\leftarrow j}$ must be Markov equivalent to $\tilde\HH$ as, also by assumption, both have exactly the same v-structures.
Therefore, $\mathbb{P}$ is Markov to $\GG_{i\leftarrow j}$. 
The result then follows from the consistency of $S(-,\mathbb{D})$ if we can show that $\mathbb{P}$ is not Markov to $\GG$. 

We begin with the case that $\pa_{\GG}(i)\setminus\pa_{\GG}(j)\neq \emptyset$.
Let $k\in \pa_{\GG}(i)\setminus\pa_{\GG}(j)$ and note that $\GG$ encodes $X_k\independent X_j\mid X_{\an_{\GG}(j)\setminus\{k\}}$.
We wish to show that this statement is not encoded by $\HH$ and hence is not true in $\mathbb{P}$. 
However, by assumption we have $i\in \an_{\GG}(j)\setminus\{k\}$ and  $k\to i\leftarrow j$ is a v-structure in $\HH$. 
Then $X_k\independent X_j\mid X_{\an_{\GG}(j)\setminus\{k\}}$ cannot be encoded by $\HH$.
In the language of \cite{lauritzen1996}, $k\to i\leftarrow j$ is a d-connecting path given $\an_{\GG}(j)\setminus\{k\}$ as $i\in \an_{\GG}(j)\setminus\{k\}$.
Thus $\mathbb{P}$ does not entail these conditional independence statements and hence $\mathbb{P}$ is not Markov to $\GG$. 

If $\pa_{\GG}(i)\setminus\pa_{\GG}(j)= \emptyset$ we must have $\pa_{\GG}(j)\setminus\pa_{\GG}(i)\neq \emptyset$ as $\GG$ and $\GG_{i\leftarrow j}$ are assumed to not be Markov equivalent. 
Let $k\in \pa_{\GG}(j)\setminus\pa_{\GG}(i)$ and note that we have $i\leftarrow j\leftarrow k$ in $\HH$. 
As we have $i\to j\leftarrow k$ in $\GG$, any distribution Markov to $\GG$ entails $X_k\independent X_i\mid X_{\an_\GG(i)\cup\an_\GG(k)}$. However, as $j\notin \an_\GG(i)\cup\an_\GG(k)$ we can repeat the previous argument with the path $i\leftarrow j\leftarrow k$.
Hence, the result follows. 
\end{proof}

What \Cref{prop: local turn consitency} says is that turning an edge such that we only add in wanted v-structures or remove unwanted v-structures increases $\BIC$. 
Thus to show that Skeletal Greedy CIM is consistent it is enough to show that we can always find an edge $i\to j$ with this property. 

%---PROPOSITION: Skeletal Greedy CIM consistent---- 
\begin{proposition}
\label{prop: skeletal greedy cim consistent}
Skeletal Greedy CIM is consistent under faithfulness given an oracle-based test for conditional independence. 
\end{proposition}

\begin{proof}
By faithfulness there exists a DAG $\HH$ to which our distribution is faithful. 
In \cite{SG91} they show that, given an oracle-based test for conditional independence, the \texttt{skeleton} algorithm will find the skeleton of $\HH$, say $G$.
Assume we are currently at the characteristic imset of a DAG $\GG$. 

Take any topological order on $[p]$ determined by $\HH$.
Let $i$ be the smallest node in this order such that $\HH|_{\cl(i)}\neq \GG|_{\cl(i)}$, notice that $\ch_\HH(i)\subseteq \ch_\GG(i)$, else we have $i'\in\ch_\HH(i')\setminus \ch_\GG(i)$ giving us $\HH|_{\cl(i')}\neq \GG|_{\cl(i')}$, contradicting that $i$ is minimal with respect to the topological sorting. 
Since $\HH$ and $\GG$ share the same skeleton we get $\pa_\HH(i)\supseteq \pa_\GG(i)$. 
Let $j$ be the smallest node in $\ch_\GG(i)\setminus\ch_\HH(i)$. 
By our choice of $j$ we have that $\GG_{i \leftarrow j}$ is a DAG. 
If $\GG$ and $\GG_{i \leftarrow j}$ are Markov equivalent we could instead have chosen $\GG_{i \leftarrow j}$ as our representative of the MEC, and thus this is a non-case. 
Any new v-structures must be of the form $k\to i\leftarrow$ for $k\in\pa_\GG(i)$, but as $\pa_\HH(i)\supseteq \pa_\GG(i)$ this v-structure is present in $\HH$ as well, and any destroyed v-structures must have been of the form $i\to j\leftarrow k$.
However, since $i$ was smaller than $j$ in the topological sorting none of these v-structures can be present in $\HH$. 
By \Cref{prop: local turn consitency} this increases the score.
Hence we can always turn an edge to increase score, unless we are in the MEC of $\HH$.
Thus, Skeletal Greedy CIM will always find the optimum.
\end{proof}
A first observation about the above proposition is that the proof directly extends to any algorithm that utilizes a fixed set of moves extending the turning phase of GES over the space of DAGs with a fixed skeleton. 
This is because the above proof shows that any such set of moves will always contain at least one move that improves the score, unless we are in the Markov equivalence class of optimal DAGs.
Indeed, we similarly get consistency of EFT when our underlying data-generating distribution is faithful to a tree.
\begin{proposition}
% Let $\mathbb{P}$ be faithful to $\GG$ and the skeleton of $\GG$ is a tree.
The EFT algorithm is consistent under faithfulness when the distribution is faithful to a DAG with skeleton $G$ where $G$ is a tree.
\end{proposition}
\begin{proof}
As was shown in \cite{RP87} the MWSP of $K_p$ where the weights are the negation of the mutual information will, in the limit of infinite data, correctly recover $G$. 
Following the proof of \Cref{prop: skeletal greedy cim consistent} we can always find an edge $i\to j$ so that we can apply \Cref{prop: local turn consitency}.
That is, if we are currently at a DAG $\GG$ that is not the optimum, there is a Markov equivalent graph $\GG'$ with an edge $i\to j$ such that reversing this edge increases the score.
Hence it is enough to show that $\{\GG, \GG_{i\leftarrow j}\}$ is an essential flip whenever $\GG$ and $\GG_{i\leftarrow j}$ are not Markov equivalent. 
However this follows directly from \Cref{thm: subtree condition}. 
\end{proof}
In \cite{WSYU17} it was shown that a version of GES, called \emph{GIES}, based on a mixture of observational and interventional data is not consistent under the faithfulness assumption. 
The above result suggests that this is specifically a feature arising when data from multiple experiments is being used as opposed to a feature of the turning phase of GES.
With these results we further the belief that the hardship of DAG discovery, in the case of purely observational data, is finding the skeleton of the graph. 

%%%%%%%%%%%%%%%%%%%%%%%%%%%%%%%%%%%%%%%%%%%%%%%%%%
%---SUBSUBSECTION: Decomposition for trees ------
% \subsubsection{Decomposition for trees}
We have now looked at how a decomposable, score equivalent and consistent scoring criterion behaves, with respect to $\CIM_G$ for arbitrary $G$. 
It then makes sense to also consider, as we did before, $\CIM_G$ when $G$ is a tree.
Let $S(-,\mathbb{D})$ be a decomposable score equivalent consistent scoring criterion, and let $G$ be a tree.
% As noted in \cite[Lemma 1]{SHL10}, $S(-,\mathbb{D})$ can be seen as a linear function over $\CIM_G$. 
Recall that we defined $N_i=\left\{S\subseteq [p]\colon i\in S\subseteq \ne_G(i)\cup\{i\}, |S|\geq 3\right\}$. 
In \cite[Lemma 1]{SHL10} the authors show that for every score equivalent decomposeable scoring criterion $S(-,\mathbb{D})$ there is a vector $S_\mathbb{D}$ such that $S(\GG,\mathbb{D})= C+S_\mathbb{D}^Tc_\GG$ for all DAGs $\GG$ and some constant $C$. 
% Then it follows by decomposability that $S(-,\mathbb{D})=
Then we have $S_\mathbb{D}=S_G+\sum_{i\in[p]}S_\mathbb{D}|_{N_i}$ for some (affine) linear function $S_G$ that is constant over $\CIM_G$. 
Note that for this decomposition we use that $G$ is a tree and thus the sets $N_i$ are mutually disjoint. 
Then maximizing $S_\mathbb{D}$ over $\CIM_G$, when $G$ is a tree, is the same as maximizing each $S_\mathbb{D}|_{N_i}$ independently. 

\begin{proposition}
\label{prop: local maximiality of tree}
Let $P$ be a (positive) distribution faithful to a DAG $\GG$ with skeleton $G$ a tree, and let $\mathbb{D}$ be a random sample from $P$.
Let $S(-,\mathbb{D})$ be a score equivalent, consistent and decomposable scoring criterion that maximizes at $c_\GG$. 
Then $c_\GG$ simultaneously maximizes $S_\mathbb{D}|_{N_i}$ for all $i\in[p]$; that is, $S_\mathbb{D}|_{N_i}^Tc_\GG=\max_{v\in \CIM_G}\{S_\mathbb{D}|_{N_i}^Tv\}$. 
\end{proposition}

\begin{proof}
For simplicity, for all $i\in [p]$, let $s_i=S_\mathbb{D}|_{N_i}$.
Let $\left\{\HH_i\right\}$ be DAGs such that $s_i^Tc_{\HH_i}=\max_{v\in \CIM_G}\{s_i^Tv\}$ for each $i\in [p]$. 
It is enough to show that $c_{\HH_i}|_{N_i}=c_{\GG}|_{N_i}$ for any given $i$. 
% \textcolor{red}{
As $s_i$ is non-zero only for entries in $N_i$ we can assume that $\HH_i$ has all edges outside of $\HH_i|_{\cl(i)}$ directed away from $j$, since the values of $S_\mathbb{D}(S)$ for $S\centernot\subseteq\cl_G(i)$ do not affect   
% it does not change the score of 
$s_i^Tc_{\HH_i}$. 
% }
Hence $\HH_i$ has no v-structures other than (possibly) the ones centered at $i$.
Assume we have a v-structure at $j\to i\leftarrow k$ in $\HH_i$ not present in $\GG$. 
Then we must have either $j\leftarrow i$ or $i\to k$ in $\GG$, by symmetry assume $j\leftarrow i$. 
Then changing the direction of $j\to i$ in $\HH_i$ can only remove v-structures not present in $\GG$, so by \Cref{prop: local turn consitency} this increases the score of $S_\mathbb{D}^Tc_{\HH_i}$. 
By our decomposition above $S_\mathbb{D}^Tc_{\HH_i}= S_G^Tc_{\HH_i}+\sum_{v\in[p]}s_v^Tc_{\HH_i}$. 
Then, as no other v-structures were created, we only change the summand $s_i^Tc_{\HH_i}$, contradicting the definition of $\HH_i$. 
Hence $\HH_i$ has no v-structures at $i$ not present in $\GG$. 

Assume we have v-structure at $\GG$ not present in $\HH_i$. 
Then we can make the same argument, utilizing the fact that all v-structures at $i$ in $\HH_i$ are present in $\GG$, as before this time adding in v-structures.
Hence $\HH_i$ and $\GG$ must have the exact same v-structures at $i$, and the result follows.
\end{proof}

The above claim asserts that we could (simultaneously) find the structure of $\cl(i)$ for all $i$ and then glue each part together. 
A priori this point of attack could lead to issues, for example if the maximum in $\cl(i)$ would have $i\to j$, but when maximizing in $\cl(j)$ we would have $i\leftarrow j$. 
The above proposition, then states that this will not happen in the limit of large data drawn from a distribution faithful to a DAG.

%%%%%%%%%%%%%%%%%%%%%%%%%%%%%%%%%%%%%%%%%%%%%%%%%%%%%
%---SUBSECTION: Experimental Results ------
\subsection{Experimental Results}
\label{subsec: computations}
We now analyze the empirical performance of the EFT algorithm (Algorithm~\ref{alg: EFT}) for learning polytrees on simulated and real data. 
The EFT algorithm can be viewed as operating in two phases: in the first phase the skeleton of the polytree is estimated.  
In the second phase, a BIC-optimal orientation of the estimated skeleton is identified.  
For the first phase, since classic CI-test approaches cannot be guaranteed to return a skeleton that is a tree, we use the approach of Rebane and Pearl \cite{RP87}, where we assign a weight to each edge $i - j$ of a complete graph: the negative mutual information $-I(X_i; X_j)$.  
We then identify a minimum weight spanning forest (MWSF) of this weighted complete graph.
(In the implementation available at \cite{github}, this step is done via Kruskal's algorithm in the \texttt{networkX} package in \texttt{Python}. 
Other algorithms for learning a MWSF are available as well via this package and can be called in the EFT implementation at \cite{github}.)
In the second phase, edges of the estimated skeleton are oriented to recover an essential graph. 
Instead of the CI-tests used by the RP-algorithm (which can be prone to error propagation), EFT uses the characterization of the edges of $\CIM_T$ in \Cref{thm: edges of trees characterization} to search for a BIC-optimal polytree on the estimated skeleton. 

%---SUBSUBSECTION: Simulations---
\subsubsection{Simulations}
\label{subsubsec: simulations}
\begin{figure}
\centering
\includegraphics[width=.7\textwidth]{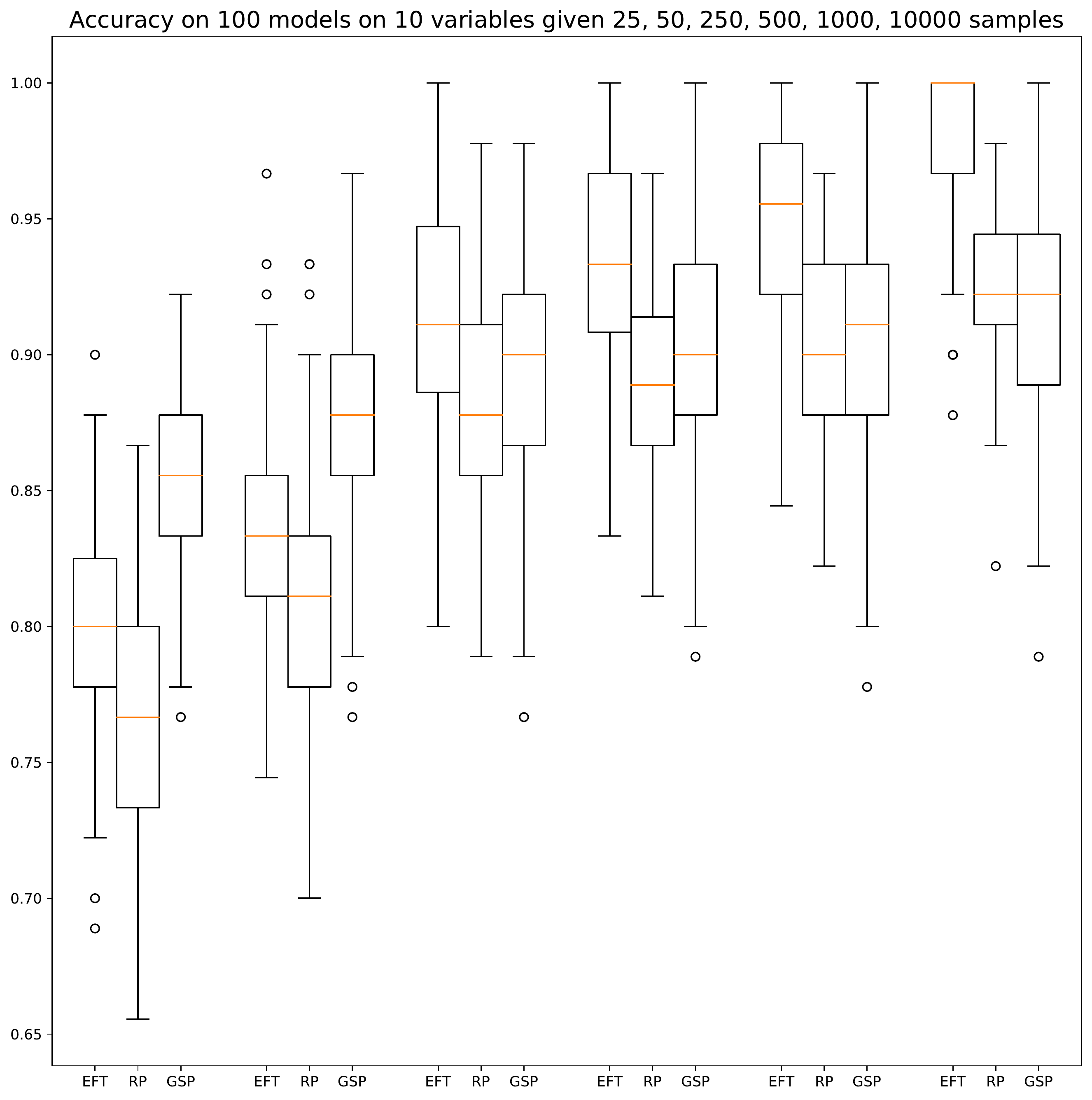}
\caption{Box plots of the accuracy of the estimated essential graphs of $100$ different random linear Gaussian polytree models by EFT, RP and GreedySP (GSP) for varying sample sizes: 25, 50, 250, 500, 1000, and 10000 going from left-to-right.}
\label{fig: EFT_RP}
\end{figure}

We compared these two methods by analyzing their performance on randomly generated linear structural equation models with Gaussian noise whose underlying DAG is a polytree on $10$ nodes.  
To generate these models, we uniformly at random generated a \emph{Pr\"ufer code} (i.e., a sequence of length $8$ with entries in $[10]$ that uniquely corresponds to an undirected tree on nodes $[10]$) to produce the skeleton $T$ of the polytree.
The edges of the skeleton were then oriented independently and uniformly at random to yield a polytree $\mathcal{T}$, and each edge $i\rightarrow j$ of $\mathcal{T}$ was assigned a weight $\lambda_{ij}$ drawn from the uniform distribution over $[-1,0)\cup(0,1]$. 
We then sampled from a multivariate Gaussian distribution over the random variables $(X_1,\ldots, X_p)$ where 
\[
X_i := \sum_{k\in\pa_{\mathcal{T}}(i)}\lambda_{ki}X_k+\varepsilon_i,
\]
where $\varepsilon_1,\ldots,\varepsilon_p$ are mutually independent standard normal random variables. 
For each $n\in\{15,20,250,500,1000,10000\}$, we generated $100$ such models and drew a random sample of size $n$.
The RP-algorithm, EFT, and GreedySP were then tasked with recovering the data-generating polytree $\mathcal{T}$ based on each sample.  
The constraint-based tests of the RP-algorithm and GreedySP were performed with a cut-off threshold of $\alpha = 0.05$, and the depth and run parameters of GreedySP were chosen to be $d = 4$, $r = 5$; the default values in the \texttt{causaldag} python package \cite{squires2018causaldag} implementation of GreedySP.
Although not specifically designed to learn polytrees, GreedySP was included to give a benchmark of the performance of the polytree-specific hybrid algorithms, EFT and RP, against a general hybrid causal discovery algorithm. 
The accuracy of each estimated DAG (computed as the fraction of off-diagonal entries in the adjacency matrix of the estimated essential graph that agree with the adjacency matrix of the essential graph of $\mathcal{T}$) was recorded.  
The results are presented in Figure~\ref{fig: EFT_RP}. 
We see that EFT outperforms both algorithms over all sample sizes greater than 50, and does increasingly better for larger sample sizes (reflecting the asymptotic consistency results observed in subsection~\ref{subsec: consistency}).  
At $10000$ samples EFT perfectly recovers at least $50\%$ of the models and $25\%$ already for $250$ samples.
Since the RP algorithm and EFT differ only in their second phases, these results suggest that meaningful gains can be had by replacing CI-tests in the second phase of RP with a greedy search over the edges of $\CIM_T$ when the data-generating distribution is approximately normal. 

%---FIGURE: EFT versus RP and GreedySP----

\begin{figure}
\centering
\includegraphics[width=.7\textwidth]{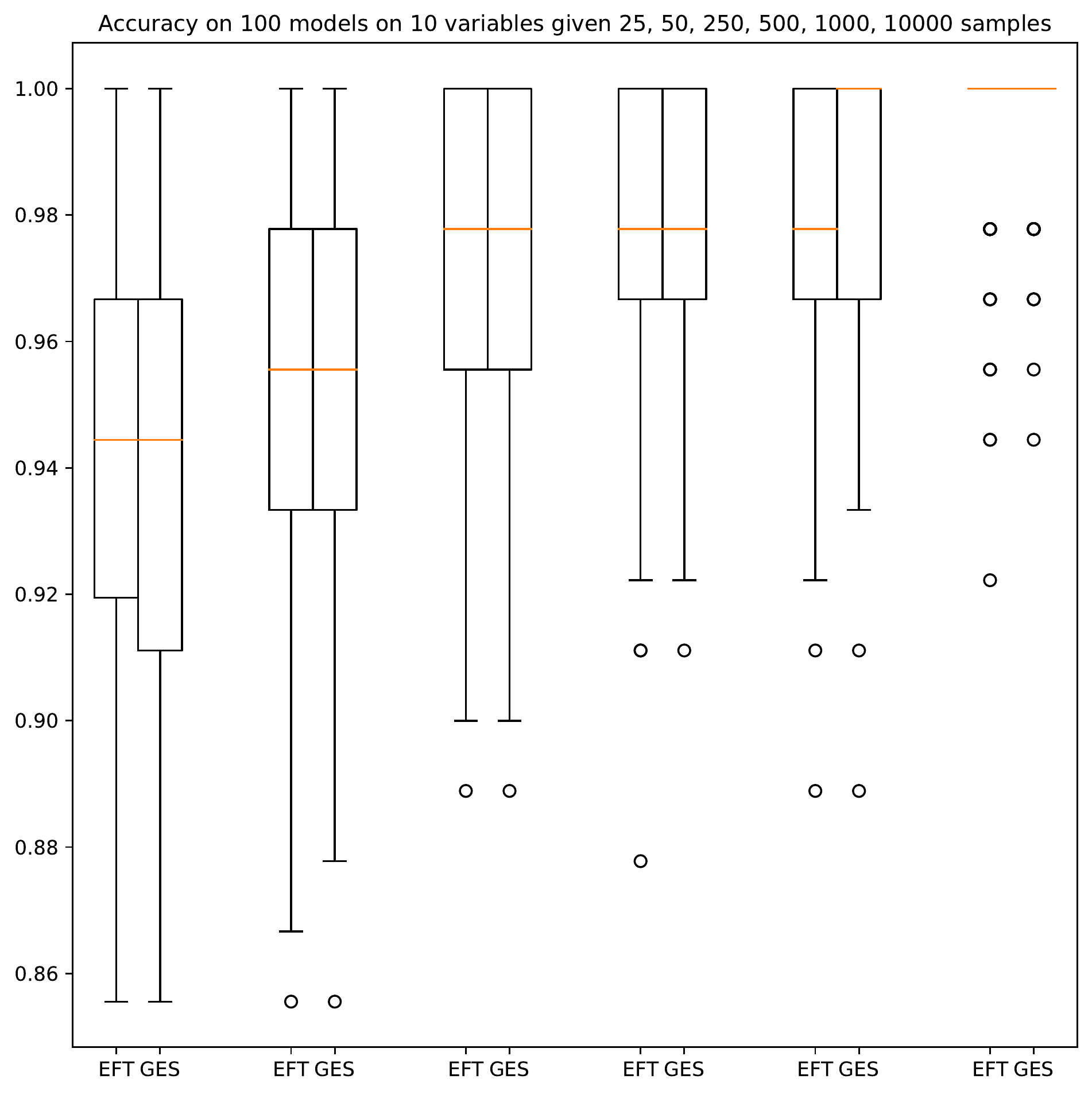}
\caption{Box plots of the accuracy of the estimated essential graphs of $100$ different random linear Gaussian polytree models by EFT and the turning phase of GES with the true skeleton given as background knowledge. Sample sizes vary as 25, 50, 250, 500, 1000, and 10000 going from left-to-right.}
\label{fig: GES_EFT}
\end{figure}

As shown in \Cref{prop: turn pair trees}, the moves used by EFT in its second phase generalize the moves of the turning phase of GES, which uses only single edge reversals.  
Hence, it is of interest to see if the additional moves provided by using all edges of $\CIM_T$ yield substantial gains over the turning phase of GES.  
To test this, we generated $100$ random linear Gaussian polytree models (as described above) and drew random samples of size $n = 15,20,250,500,1000,10000$ from each model.  
We then had EFT and the turning phase of GES estimate the data-generating polytree for each model with the true skeleton given as background knowledge.  
The accuracies of the estimated essential graphs is presented in the box plots in Figure~\ref{fig: GES_EFT}.  

%---FIGURE: EFT versus GES----

Intuitively, we expect that the additional moves provided to EFT over GES by using all edges of $\CIM_T$ when searching (as opposed to only turn pairs) does not substantially boost performance.  
This is because GES is known to perform quite well when the data-generating system is sparse.
The results presented in Figure~\ref{fig: GES_EFT} confirm this intuition. 
The accuracy rates are highly comparable, especially as the sample size grows.
This suggests that one could replace phase two of EFT with the turning phase of GES to obtain a comparable, but perhaps slightly more efficient hybrid algorithm for estimating polytrees. 
The results also suggest interesting features of the geometry of $\CIM_G$.  
For instance, the subgraph consisting only of (polytope) edges corresponding to (essential graph) edge reverals of the edge graph of $\CIM_G$ appears to provide as good of a search mechanism as the entire graph.  
It would be interesting to have a deeper understanding of how the essential graph edge reversals are distributed along the edge graph of $\CIM_G$ for $G$ a tree.  
It would also be of interest to understand if this phenomenon occurs for denser choices of $G$, in which one would generalize EFT by replacing the search phase based on essential flips with a corresponding characterization of edges of $\CIM_G$ for the given skeleton.
Since GES tends to perform increasingly worse as the density of the graph increases, one might observe that GES's competitiveness with a such generalized version EFT decreases as the graphs become increasingly dense.
Since the set of moves used by the search phase of such a generalization of EFT would necessarily include edge reversals, such an algorithm would be expected to be consistent at least as often as GES.  
It is possible that, for certain skeletons, we see improvements in accuracy over GES in the finite sample setting.  
As well, such moves could possibly improve efficiency if they offer substantially faster paths between two graphs than any corresponding sequence of edge reversals.  
Such complexity questions could be better understood via an analysis of the diameter of the edge graph of $\CIM_G$ versus the diameter of its subgraph with edges corresponding to edge reversals. 

%---SUBSUBSECTION: Real Data---
\subsubsection{Real Data Analysis}
\label{subsubsec: real data}
The protein signaling dataset of Sachs et al. \cite{sachs05} is a standard benchmark dataset in causal inference consisting of $7466$ abundance measurements of phosolipids and phosphoproteins taken under varying experimental conditions in primary human immune system cells. 
The different experimental conditions are produced by reagents that inhibit or activate different sets of signaling proteins within the system before measurements are taken.  
Such measurements are termed ``samples from an \emph{interventional distribution}'' \cite{pearl2009causality} since the reagents are altering the data-generating distribution. 
To learn a (probabilistic) DAG model for the system we work with only \emph{observational data}; i.e., samples drawn from the unaltered data-generating distribution.  
An observational dataset can be extracted from the Sachs data as described in \cite{WSYU17}, which results in a sample size of $1755$.

%---FIGURE: Real Data---
\begin{figure}
\centering
\begin{subfigure}{.5\textwidth}
\includegraphics[width =\textwidth]{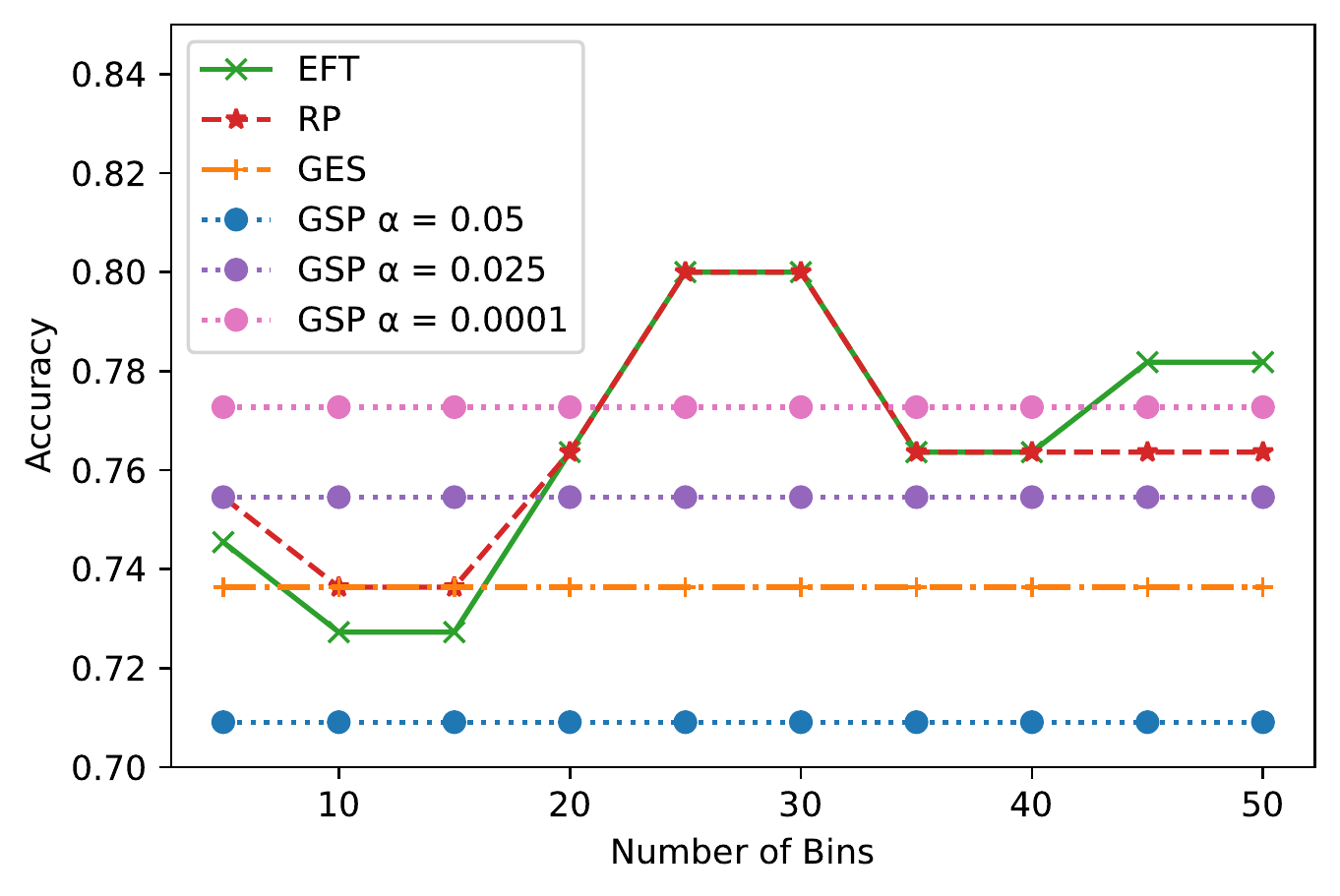}
\caption{Number of bins versus accuracy rate of estimated essential graphs.  Note that GES and GreedySP (GSP) have constant accuracies in these plots as binning is not required for these algorithms.}
\label{subfig: accuracy}
\end{subfigure}
\hfill
% \hspace{-30pt}
\begin{subfigure}{.45\textwidth}
\includegraphics[width =\textwidth]{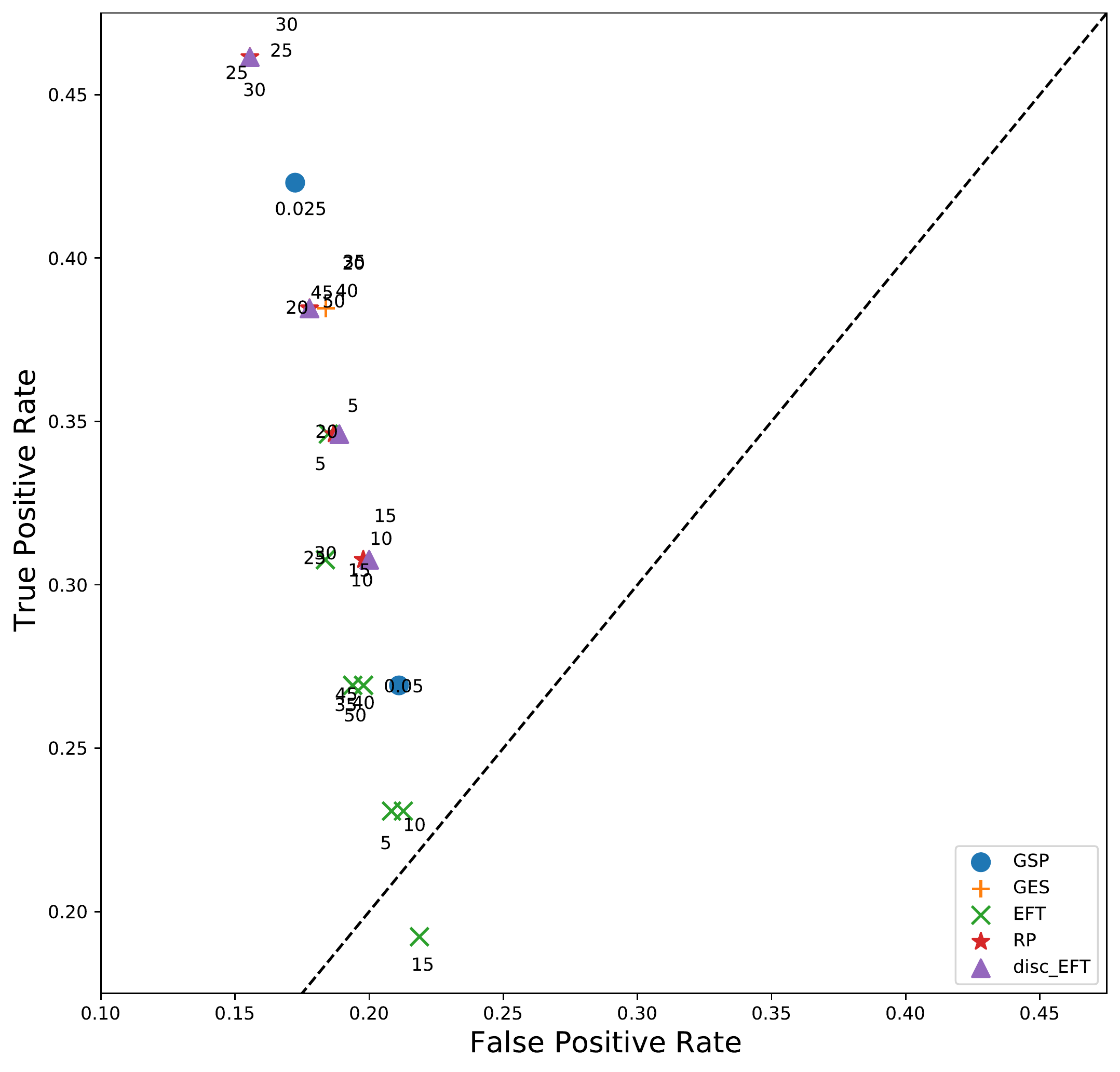}
\caption{ROC plot for the estimated essential graphs.  The points for EFT and the RP-algorithm are labeled by their bin count. GSP points are labeled by $\alpha$-value.  The dashed line indicates random guessing.}
\label{subfig: roc}
\end{subfigure}
\caption{Performance of EFT, RP, GES, and GSP on Sachs et al. protein signaling dataset.}
\label{fig: real data}
\end{figure}

\begin{figure}
\centering
\begin{subfigure}{.3\textwidth}
\begin{tikzpicture}[scale=0.15]

\def\n{11}
\def \radius {12cm}
\def \margin {3} % margin in angles, depends on the radius

  \node[circle, draw, inner sep=1pt, minimum width=1pt] (PIP2) at ({360/\n * (3)}:\radius) {\tiny PIP2};
  \node[circle, draw, inner sep=1pt, minimum width=1pt] (Plcg) at ({360/\n * (2)}:\radius) {\tiny Plc$\gamma$};
  \node[circle, draw, inner sep=1pt, minimum width=1pt] (Mek) at ({360/\n * (1)}:\radius) {\tiny Mek};
  \node[circle, draw, inner sep=1pt, minimum width=1pt] (Raf) at ({360/\n * (11)}:\radius) {\tiny Raf};
  \node[circle, draw, inner sep=1pt, minimum width=1pt] (Jnk) at ({360/\n * (10)}:\radius) {\tiny Jnk};
  \node[circle, draw, inner sep=1pt, minimum width=1pt] (p38) at ({360/\n * (9)}:\radius) {\tiny p38};
  \node[circle, draw, inner sep=1pt, minimum width=1pt] (PKC) at ({360/\n * (8)}:\radius) {\tiny PKC};
  \node[circle, draw, inner sep=1pt, minimum width=1pt] (PKA) at ({360/\n * (7)}:\radius) {\tiny PKA};
  \node[circle, draw, inner sep=1pt, minimum width=1pt] (Akt) at ({360/\n * (6)}:\radius) {\tiny Akt};
  \node[circle, draw, inner sep=1pt, minimum width=1pt] (Erk) at ({360/\n * (5)}:\radius) {\tiny Erk};
  \node[circle, draw, inner sep=1pt, minimum width=1pt] (PIP3) at ({360/\n * (4)}:\radius) {\tiny PIP3};
  
  \draw[->, blue] (Raf) -- (Mek) ;
  \draw[<-, blue] (Raf) -- (Mek) ;
  \draw[->, blue] (PKC) -- (p38) ;
  \draw[->, blue] (PKC) -- (Jnk) ;
  \draw[->, blue] (PKC) -- (Raf) ;
  \draw[->, blue] (PKC) -- (Mek) ;
  \draw[->, blue] (PKA) -- (p38) ;
  \draw[->, blue] (PKA) -- (Jnk) ;
  \draw[->, blue] (PKA) -- (Raf) ;
  \draw[->, blue] (PKA) -- (Mek) ;
  \draw[->, blue] (PKA) -- (Erk) ;
  \draw[->, blue] (PKA) -- (Akt) ;
  \draw[->, blue] (PIP3) -- (Akt) ;
  \draw[->, blue] (PIP3) -- (Plcg) ;
  \draw[<-, blue] (PIP3) -- (Plcg) ;
  \draw[->, blue] (PIP3) -- (PIP2) ;
  \draw[<-, blue] (PIP3) -- (PIP2) ;
  \draw[->, blue] (Plcg) -- (PIP2) ;
  \draw[<-, blue] (Plcg) -- (PIP2) ;
  \draw[->, blue] (Plcg) -- (PKC) ;
  \draw[<-, blue] (Plcg) -- (PKC) ;
  \draw[->, blue] (Mek) -- (Erk) ;
\end{tikzpicture}
\caption{Ground truth.}
\label{subfig: ground-truth}
\end{subfigure}
\hfill
% \hspace{-30pt}
\begin{subfigure}{.3\textwidth}
\begin{tikzpicture}[scale=0.15]

\def \n {11}
\def \radius {12cm}
\def \margin {3} % margin in angles, depends on the radius

  \node[circle, draw, inner sep=1pt, minimum width=1pt] (PIP2) at ({360/\n * (3)}:\radius) {\tiny PIP2};
  \node[circle, draw, inner sep=1pt, minimum width=1pt] (Plcg) at ({360/\n * (2)}:\radius) {\tiny Plc$\gamma$};
  \node[circle, draw, inner sep=1pt, minimum width=1pt] (Mek) at ({360/\n * (1)}:\radius) {\tiny Mek};
  \node[circle, draw, inner sep=1pt, minimum width=1pt] (Raf) at ({360/\n * (11)}:\radius) {\tiny Raf};
  \node[circle, draw, inner sep=1pt, minimum width=1pt] (Jnk) at ({360/\n * (10)}:\radius) {\tiny Jnk};
  \node[circle, draw, inner sep=1pt, minimum width=1pt] (p38) at ({360/\n * (9)}:\radius) {\tiny p38};
  \node[circle, draw, inner sep=1pt, minimum width=1pt] (PKC) at ({360/\n * (8)}:\radius) {\tiny PKC};
  \node[circle, draw, inner sep=1pt, minimum width=1pt] (PKA) at ({360/\n * (7)}:\radius) {\tiny PKA};
  \node[circle, draw, inner sep=1pt, minimum width=1pt] (Akt) at ({360/\n * (6)}:\radius) {\tiny Akt};
  \node[circle, draw, inner sep=1pt, minimum width=1pt] (Erk) at ({360/\n * (5)}:\radius) {\tiny Erk};
  \node[circle, draw, inner sep=1pt, minimum width=1pt] (PIP3) at ({360/\n * (4)}:\radius) {\tiny PIP3};
  
  \draw[<-, red] (PKC) -- (PIP3) ;
  \draw[->, red] (Akt) -- (Erk) ;
  \draw[<-, red] (Akt) -- (Erk) ;
  \draw[<-, red] (PKA) -- (Akt) ;
  \draw[->, red] (PKC) -- (PIP3) ;
  \draw[<-, red] (PKC) -- (p38) ;
  \draw[->, red] (Jnk) -- (Mek) ;
  \draw[->, red] (Jnk) -- (PKA) ;
  \draw[->, red] (Jnk) -- (PKC) ;
  \draw[->, green] (Raf) -- (Mek) ;
  \draw[->, green] (Plcg) -- (PIP2) ;
  \draw[<-, green] (Plcg) -- (PIP2) ;
  \draw[->, green] (PIP2) -- (PIP3) ;
  \draw[<-, green] (PIP2) -- (PIP3) ;
  \draw[->, green] (PKA) -- (Akt) ;
  \draw[->, green] (PKA) -- (Jnk) ;
  \draw[->, green] (PKC) -- (Jnk) ;
  \draw[->, green] (PKC) -- (p38) ;

\end{tikzpicture}
\label{subfig: EFT}
\caption{EFT\\ (Gaussian BIC, bins $= 20$)}
\end{subfigure}
\hfill
% \hspace{-30pt}
\begin{subfigure}{.3\textwidth}
\begin{tikzpicture}[scale=0.15]

\def \n {11}
\def \radius {12cm}
\def \margin {3} % margin in angles, depends on the radius

  \node[circle, draw, inner sep=1pt, minimum width=1pt] (PIP2) at ({360/\n * (3)}:\radius) {\tiny PIP2};
  \node[circle, draw, inner sep=1pt, minimum width=1pt] (Plcg) at ({360/\n * (2)}:\radius) {\tiny Plc$\gamma$};
  \node[circle, draw, inner sep=1pt, minimum width=1pt] (Mek) at ({360/\n * (1)}:\radius) {\tiny Mek};
  \node[circle, draw, inner sep=1pt, minimum width=1pt] (Raf) at ({360/\n * (11)}:\radius) {\tiny Raf};
  \node[circle, draw, inner sep=1pt, minimum width=1pt] (Jnk) at ({360/\n * (10)}:\radius) {\tiny Jnk};
  \node[circle, draw, inner sep=1pt, minimum width=1pt] (p38) at ({360/\n * (9)}:\radius) {\tiny p38};
  \node[circle, draw, inner sep=1pt, minimum width=1pt] (PKC) at ({360/\n * (8)}:\radius) {\tiny PKC};
  \node[circle, draw, inner sep=1pt, minimum width=1pt] (PKA) at ({360/\n * (7)}:\radius) {\tiny PKA};
  \node[circle, draw, inner sep=1pt, minimum width=1pt] (Akt) at ({360/\n * (6)}:\radius) {\tiny Akt};
  \node[circle, draw, inner sep=1pt, minimum width=1pt] (Erk) at ({360/\n * (5)}:\radius) {\tiny Erk};
  \node[circle, draw, inner sep=1pt, minimum width=1pt] (PIP3) at ({360/\n * (4)}:\radius) {\tiny PIP3};
  
  \draw[->, red] (Mek) -- (Jnk) ;
  \draw[->, red] (Erk) -- (Akt) ;
  \draw[<-, red] (Erk) -- (Akt) ;
  \draw[->, red] (Akt) -- (PKA) ;
  \draw[<-, red] (PKC) -- (p38) ;
  \draw[->, red] (Jnk) -- (Mek) ;
  \draw[->, red] (Jnk) -- (PKA) ;
  \draw[->, red] (Jnk) -- (PKC) ;
  \draw[->, green] (Raf) -- (Mek) ;
  \draw[<-, green] (Raf) -- (Mek) ;
  \draw[->, green] (PKA) -- (Akt) ;
  \draw[->, green] (PKA) -- (Jnk) ;
  \draw[<-, green] (Plcg) -- (PKC) ;
  \draw[->, green] (PKC) -- (p38) ;
  \draw[->, green] (PKC) -- (Jnk) ;
  \draw[->, green] (Plcg) -- (PIP2) ;
  \draw[<-, green] (Plcg) -- (PIP2) ;
  \draw[->, green] (Plcg) -- (PKC) ;
  \draw[->, green] (PIP3) -- (PIP2) ;
  \draw[<-, green] (PIP3) -- (PIP2) ;
  
\end{tikzpicture}
\caption{RP and discEFT\\ (bins $= 30$)}
\label{subfig: RP}
\end{subfigure}
\linebreak
% \hspace{-30pt}
\begin{subfigure}{.3\textwidth}
\begin{tikzpicture}[scale=0.15]

\def \n {11}
\def \radius {12cm}
\def \margin {3} % margin in angles, depends on the radius

  \node[circle, draw, inner sep=1pt, minimum width=1pt] (PIP2) at ({360/\n * (3)}:\radius) {\tiny PIP2};
  \node[circle, draw, inner sep=1pt, minimum width=1pt] (Plcg) at ({360/\n * (2)}:\radius) {\tiny Plc$\gamma$};
  \node[circle, draw, inner sep=1pt, minimum width=1pt] (Mek) at ({360/\n * (1)}:\radius) {\tiny Mek};
  \node[circle, draw, inner sep=1pt, minimum width=1pt] (Raf) at ({360/\n * (11)}:\radius) {\tiny Raf};
  \node[circle, draw, inner sep=1pt, minimum width=1pt] (Jnk) at ({360/\n * (10)}:\radius) {\tiny Jnk};
  \node[circle, draw, inner sep=1pt, minimum width=1pt] (p38) at ({360/\n * (9)}:\radius) {\tiny p38};
  \node[circle, draw, inner sep=1pt, minimum width=1pt] (PKC) at ({360/\n * (8)}:\radius) {\tiny PKC};
  \node[circle, draw, inner sep=1pt, minimum width=1pt] (PKA) at ({360/\n * (7)}:\radius) {\tiny PKA};
  \node[circle, draw, inner sep=1pt, minimum width=1pt] (Akt) at ({360/\n * (6)}:\radius) {\tiny Akt};
  \node[circle, draw, inner sep=1pt, minimum width=1pt] (Erk) at ({360/\n * (5)}:\radius) {\tiny Erk};
  \node[circle, draw, inner sep=1pt, minimum width=1pt] (PIP3) at ({360/\n * (4)}:\radius) {\tiny PIP3};
  
  \draw[->, red] (Plcg) -- (Mek) ;
  \draw[->, red] (PIP2) -- (Mek) ;
  \draw[->, red] (PIP2) -- (PKA) ;
  \draw[->, red] (PIP3) -- (Mek) ;
  \draw[->, red] (PIP3) -- (PKA) ;
  \draw[->, red] (Erk) -- (Akt) ;
  \draw[->, red] (Erk) -- (PKA) ;
  \draw[->, red] (Akt) -- (Erk) ;
  \draw[->, red] (Akt) -- (PKA) ;
  \draw[->, red] (PKA) -- (PIP2) ;
  \draw[->, red] (PKA) -- (PIP3) ;
  \draw[->, red] (p38) -- (PKC) ;
  \draw[->, red] (p38) -- (Jnk) ;
  \draw[<-, red] (p38) -- (Jnk) ;
  \draw[->, green] (Raf) -- (Mek) ;
  \draw[->, green] (Plcg) -- (PIP2) ;
  \draw[->, green] (PIP2) -- (Plcg) ;
  \draw[->, green] (PIP2) -- (PIP3) ;
  \draw[->, green] (PIP3) -- (PIP2) ;
  \draw[->, green] (PKA) -- (Mek) ;
  \draw[->, green] (PKA) -- (Erk) ;
  \draw[->, green] (PKA) -- (Akt) ;
  \draw[->, green] (PKC) -- (p38) ;

\end{tikzpicture}
\caption{GES}
\label{subfig: GES}
\end{subfigure}
% \hfill
\hspace{30pt}
\begin{subfigure}{.3\textwidth}
\begin{tikzpicture}[scale=0.15]

\def \n {11}
\def \radius {12cm}
\def \margin {3} % margin in angles, depends on the radius

  \node[circle, draw, inner sep=1pt, minimum width=1pt] (PIP2) at ({360/\n * (3)}:\radius) {\tiny PIP2};
  \node[circle, draw, inner sep=1pt, minimum width=1pt] (Plcg) at ({360/\n * (2)}:\radius) {\tiny Plc$\gamma$};
  \node[circle, draw, inner sep=1pt, minimum width=1pt] (Mek) at ({360/\n * (1)}:\radius) {\tiny Mek};
  \node[circle, draw, inner sep=1pt, minimum width=1pt] (Raf) at ({360/\n * (11)}:\radius) {\tiny Raf};
  \node[circle, draw, inner sep=1pt, minimum width=1pt] (Jnk) at ({360/\n * (10)}:\radius) {\tiny Jnk};
  \node[circle, draw, inner sep=1pt, minimum width=1pt] (p38) at ({360/\n * (9)}:\radius) {\tiny p38};
  \node[circle, draw, inner sep=1pt, minimum width=1pt] (PKC) at ({360/\n * (8)}:\radius) {\tiny PKC};
  \node[circle, draw, inner sep=1pt, minimum width=1pt] (PKA) at ({360/\n * (7)}:\radius) {\tiny PKA};
  \node[circle, draw, inner sep=1pt, minimum width=1pt] (Akt) at ({360/\n * (6)}:\radius) {\tiny Akt};
  \node[circle, draw, inner sep=1pt, minimum width=1pt] (Erk) at ({360/\n * (5)}:\radius) {\tiny Erk};
  \node[circle, draw, inner sep=1pt, minimum width=1pt] (PIP3) at ({360/\n * (4)}:\radius) {\tiny PIP3};
  
  \draw[->, red] (Raf) -- (PKA) ;
  \draw[->, red] (Mek) -- (PKA) ;
  \draw[->, red] (Plcg) -- (Mek) ;
  \draw[->, red] (Plcg) -- (p38) ;
  \draw[->, red] (PIP2) -- (Mek) ;
  \draw[->, red] (PIP2) -- (PKA) ;
  \draw[->, red] (PIP2) -- (Jnk) ;
  \draw[->, red] (PIP3) -- (Mek) ;
  \draw[->, red] (PIP3) -- (PKA) ;
  \draw[->, red] (Erk) -- (Akt) ;
  \draw[->, red] (Akt) -- (Erk) ;
  \draw[->, red] (p38) -- (PKC) ;
  \draw[->, red] (Jnk) -- (PIP2) ;
  \draw[->, red] (Jnk) -- (p38) ;
  \draw[->, green] (Raf) -- (Mek) ;
  \draw[->, green] (Plcg) -- (PIP2) ;
  \draw[->, green] (Plcg) -- (PIP3) ;
  \draw[->, green] (PIP2) -- (Plcg) ;
  \draw[->, green] (PIP2) -- (PIP3) ;
  \draw[->, green] (PIP3) -- (Plcg) ;
  \draw[->, green] (PIP3) -- (PIP2) ;
  \draw[->, green] (PKA) -- (Erk) ;
  \draw[->, green] (PKA) -- (Akt) ;

\end{tikzpicture}
\caption{GreedySP ($\alpha = 0.025$)}
\label{subfig: GreedySP}
\end{subfigure}
\caption{The ground-truth essential graph for the protein signaling system and the optimal learned essential graphs for each algorithm according to the ROC plot in Figure~\ref{subfig: roc}. Here, undirected edges in essential graphs are represented by bidirected edges.}
\label{fig: networks}
\end{figure}
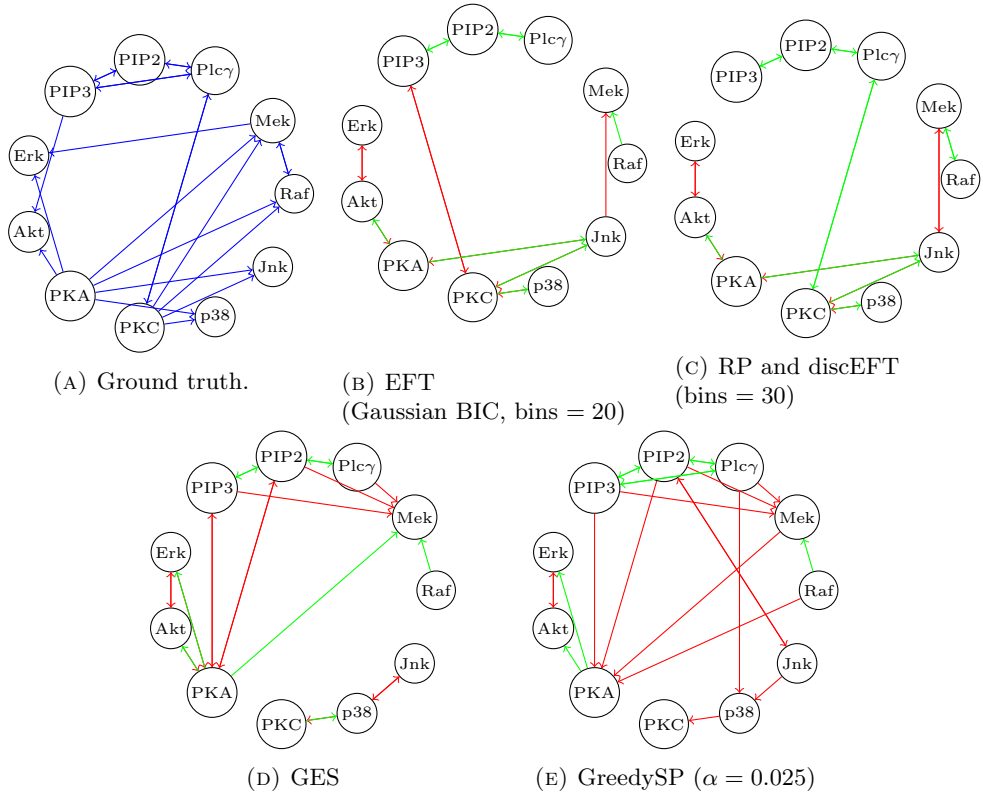

We compared the accuracy of the essential graphs estimated by EFT, the RP algorithm, GES and GreedySP for the protein signaling system based on these $1755$ samples in Figure~\ref{fig: real data}. 
Here, accuracy of the estimated essential graphs was computed relative to the accepted ground-truth network depicted in Figure~\ref{subfig: ground-truth}. 
Since the ground-truth network is not a polytree, only GES and GreedySP could potentially learn the network exactly.  
However, the results presented in Figure~\ref{subfig: accuracy} show that both EFT and the RP-algorithm tend to perform best in regards to accuracy.  
In this analysis, the cut-off threshold for CI-testing in the RP-algorithm was set to $\alpha = 0.05$, and GES and the second phase of EFT were conducted using a Gaussian MLE in the BIC computations.
GreedySP (GSP) was implemented with depth and run parameters $d = 4$ and $r = 5$ and a Gaussian (partial correlation) test for conditional independence.  
Since the data is believed to be highly non-Gaussian \cite{sachs05}, the Gaussianity assumptions on EFT, GES and GreedySP can result in suboptimal performance.  
Hence, we also assessed the performance of GreedySP with varying cut-off thresholds for CI-testing ranging over $\alpha\in\{0.05, 0.025, 0.001, 0.005, 0.0025, 0.001, 0.0005, 0.00025, 0.0001\}$. 
In the first phase of EFT and the RP-algorithm, where the skeleton is estimated, the data (which is drawn from a continuous distribution) is binned to yield a discrete approximation of the sampling distribution that is used to calculate the empirical mutual information weights.  
The results for different numbers of bins were analyzed, ranging over $\{5,10,15,\ldots,50\}$.
We see in Figure~\ref{subfig: accuracy} that the best accuracy of both EFT and the RP-algorithm is achieved simultaneously by $25$ and $30$ bins, with both achieving identical accuracy measurements.
In this plot, we report GreedySP's highest accuracy rate, which was achieved with $\alpha=0.0001$ along with the results for $\alpha = 0.05, 0.025$ (the default choice in \texttt{causaldag} \cite{squires2018causaldag} and the best performing threshold value in Figure~\ref{subfig: roc}, respectively). 
Figure~\ref{subfig: roc} presents an ROC plot, which shows the false positive rates versus the true positive rates, of the different experiments.  
(Here, for the sake of space, we only report GreedySP's best performing value $\alpha = 0.025$ and the default value in \texttt{causaldag} $\alpha = 0.05$.)
Despite tying with EFT in optimal accuracy scores in Figure~\ref{subfig: accuracy}, we see in Figure~\ref{subfig: roc} that the RP-algorithm outperforms EFT, GES and GreedySP, having lower false positive rates and higher true positive rates.  
% From an accuracy perspective, EFT is competitive with the RP algorithm on both simulated Gaussian data (Figure~\ref{fig: EFT_RP}) and real data (Figure~\ref{subfig: accuracy}). 
% On the other hand, the ROC plot in Figure~\ref{subfig: roc} has the RP algorithm outperforming EFT, GES and GreedySP. 
One possible explanation for this is that the data is known to be highly non-Gaussian, contrary to the implemented assumptions for EFT, GES and GreedySP.  
On the other hand, the RP algorithm works directly with a discrete approximation of the data. 
To test this, we implemented a version of EFT, denoted discEFT in Figure~\ref{fig: real data}, that discretizes the data according to the specified number of bins and computes the BIC score of the discretized data set.
Due to computational complexity limitations, this was done only for numbers of bins ranging over $\{5,10,15,20,25,30\}$.
We see that EFT assuming discrete data (discEFT) then performs equally as well as the RP-algorithm in Figure~\ref{subfig: roc} for the optimal numbers of bins (25 and 30) identified in Figure~\ref{subfig: accuracy}. 
We note that the accuracy of EFT on discretized data for these numbers of bins is $0.8$, tying with that of Gaussian EFT and the RP-algorithm.
All details of this real data analysis and the simulations are available for reproduction at \cite{github} (\url{https://github.com/soluslab/causalCIM}).
Overall, these analyses suggest that EFT performs equally well as the RP-algorithm on the real data example and outperforms it in the linear Gaussian model regime. 

%---FIGURE: Sachs Networks---

%---SUBSECTION: Acknowledgements---
\subsection*{Acknowledgements}
All three authors were partially supported by the Wallenberg AI, Autonomous Systems and Software Program (WASP) funded by the Knut and Alice Wallenberg Foundation. Liam Solus was additionally supported by the G\"oran Gustafsson Prize for Young Researchers, and Starting Grant No.\ 2019-05195 from The Swedish Research Council (Vetenskapsr\aa{}det).

%%%%%%%%%%%%%%%%%%%%%%%%%%
%---BIBLIOGRAPHY
\bibliographystyle{plain}
\bibliography{references}

\end{document}